\documentclass[preprint,superscriptaddress]{revtex4-1}
\usepackage{amsfonts}
\usepackage{graphicx}
\usepackage{epstopdf}
\usepackage{epsfig}
\usepackage{amsmath}
\usepackage[normalem]{ulem}
\usepackage{color}
\usepackage{makecell}
\usepackage{array}
\usepackage{booktabs}
\usepackage{mathtools}
\usepackage{bm}
\usepackage{color}
\usepackage{array}
\usepackage{booktabs}
\usepackage{tabularx}

\usepackage{xr}

\definecolor{blue}{rgb}{0,0,1}

\newcolumntype{Y}{>{\centering\arraybackslash}X}

\newcommand{\Rmnum}[1]{\expandafter\@slowromancap\romannumeral #1@}
\makeatother
\begin{document}

\title{Finding nonlinear system equations and complex network structures from data: a sparse optimization approach}

\date{\today}

\author{Ying-Cheng Lai} 
\affiliation{School of Electrical, Computer and Energy Engineering, Arizona State University, Tempe, AZ 85287, USA}
\affiliation{Department of Physics, Arizona State University,
Tempe, Arizona 85287, USA}

\begin{abstract}

In applications of nonlinear and complex dynamical systems, a common situation
is that the system can be measured but its structure and the detailed rules of 
dynamical evolution are unknown. The inverse problem is to determine the 
system equations and structure based solely on measured time series. Recently,
methods based on sparse optimization have been developed. For example, the 
principle of exploiting sparse optimization such as compressive sensing to 
find the equations of nonlinear dynamical systems from data was articulated 
in 2011~\footnote{The idea of using sparse optimization to discover system equations was also published in 2016 [S. L. Brunton, J. L. Proctor, and J. Nathan Kutz, ``Discovering governing equations from data by sparse identification of nonlinear dynamical systems,'' Proc. Nat. Acad. Sci. {\bf 113}, 3932-3937 (2016)] by Prof. Kutz's group at University of Washington. On October 8, 2015, Prof. Kutz gave a talk at an AFOSR (Air Force Office of Scientific Research) Program Review meeting about this idea. The author of the present manuscript (Y.-C. Lai) was in the audience and rose to point out politely that the idea had already been published in 2011 by the ASU group, and he immediately forwarded the 2011 paper to Prof.~Kutz.}.
This article presents a brief review of the recent progress in this area. The 
basic idea is to expand the equations governing the dynamical evolution of the 
system into a power series or a Fourier series of a finite number of terms and 
then to determine the vector of the expansion coefficients based solely on data 
through sparse optimization. Examples discussed here include discovering the 
equations of stationary or nonstationary chaotic systems to enable prediction 
of dynamical events such as critical transition and system collapse, inferring 
the full topology of complex networks of dynamical oscillators and social 
networks hosting evolutionary game dynamics, and identifying partial 
differential equations for spatiotemporal dynamical systems. Situations where 
sparse optimization is effective and those in which the method fails are 
discussed. Comparisons with the traditional method of delay coordinate 
embedding in nonlinear time series analysis are given and the recent 
development of model-free, data driven prediction framework based on machine 
learning is briefly introduced.   

\end{abstract}

\maketitle

\section{Introduction} \label{sec:intro}

In nonlinear dynamics, the traditional solution to the inverse problem, i.e., 
to analyze time series to probe into the inner ``gears'' of the system, is 
based on the paradigm of delay-coordinate 
embedding~\cite{Takens:1981,PCFS:1980}. The research started about four 
decades ago when Takens~\cite{Takens:1981} proved that the underlying 
dynamical system can be faithfully reconstructed from time series with a
one-to-one correspondence between the reconstructed and the true but unknown 
dynamical systems. From the reconstructed system, statistical quantities 
characterizing the dynamical invariant set of the original system can be 
assessed~\cite{KS:book,Hegger:book}. For example, from time series,
the fractal dimensions of the underlying chaotic attractor can be 
estimated~\cite{GP:1983,Grass:1986,Proc:1988,OP:1989,Lorenz:1991,DGOSY:1993,
LLH:1996,LL:1998}, as well as the Lyapunov exponents~\cite{WSSV:1985,Sano:1985,
ER:1985,EKRC:1986,BBA:1991,STY:1998,SY:1999} and some unstable periodic 
orbits~\cite{LK:1989,BBFFHPRS:1994,PM:1995,PM:1996a,SOSKSG:1996,AM:1997}.
The continuity and differentiability of the original dynamical system can 
be tested~\cite{PCH:1995,PC:1996,PCH:1997,GPCR:2001,PMNC:2007}. Practical 
issues on determining the basic parameters of delay-coordinate
embedding such as the proper time delay~\cite{Theiler:1986,LS:1989,LPS:1991,
BP:1992,KF:1993,RCdL:1994,LLH:1996,LL:1998} and the embedding 
dimension~\cite{SYC:1991} were addressed. The Takens' paradigm was also 
extended to dynamical systems in the regime of transient 
chaos~\cite{LT:book,JFT:1994,JT:1994,DLK:2000,DLK:2001,TBS:2003} and to
systems with a time delay~\cite{TC:2007}. 

There were previous methods on data based identification and forecasting of 
nonlinear dynamical systems~\cite{FS:1987,Casdagli:1989,SGMCPW:1990,KR:1990,GS:1990,Gouesbet:1991,TE:1992,BBBB:1992,Longtin:1993,Murray:1993,Sauer:1994,Sugihara:1994,FK:1995,Parlitz:1996,SSCBS:1996,Szpiro:1997,HKS:1999,HKMS:2000,Sello:2001,MNSSH:2001,Smith:2002,Judd:2003,Sauer:2004,TZJ:2007}.
One approach is to approximate a nonlinear system by a large collection of 
linear equations in different regions of the phase space to reconstruct the 
Jacobian matrices on a proper grid~\cite{FS:1987,Gouesbet:1991,Sauer:1994} or 
to fit ordinary differential equations to chaotic data~\cite{BBBB:1992}. 
Approaches based on chaotic synchronization~\cite{Parlitz:1996} or genetic 
algorithms~\cite{Szpiro:1997,TZJ:2007} to system-parameter estimation were 
also investigated. In most studies, short-term predictions can be achieved.
For nonstationary systems, the method of over-embedding was 
introduced~\cite{HKMS:2000} in which the time-varying parameters were treated 
as independent dynamical variables so that the essential aspects of determinism
of the underlying system can be identified. 

Takens' embedding paradigm, while having evolved into a powerful and effective
framework in the past forty years to address the inverse problem in nonlinear 
dynamical systems, gives only a topological equivalent of the system of 
interest: it does not give the equations of motion of the original system. As 
such, the state evolution cannot be predicted, nor critical transitions 
leading to a possible system collapse upon parameter variations. The inverse 
problem of determining the system equations from measurements was addressed 
in 1987 by Crutchfield and McNamara~\cite{CM:1987}, who extended the notion 
of qualitative information contained in a sequence of observations to deduce 
the effective equations of motion that give the deterministic portion of the 
observed behavior. The effective equations, however, may contain additional 
terms that are not present in the original equations. The idea of using the 
inverse Frobenius–Perron problem through $L^{\infty}$ to design a dynamical 
system that is ``near'' the original system with a desired invariant density 
was proposed in 2000 by Bollt~\cite{Bollt:2000}. Modeling and nonlinear 
parameter estimation using the least-squares ($L^2$) best approximation 
(Kronecker product representation) were articulated and analyzed in 2007 by 
Yao and Bollt~\cite{YB:2007}. In the past decade, the idea of using sparse 
($L^1$) optimization methods such as compressive sensing~\cite{CRT:2006a,
CRT:2006b,Candes:2006,Donoho:2006,Baraniuk:2007,CW:2008} was articulated~\cite{WYLKG:2011,WLGY:2011,WYLKH:2011,SNWL:2012,SWL:2012,SWL:2014,SLWD:2014,SWFDL:2014,SWWL:2016}
for discovering the exact equations of motion for certain class of nonlinear 
and complex dynamical systems~\cite{EMB}. Quite recently, entropic
regression for overcoming the problem of outliers in nonlinear system 
identification was articulated~\cite{ASB:2020}. 

The principle of sparse optimization for finding the equations of nonlinear 
dynamical systems from data was first published~\cite{WYLKG:2011} 
in 2011. The basic idea is that many nonlinear dynamical systems in nature and 
engineering are governed by smooth functions that can be approximated by 
series expansions. The inverse problem of determining the system equations 
then boils down to estimating the coefficients in the series. If the series 
contain many high-order terms, the total number of coefficients to be 
estimated will be large. In this case, there is no advantage to use the series 
expansions and the problem remains as difficult as with the original equations.
However, if the system equations are relatively simple in the sense that most 
coefficients in the series expansions are zero, as in many classical nonlinear 
dynamical systems, the vector of all the coefficients to be determined will be 
sparse, rendering applicable and effective sparse optimization methods such as
compressive sensing~\cite{CRT:2006a,CRT:2006b,Donoho:2006,Baraniuk:2007,CW:2008}
originally developed in the field of signal processing in engineering and 
applied mathematics. A virtue common to sparse optimization methods is the 
low requirement of observation data. In addition to enabling finding equations 
of nonlinear dynamical systems from data~\cite{WYLKG:2011,YLG:2012}, 
compressive sensing has also been exploited for reconstruction of complex 
networks with discrete and continuous time nodal 
dynamics~\cite{WYLKG:2011,WYLKH:2011} and evolutionary game 
dynamics~\cite{WLGY:2011}, for detecting hidden 
nodes~\cite{SWL:2012,SLWD:2014}, for predicting and controlling network 
synchronization dynamics~\cite{SNWL:2012}, and for reconstructing spreading 
dynamics based on binary data~\cite{SWFDL:2014}.

\section{Principle of discovering system equations from data based on sparse optimization}

Compressive sensing solves the following convex optimization problem:
\begin{equation} \label{eq:L30_CS}
\min \|\mathbf{a}\|_1 \quad \mbox{subject \ to} \quad
\mathcal{G}\cdot \mathbf{a}=\mathbf{X},
\end{equation}
where $\mathbf{a}$ is a sparse vector to be determined, $\mathcal{G}$ is a 
known random projection matrix, $\mathbf{X}$ is a measurement vector 
constructed from the available data, and 
$\|\mathbf{a}\|_1=\sum_{i=1}^{N}|\mathbf{a}_i|$ denotes the $L_1$ norm of 
the vector $\mathbf{a}$. Compressive sensing is a paradigm of high-fidelity 
signal reconstruction using only sparse data~\cite{CRT:2006a,CRT:2006b,
Donoho:2006,Baraniuk:2007,CW:2008}, which was originally developed to solve 
the problem of transmitting massive data sets, such as those collected from a
large-scale sensor network. Because of the high dimensionality, direct 
transmission of such data sets would require a broad bandwidth. However, a 
common situation with sensor networks is that, most of the time, majority 
of the sensors are inactive so that the data set collected from the entire 
network at a time is sparse. For example, say a data set of $N$ points is 
represented by an $N\times 1$ vector ${\bf a}$, where $N$ is a large integer. 
Since ${\bf a}$ is sparse, most of its entries are zero and only a small number
of $k$ entries are non-zero, where $k \ll N$. One can use a Gaussian random 
matrix $\mathcal{G}$ of dimension $M\times N$ to obtain an $M\times 1$ vector 
${\bf X}$: ${\bf X} = \mathcal{G}\cdot {\bf a}$, where $M \sim k$. Because 
the dimension of ${\bf X}$ is much smaller than that of the original vector 
${\bf a}$, transmitting ${\bf X}$ would require a much smaller bandwidth, 
provided that ${\bf a}$ can be reconstructed at the receiver end of the 
communication channel.

The problem of reconstructing the equations of a nonlinear dynamical system
from data can be formulated in the compressive sensing 
framework~\cite{WYLKG:2011,YLG:2012}. Consider a general dynamical system 
described by 
\begin{equation} \label{eq:L30_model}
\frac{ d\mathbf{x}}{dt} = \mathbf{F}(\mathbf{x}),
\end{equation}
where $\mathbf{x}$ is an $m$-dimensional vector: 
$\mathbf{x} \equiv \left(x_1,x_2,\ldots,x_m\right)^T$, 
and $\mathbf{F}(\mathbf{x})$ represents the velocity (vector) field of the 
system with $m$ components: 
$\mathbf{F}(\mathbf{x}) = \left[F_1(\mathbf{x}),F_2(\mathbf{x}),\ldots,F_m(\mathbf{x})\right]^T$. 
The goal is to determine $\mathbf{F}(\mathbf{x})$ from limited measured time 
series $\mathbf{x}(t)$. The basic idea~\cite{WYLKG:2011} is to expand the
velocity field as a multi-dimensional power series. In particular, the $j$th
component of the velocity field can be expanded to order $q$ as  
\begin{equation} \label{eq:L30_F_j}
F_j(\mathbf{x}) = \sum^{q}_{l_1=0}\sum^q_{l_2=0}\ldots\sum^q_{l_m=0} (a_j)_{l_1l_2\ldots l_m}x_1^{l_1}x_2^{l_2}\ldots x_m^{l_m},
\end{equation}
where $(a_j)_{l_1l_2\ldots l_m}$ ($l_1,\ldots,l_m = 1,\ldots,q$) constitute 
the set of $(1+q)^m$ coefficients to be determined from the measurements. 
If this set of coefficients is dense in the sense that most of them are 
nontrivial, then resorting to the power series expansion as in 
Eq.~(\ref{eq:L30_F_j}) does not lead to any step closer to the solution 
because the coefficients cannot be determined based on limited measurements. 
However, if the coefficient set is sparse in that most of its elements are 
zero, then it would be possible to use compressive sensing to uniquely solve 
for the nontrivial coefficients. Some well studied nonlinear dynamical 
systems, such as the classic Lorenz~\cite{Lorenz:1963} and  
R\"{o}ssler~\cite{Rossler:1976} oscillators, fall into this ``sparse'' 
category because their vector fields contain only a few power series terms.  

\begin{figure}[ht!]
\centering
\includegraphics[width=0.8\linewidth]{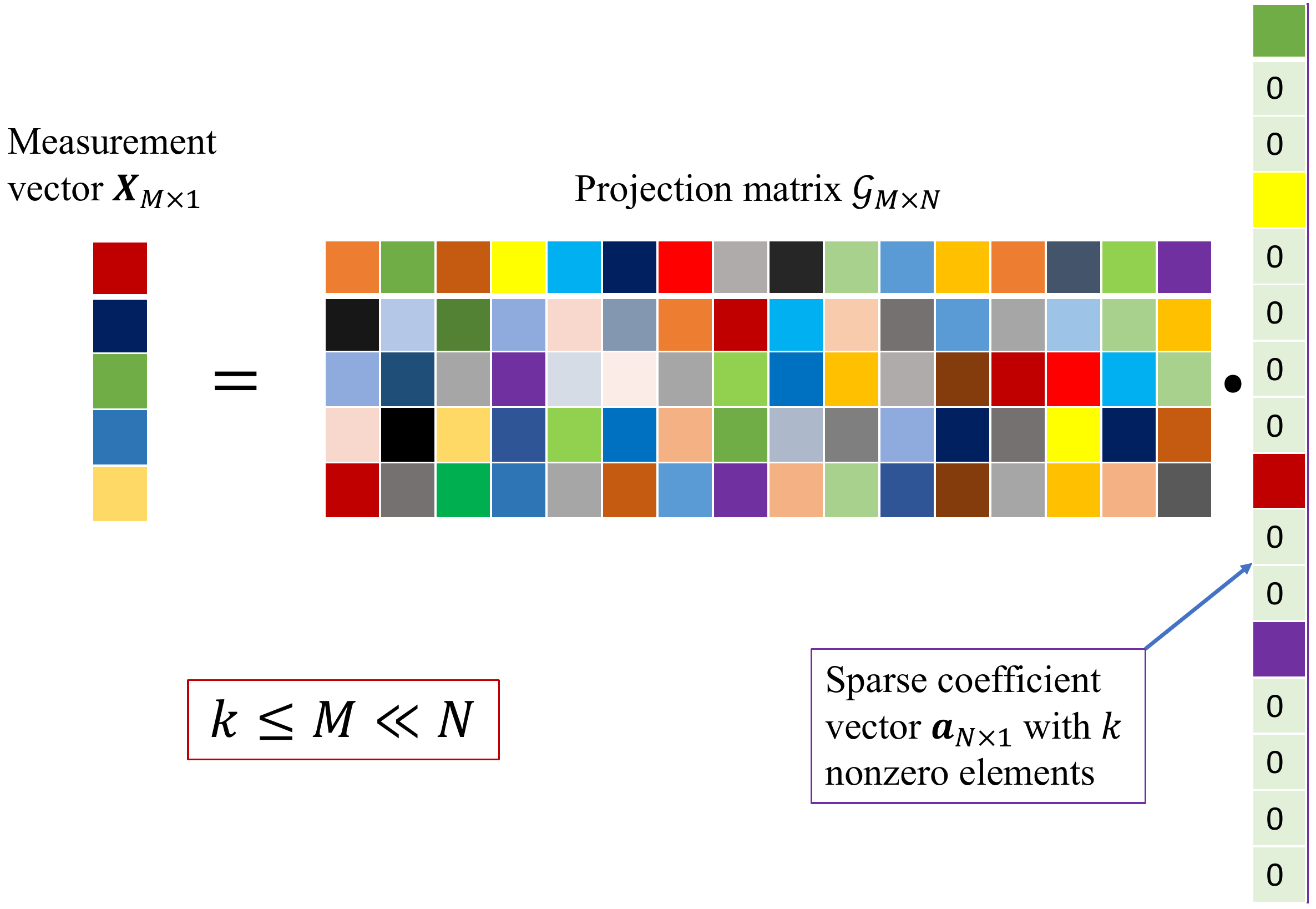}
\caption{ Standard form of compressive sensing. The goal is to obtain the 
optimal solution of the $N$-dimensional coefficient vector $\mathbf{a}$ from 
the $M$-dimensional measurement vector $\mathbf{X}$ through the projection 
matrix $\mathcal{G}$ that is $M\times N$ dimensional, where $M \ll N$. The 
mathematical framework of compressive sensing guarantees an optimal solution
insofar as the coefficient vector to be solved is sparse: it has only $k$ 
nontrivial elements, where $k \le M$, provided that the projection matrix is
random.}
\label{fig:L30_CS_form}
\end{figure}

To better understand the mathematical structure of the problem formulation,
consider the concrete case of a three-dimensional phase space ($m=3$) and 
power series expansion up to order three ($n=3$). In this case, the total 
number of unknown coefficients is $(1+n)^m = 64$. For convenience, let the 
dynamical variables be $\mathbf{x}(t)\equiv [x(t),y(t),z(t)]^T$. The first 
component of the vector field can be written as
\begin{equation} \label{eq:L30_F_1_33}
F_1(\mathbf{x}) = (a_1)_{000}x^0y^0z^0 + (a_1)_{100}x^1y^0z^0 + \ldots + (a_1)_{333}x^3y^3z^3.
\end{equation}
The $N = 64$ coefficients to be determined can be organized into a vector, a
$N\times 1$-dimensional column vector:
\begin{equation} \label{eq:L30_a_1}
\mathbf{a}_1 \equiv \left( \begin{array}{c}
(a_1)_{000} \\ (a_1)_{100} \\ \ldots \\ (a_1)_{333} \end{array} 
\right).
\end{equation} 
Defining all the combinations of the powers of the dynamical variables in 
Eq.~(\ref{eq:L30_F_1_33}) as a $1\times 64$-dimensional row vector:
\begin{equation} \label{eq:L30_g}
\mathbf{g}(t) \equiv \left[ x^0(t)y^0(t)z^0(t),x^1(t)y^0(t)z^0(t),\ldots,x^3(t)y^3(t)z^3(t) \right],
\end{equation}
we can write Eq.~(\ref{eq:L30_F_1_33}) as
\begin{equation} \label{eq:L30_F_1_33_2}
F_1[\mathbf{x}(t)] = \mathbf{g}(t)\cdot \mathbf{a}_1.
\end{equation}
Suppose measurements of the dynamical variables $\mathbf{x}(t)$ are available
at $(M+1)$ time instants $t_0,t_1,\ldots,t_{M}$, where $M\ll N$.  We have
\begin{equation} \label{eq:L30_x_measurements}
\begin{array}{c}
dx(t_1)/dt = F_1[\mathbf{x}(t_1)] = \mathbf{g}(t_1)\cdot \mathbf{a}_1  \\ 
dx(t_2)/dt = F_1[\mathbf{x}(t_2)] = \mathbf{g}(t_2)\cdot \mathbf{a}_1  \\ 
 \ldots  \\ 
dx(t_{M})/dt  =  F_1[\mathbf{x}(t_{M})] = \mathbf{g}(t_{M})\cdot \mathbf{a}_1. \end{array} 
\end{equation}
The derivatives can be estimated from the measurements: 
\begin{displaymath}
	\frac{dx(t_i)}{dt} \approx \frac{x(t_i)-x(t_{i-1})}{\delta t},
\end{displaymath}
for $i = 1, \ldots, M$. All the $M$ derivatives can be organized into an
{\em measurement} vector that is $M\times 1$-dimensional:
\begin{equation} \label{eq:L30_measurement_vector}
	\mathbf{X} = \left( \begin{array}{c}
		dx(t_1)/dt \\
		dx(t_2)/dt \\
		\ldots  \\
	dx(t_{M})/dt \end{array} \right).
\end{equation} 
Likewise, the $M$ row vectors $\mathbf{g}(t_1),\ldots,\mathbf{g}(t_M)$, each
being $1\times N$ dimensional, can be organized into an $M\times N$ dimensional
matrix as
\begin{equation} \label{eq:L30_measurement_matrix}
        \mathcal{G} = \left( \begin{array}{c}
                \mathbf{g}(t_1) \\
                \mathbf{g}(t_2) \\
                \ldots  \\
        \mathbf{g}(t_M) \end{array} \right).
\end{equation}
Equation~(\ref{eq:L30_x_measurements}) can then be written as
\begin{equation} \label{eq:L30_CS_form_x}
\mathbf{X}_{M\times 1} = \mathcal{G}_{M\times N}\cdot (\mathbf{a}_1)_{N\times 1}.
\end{equation}
The structure of Eq.~(\ref{eq:L30_CS_form_x}) is shown in 
Fig.~\ref{fig:L30_CS_form}, where $\mathcal{G}$ is the projection 
matrix. Suppose that the coefficient vector $\mathbf{a}_1$ is sparse: it has 
$k$ nontrivial elements, where $k \ll N$. If the inequality $k \le M \ll N$ 
holds, then Eq.~(\ref{eq:L30_CS_form_x}) is in the standard form of compressive 
sensing~\cite{CRT:2006a,CRT:2006b,Donoho:2006,Baraniuk:2007,CW:2008}, which
is Eq.~(\ref{eq:L30_CS}). 

Two remarks are in order. 

Firstly, the coefficient vector $\mathbf{a}_1$ is associated with the velocity
field of the first dynamical variable $x$. For any other dynamical variable
in Eq.~(\ref{eq:L30_model}), a similar expansion procedure can be carried out.
For example, for the second dynamical variable $y$, the coefficient vector 
$\mathbf{a}_2$ of its velocity field can be solved through:
\begin{equation} \label{eq:L30_CS_form_y}
        \mathbf{Y}_{M\times 1} = \mathcal{G}_{M\times N}\cdot (\mathbf{a}_2)_{N\times 1},
\end{equation}
where $\mathbf{Y}_{M\times 1}$ is the measurement vector constituting the 
derivatives of $y$ at the measurement points. Note that the projection matrix
$\mathcal{G}_{M\times N}$ has the same form in Eqs.~(\ref{eq:L30_CS_form_x})
and (\ref{eq:L30_CS_form_y}). 

Secondly, in the mathematical framework of compressive sensing, a requirement
is that the projection matrix $\mathcal{G}$ be random with zero correlation
among its elements. However, in the power series expansion formulation, the 
elements of this matrix are distinct combinations of the powers of all the 
dynamical variables in the system. Since the system is deterministic, nonzero 
correlations among the matrix elements are inevitable. For chaotic systems,
this violation of the randomness condition may not be as severe, since the
time evolution of the dynamical variables is effectively random, insofar as 
the time interval between the adjacent measurement points is reasonably large.
Nonetheless, there is no mathematical guarantee that an optimal solution can
be obtained for solving the inverse problem in nonlinear dynamical systems 
through the compressive sensing approach, in spite of previously demonstrated
successes~\cite{WYLKG:2011,YLG:2012,WYLKH:2011,WLGY:2011,SWL:2012,SLWD:2014,
SNWL:2012,SWFDL:2014}.

\section{Applications}

A number of applications of the compressive-sensing based solutions of 
inverse problems in nonlinear and complex dynamical systems are described.

\subsection{Predicting system collapse} \label{subsec:L30_pred_collapse}  

\begin{figure}[ht!]
\centering
\includegraphics[width=0.8\linewidth]{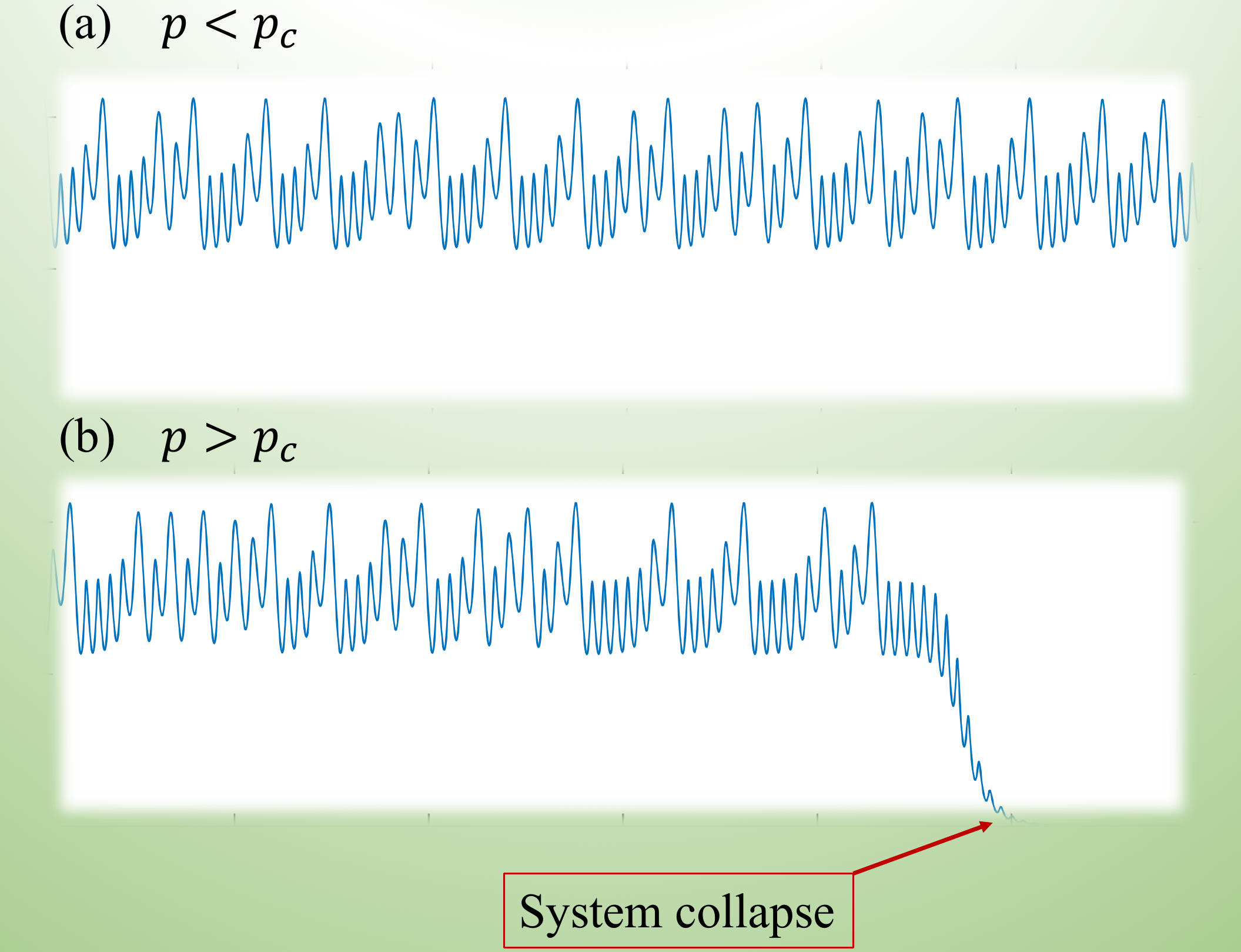}
\caption{ Dynamical behaviors before and after a crisis bifurcation. The
bifurcation occurs at the critical parameter value $p_c$. (a) Before the
crisis ($p < p_c$), the system functions normally, as the average values of 
its dynamical variables maintain at a healthy level, in spite of chaotic 
fluctuations. (b) After the bifurcation, the system exhibits transient chaos 
and eventually collapses.}
\label{fig:L30_crisis}
\end{figure}

When some parameter of a nonlinear dynamical system changes, a bifurcation 
that leads to the collapse of the system can occur. For example, in a global
bifurcation called crisis~\cite{GOY:1983a}, at the bifurcation point a chaotic 
attractor collides with its own basin boundary and is destroyed. Let $p$ be
the bifurcation parameter and $p_c$ be the crisis point. Before the crisis,
i.e., $p < p_c$, the system functions normally with a sustained chaotic 
behavior in its time evolution, as shown in Fig.~\ref{fig:L30_crisis}(a),
where the ordinate represents a typical dynamical variable of the system. 
Beyond the bifurcation point, the system collapses eventually after exhibiting 
transient chaos, as shown in Fig.~\ref{fig:L30_crisis}(b). The bifurcation is 
thus a catastrophe that must be prevented. In natural and engineering systems, 
catastrophic collapse is always a possibility. For example, in electrical 
power systems, voltage collapse~\cite{DL:1999} can occur after the system 
enters into the state of transient chaos~\cite{GOY:1983a,LT:book}. In ecology, 
slow parameter drift caused by environmental deterioration can induce a 
transition into transient chaos, followed by species 
extinction~\cite{McY:1994,HACFGLMPSZ:2018}. 
For a dynamical system of interest, predicting a catastrophe in advance of its 
occurrence is of paramount importance. This is a challenging problem when the 
system equations are unknown and the only available information that one can 
rely on to make the prediction is time series measured while the system 
still functions normally. 

If the underlying equations of the system have a simple mathematical structure,
e.g., if its velocity field consists of a few power-series or Fourier-series
terms only, then sparse optimization methods such as compressive sensing can
be exploited to identify the system equations~\cite{WYLKG:2011,WLG:2016} and
consequently to predict transitions. In particular, based on the time series,
one first predicts the system equations and then identifies the pertinent 
parameter of the system that can potentially lead to a collapse. With such
information, one can perform a computational bifurcation analysis to locate 
potential catastrophic events in the parameter space so as to determine the 
likelihood of system's drifting into a catastrophic regime. For example, if 
it is determined that the system currently operates in a parameter region 
close to the crisis bifurcation, a catastrophe may be imminent as a small 
parameter change or a random perturbation can push the system into the regime 
of transient chaos where collapse is inevitable.

The principle of compressive-sensing based prediction of system collapse
has been demonstrated with a number of nonlinear dynamical 
systems~\cite{WYLKG:2011,WLG:2016}.

\subsection{Predicting future attractors in time-varying dynamical systems}

Physical systems are constantly under external influences that can lead to 
parameter drifting. If the the time scale of the internal dynamics of the 
system is much faster than that of the external perturbation, the resulting
drift in the parameter can be regarded as adiabatic. In this case, in a 
time window of duration much longer than the internal but shorter than the 
external time scale, the system can be viewed to be in a kind of ``asymptotic 
state'' or an ``attractor,'' as for the case of a stationary dynamical system
with no time-varying parameters. However, in a long time scale, the attractor  
does depend on time. A problem of significant interest is to forecast the 
``future'' attractor of the system. For example, consider the climate system.
It is under random disturbances all the time but adiabatic parameter drifting
is also present, such as that induced by the injection of CO$_2$ into the 
atmosphere due to human activities. The time scale for any appreciable 
increase in the CO$_2$ level (e.g., months or years) is typically much slower 
than the intrinsic time scales of the system (e.g., days). The climate system 
can thus be regarded as an adiabatically time-varying, nonlinear dynamical
system. To forecast the possible future attractor of the system is key to 
sustainability, an issue that is critical to many other natural and 
engineering systems.

It was demonstrated~\cite{YLG:2012} that predicting the future attractor of 
an adiabatically time-varying dynamical system can be formulated as a problem
solvable by compressive sensing. In particular, let the system be described by
\begin{equation} \label{eq:L30_TV_system}
	d\mathbf{x}/dt = \mathbf{F}[\mathbf{x}, \mathbf{p}(t)],
\end{equation}
where $\mathbf{x}$ is the set of dynamical variables of the system in the 
$m$-dimensional phase space and $\mathbf{p}(t) \equiv [p_1(t),\ldots,p_K(t)]$ 
denotes $K$ independent, time-varying parameters. The tacit assumption is 
that both the velocity field $\mathbf{F}$ and the vector parameter function 
$\mathbf{p}(t)$ are unknown, but only time series $\mathbf{x}(t)$ measured
from the system in the time interval $t_M - T_M \le t \le t_M$ are available,
where $t_M$ is the current time. The goal is to determine the precise
mathematical forms of $\mathbf{F}$ and $\mathbf{p}(t)$ from the available 
time series at $t_M$ so that the dynamical evolution of the system and the 
likely attractors for $t > t_M$ can be computationally assessed.

As for the case of predicting system collapse discussed in 
Sec.~\ref{subsec:L30_pred_collapse}, the first step is to expand all components
of the time-dependent vector field $\mathbf{F}[\mathbf{x},\mathbf{p}(t)]$ 
into a power series in terms of both dynamical variables $\mathbf{x}$ and 
time $t$. The $i$th component $\mathbf{F}[\mathbf{x},\mathbf{p}(t)]_i$ 
of the vector field can be written as
\begin{equation} \label{eq:L30_TV_series}
\sum_{l_1,\cdots,l_m=1}^{n} [(\alpha_i)_{l_1,\cdots,l_m}x_1^{l_1} \cdots x_m^{l_m}\cdot\sum_{w=0}^{v}(\beta_i)_w t^w] \equiv \sum_{l_1,\cdots,l_m=1}^{n}\sum_{w=0}^{v} (a_i)_{l_1,\cdots,l_m;w}x_1^{l_1} \cdots x_m^{l_m}\cdot t^w
\end{equation}
where $x_k$ ($k=1,\cdots, m$) is the $k$th component of the dynamical variable, 
$a_i$ is the $i$th component of the coefficient vector to be determined, and 
the time evolution of each term can be approximated by the power series 
expansion in time, i.e., $\sum_{w=0}^{v}(\beta_i)_w t^w$. The power-series 
expansion can then be cast into the standard form of compressive sensing
Eq.~(\ref{eq:L30_CS}). If every combined scalar coefficient 
$(a_i)_{l_1,\cdots,l_m;w}$ associated with the corresponding term in 
Eq.~(\ref{eq:L30_TV_series}) can be determined from time series for 
$t \le t_M$, the vector-field component $[{\bf F}({\bf x},{\bf p}(t))]_i$ 
becomes known. Repeating the procedure for all components, the entire vector 
field for $t > t_M$ can be identified.

Note that, the predicted form of $\mathbf{F}$ and $\mathbf{p}(t)$ at time
$t_M$ would contain errors that in general will increase with time. In 
addition, for $t > t_M$ new perturbations can occur to the system so that 
the forms of $\mathbf{F}$ and $\mathbf{p}(t)$ may be further changed. It is 
thus necessary to execute the prediction algorithm frequently using time 
series available at the time. For example, the system could be monitored at 
all times so that time series can be collected, and predictions can be carried 
out at $t_i$'s, where $\ldots > t_i > \ldots > t_{M+2} > t_{M+1} > t_M$. For 
any $t_i$, the prediction algorithm is to be performed based on available 
time series in a suitable window prior to $t_i$.

\subsection{Finding complex network structure from data}

A class of inverse problems in complex dynamical networks can be stated
as follows. Suppose that the connecting structure of a sparse network and the 
nodal dynamical equations are not known but oscillatory time series can be 
measured from all nodes in the network. Can the network structure and all
the equations of motion of the system be inferred from the measurements?
Complex networks in the real world are typically sparse~\cite{Newman:book}.
It turns out that the compressive sensing based approach can be exploited for 
deciphering the connection structure and the equations of complex oscillator 
networks~\cite{WYLKH:2011}.   

An oscillator network can be viewed as a high-dimensional dynamical system 
that generates oscillatory time series at various nodes, where the local
dynamics at a node are described by
\begin{equation} \label{eq:L30_CON}
\dot{\mathbf{x}}_i = \mathbf{F}_i(\mathbf{x}_i) + \sum_{j=1,j\neq i}^{n}
{\bf C}_{ij} (\mathbf{x}_j - \mathbf{x}_i), \ \ (i=1,\cdots ,n),
\end{equation}
where $\mathbf{x}_i \in R^d$ represents the set of externally accessible 
dynamical variables of node $i$, $n$ is the number of nodes, and 
$\mathcal{C}_{ij}$ is the $d\times d$ coupling matrix between the dynamical 
variables at nodes $i$ and $j$ denoted by
\begin{equation}
\mathcal{C}_{ij} =\left(
  \begin{array}{cccc}
    c_{ij}^{1,1} & c_{ij}^{1,2} & \cdots & c_{ij}^{1,d} \\
    c_{ij}^{2,1} & c_{ij}^{2,2} & \cdots & c_{ij}^{2,d} \\
    \cdots & \cdots  & \cdots  & \cdots  \\
    c_{ij}^{d,1} & c_{ij}^{d,2} & \cdots & c_{ij}^{d,d} \\
  \end{array}
\right).
\end{equation}
In $\mathcal{C}_{ij}$, the superscripts $kl$ ($k,l=1,2,...,d$) stand for the 
coupling from the $k$th component of the dynamical variable at node $i$ to the 
$l$th component of the dynamical variable at node $j$. For any two nodes, the 
number of possible coupling terms is $d^2$. If there is at least one nonzero 
element in the matrix $\mathcal{C}_{ij}$, nodes $i$ and $j$ are coupled and,
as a result, there is a link (or an edge) between them in the network. 
Generally, more than one element in $\mathcal{C}_{ij}$ can be nonzero. On
the contrary, if all the elements of $\mathcal{C}_{ij}$ are zero, there is no 
coupling between nodes $i$ and $j$. The connecting structure and the 
interaction strengths among various nodes of the network can be identified if
the coupling matrix $\mathcal{C}_{ij}$ can be determined from time-series 
measurements.

To formulate a solution of the network inverse problem in terms of sparse
optimization, the first step is to rewrite Eq.~(\ref{eq:L30_CON}) as 
\begin{equation} \label{eq:L30_CON_2}
\dot{\mathbf{x}}_i = [\mathbf{F}_i(\mathbf{x}_i) - \sum_{j=1,j\neq i}^{n}
\mathcal{C}_{ij} \mathbf{x}_i ]  + \sum_{j=1,j\neq i}^{n} \mathcal{C}_{ij}
\mathbf{x}_j ,
\end{equation}
where the first term on the right side is exclusively a function of 
$\mathbf{x}_i$ and the second term is a function of variables of other nodes 
(couplings). The first term can be conveniently denoted as 
$\mathbf{\Gamma}_i(\mathbf{x}_i)$, which is unknown. The $k$th component of 
$\mathbf{\Gamma}_i(\mathbf{x}_i)$ can be represented by a power series of 
order up to $q$:
\begin{equation} \label{eq:L30_2_1}
\left[ \mathbf{\Gamma}_i(\mathbf{x}_i)\right]_k \equiv \left[\mathbf{F}_i(\mathbf{x}_i) - \sum_{j=1,j\neq i}^{n} \mathcal{C}_{ij} \mathbf{x}_i \right]_k = \sum_{l_1=0}^{q}\sum_{l_2=0}^{q} \cdots\sum_{l_d=0}^{q} [(\mathbf{\alpha}_i)_k]_{l_1,\cdots,l_d} [(\mathbf{x}_i)_1]^{l_1} [(\mathbf{x}_i)_2]^{l_2} \cdots [(\mathbf{x}_i)_d]^{l_d},
\end{equation}
where $(\mathbf{x}_i)_{k}$ ($k=1,\cdots, d$) is the $k$th component of the 
dynamical variable at node $i$, the total number of products is $(1+q)^d$, 
and $[(\mathbf{\alpha}_i)_k]_{l_1,\cdots,l_m}\in R^m$ is the coefficient scalar
of each product term, which is to be determined from measurements as well. 
Note that terms in Eq.~(\ref{eq:L30_2_1}) are all possible products of 
different components with different power of exponents. As an example, for 
$d=2$ (the components are $x$ and $y$) and $q=2$, the power series expansion 
is 
\begin{displaymath}
\alpha_{0,0} + \alpha_{1,0}x + \alpha_{0,1}y + \alpha_{2,0}x^2 + \alpha_{0,2}y^2
+ \alpha_{1,1}xy +\alpha_{2,1}x^2y + \alpha_{1,2} xy^2 + \alpha_{2,2}x^2y^2.
\end{displaymath}

The second step is to rewrite Eq.~(\ref{eq:L30_CON_2}) as 
\begin{equation} \label{eq:L30_CON_3}
\dot{\mathbf{x}}_i = \mathbf{\Gamma}_i(\mathbf{x}_i) + \mathcal{C}_{i1} \mathbf{x}_1 + \mathcal{C}_{i2} \mathbf{x}_2 + \cdots + \mathcal{C}_{in} \mathbf{\bf x}_n.
\end{equation}
The goal is to estimate the various coupling matrices 
$\mathcal{C}_{ij}$ $(j=1,\cdots, i-1,i+1,\cdots,n)$ and the coefficients of
$\mathbf{\Gamma}_i({\bf x}_i)$ from sparse measurements. The sparsity 
requirement of the compressive sensing theory stipulates that, to reconstruct 
the coefficients of Eq.~(\ref{eq:L30_CON_3}) from measurements, most 
coefficients must be zero. To include as many coupling forms as possible, 
one can write each term $\mathcal{C}_{ij} \mathbf{x}_j$ in 
Eq.~(\ref{eq:L30_CON_3}) as a power series in the same form of 
$\mathbf{\Gamma}_i(\mathbf{x}_i)$ but with different coefficients:
\begin{equation} \label{eq:L30_CON_4}
\dot{\mathbf{x}}_i = \mathbf{\Gamma}_1(\mathbf{x}_1) + \mathbf{\Gamma}_2(\mathbf{x}_2)+ \cdots + \mathbf{\Gamma}_n(\mathbf{x}_n). 
\end{equation}
This setting not only includes many possible coupling forms but also ensures 
that the sparsity condition is satisfied so that the prediction problem can be 
formulated in the compressive-sensing framework. For an arbitrary node $i$, 
information about the node-to-node coupling, or about the network connectivity, 
is contained completely in $\mathbf{\Gamma}_j(j\neq i)$. For example, if in 
the equation of $i$, a term in $\mathbf{\Gamma}_j(j\neq i)$ is not zero, there 
then exists a coupling between $i$ and $j$ with the strength given by the
coefficient of the term. Subtracting the coupling terms 
$-\sum_{j=1,j\neq i}^{n} \mathcal{C}_{ij}{\bf x}_i$ from $\mathbf{\Gamma}_i$ 
in Eq.~(\ref{eq:L30_2_1}), which is the sum of the coupling coefficients of 
all $\mathbf{\Gamma}_j$ $(j\neq i)$, the local system equations 
$\mathbf{F}_i(\mathbf{x}_i)$ can be obtained. Therefore, once the coefficients 
of Eq.~(\ref{eq:L30_CON_4}) have been determined, the nodal dynamical equations
and the couplings among the nodes are all known. As explained in 
Sec.~(\ref{subsec:L30_pred_collapse}), the power series expansion coefficients 
of Eq.~(\ref{eq:L30_CON_4}) can be determined through compressive sensing. 

\subsection{Finding social network structure from evolutionary game data}

In the summer of 2011, a small social network experiment was conducted at 
Arizona State University, where 22 students from different Schools were 
invited to test the effectiveness of a sparse-optimization based approach
to mapping out the structure of social networks. In particular, the 22
participants constituted a social network, where any individual has a few
acquaintances/friends in the group, but many participants had not known 
each other prior to the experiment. The network is thus sparse. During 
the experiment, all participants were asked to play an evolutionary game,
the prisoner's dilemma game, with their friends for about 30 runs. That is, 
each and every agent (node) in the network was asked to play the game but only 
with his/her direct neighbors. For each run of the play, the strategies used 
by the opponents of each and every pair of players were recorded, together 
with the outcome (i.e., the winner and loser). Based on the data from all 30 
runs, a compressive sensing based algorithm was executed, yielding a network 
structure that matches exactly with that of the actual social network. In 
fact, it was demonstrated that data from about 15 runs were already sufficient 
to infer the underlying network structure with $100\%$ 
accuracy~\cite{WLGY:2011}.

The theoretical underpinning of this successful experiment lies in formulating
the dynamical process of evolutionary game on a network as a problem of 
sparse optimization. Specifically, in an evolutionary game, the players use 
different strategies in order to gain the maximum payoff which, in general,
can be divided into two types: cooperation and defection. It was 
shown~\cite{WLGY:2011} that, with limited data on each player's strategy and 
payoff, a compressive-sensing based framework can be developed to yield precise 
knowledge of the node-to-node interaction patterns in an efficient manner.    

A sketch of the principle underlying the compressive sensing framework is as
follows. In an evolutionary game, at any time a player can choose one of the 
two strategies $S$: cooperation (C) or defection (D), which can be expressed as 
$\mathbf{S}(C) = (1, 0)^T$ and $\mathbf{S}(D) = (0, 1)^T$. The payoffs of two 
players in a game are determined by their strategies and the payoff matrix of 
the specific game. For example, for the prisoner's dilemma game 
(PDG)~\cite{NM:1992} and the snowdrift games (SG)~\cite{HD:2004}, the payoff 
matrices are given, respectively, by
\begin{eqnarray}
\mathcal{P}_{PDG} = \left (%
\begin{array}{cc}
  1 & 0 \\
  b & 0 \\
\end{array}
\right) \ \ \hbox{or} \ \
\mathcal{P}_{SG} = \left (%
\begin{array}{cc}
  1 & 1-r \\
  1+r & 0 \\
\end{array}
\right),
\end{eqnarray}
where $b$ ($1<b<2$) and $r$ ($0<r<1$) are parameters characterizing the 
temptation to defect. When a defector encounters a cooperator, the defector 
gains payoff $b$ in the PDG and payoff $1+r$ in the SG, but the cooperator 
gains the sucker payoff 0 in the PDG and payoff $1-r$ in the SG. At each time
step, all agents play the game with their neighbors and gain payoffs. For 
agent $i$, the payoff is
\begin{equation}
P_i = \sum_{j \in \Gamma_i}\mathbf{S}_i^T \cdot \mathcal{P} \cdot \mathbf{S}_j,
\end{equation}
where $\mathbf{S}_i$ and $\mathbf{S}_j$ are the strategies of agents $i$ and 
$j$ at the time and the sum is over the neighboring set $\Gamma_i$ of $i$. 
After obtaining its payoff, an agent updates its strategy according to its own 
and neighbors' payoffs, attempting to maximize its payoff at the next round. 
Possible mathematical prescriptions to describe quantitatively an agent's 
decision making process include the best-take-over rule~\cite{NM:1992}, the 
Fermi rule~\cite{ST:1998}, and one based on the payoff-difference-determined 
updating probability~\cite{SSP:2008}. For example, the Fermi rule is defined,
as follows. After player $i$ randomly chooses a neighbor $j$, $i$ adopts 
$j$'s strategy $\mathbf{S}_j$ with the probability~\cite{ST:1998}:
\begin{equation}
W(\mathbf{S}_i\leftarrow
\mathbf{S}_j)=\frac{1}{1+\exp{[(P_i-P_j)/\kappa]}},
\end{equation}
where $\kappa$ characterizes the stochastic uncertainties in the evolutionary
game dynamics: $\kappa=0$ corresponds to absolute rationality where the 
probability is zero if $P_j < P_i$ and one if $P_i < P_j$, and 
$\kappa \rightarrow \infty$ indicates completely random decision making. The 
probability $W$ thus characterizes the bounded rationality of agents in 
society and the natural selection based on the relative fitness in evolution.

The key to solving the inverse problem of network reconstruction is the
relationship between agents' payoffs and strategies. The interactions among 
the agents in the network can be characterized by an $n \times n$ adjacency 
matrix $\mathcal{A}$ with elements $a_{ij}=1$ if agents $i$ and $j$ are 
connected and  $a_{ij}=0$ otherwise. The payoff of agent $x$ can be expressed 
by
\begin{eqnarray}
\nonumber 
P_x(t) & = & a_{x1}\mathbf{S}_x^T(t) \cdot \mathcal{P} \cdot \mathbf{S}_1(t) + \cdots  + a_{x,x-1}\mathbf{S}_x^T(t)\cdot \mathcal{P} \cdot \mathbf{S}_{x-1}(t) + a_{x,x+1}\mathbf{S}_x^T(t) \cdot \mathcal{P} \cdot \mathbf{S}_{x+1}(t) \\ \label{eq:_L30_Game_main}
	& + & \cdots + \ a_{xn}\mathbf{S}_x^T(t) \cdot \mathcal{P} \cdot \mathbf{S}_n(t),
\end{eqnarray}
where $a_{xi}$ ($i=1, \cdots, x-1, x+1, \cdots, n$) represents a possible 
connection between agent $x$ and its neighbor $i$,
$a_{xi}\mathbf{S}_x^T(t) \cdot \mathbf{P} \cdot \mathbf{S}_i(t)$
($i=1, \cdots, x-1, x+1, \cdots, n$) stands for the possible payoff of agent 
$x$ from playing the game with $i$ (if there is no connection between 
$x$ and $i$, the payoff is zero because $a_{xi}=0$), and $t=1, \cdots, M$ is 
the number of rounds that all agents play the game with their neighbors. This 
relation provides a base to construct the vector $\mathbf{X}_x$ and matrix 
$\mathcal{G}_x$ in a proper compressive-sensing framework to obtain the 
solution of the neighboring vector $\mathbf{A}_x$ of agent $x$. In particular,
one can define 
\begin{eqnarray}
\mathbf{X}_x &\equiv& (P_x(t_1),P_x(t_2),\cdots,P_x(t_M))^T,
\nonumber \\
\mathbf{A}_x &\equiv& (a_{x1},\cdots ,
a_{x,x-1},a_{x,x+1}, \cdots , a_{xn})^T,
\end{eqnarray}
and 
\begin{eqnarray}
\mathcal{G}_x \equiv \left(
\begin{array}{cccccc}
 F_{x1}(t_1)& \cdots & F_{x,x-1}(t_1) & F_{x,x+1}(t_1) &\cdots & F_{xn}(t_1) \\
 F_{x1}(t_2)& \cdots & F_{x,x-1}(t_2) & F_{x,x+1}(t_2) &\cdots & F_{xn}(t_2) \\
 \vdots & \cdots & \vdots & \vdots & \vdots & \vdots \\
 F_{x1}(t_M)& \cdots & F_{x,x-1}(t_M) & F_{x,x+1}(t_M) &\cdots & F_{xn}(t_M) \\
\end{array}
\right), \nonumber
\end{eqnarray}
where $F_{xy}(t_i)=\mathbf{S}_x^T(t_i)\cdot\mathcal{P}\cdot\mathbf{S}_y(t_i)$. 
The relation among the vectors $\mathbf{X}_x$, $\mathbf{A}_x$, and the matrix 
$\mathcal{G}_x$ is given by exactly the same form as of compressive sensing
Eq.~(\ref{eq:L30_CS}):  
\begin{equation}
\mathbf{X}_x = \mathcal{G}_x \cdot \mathbf{A}_x,
\end{equation}
where $\mathbf{A}_x$ is sparse due to the sparsity of the underlying network, 
making the compressive-sensing framework applicable. Since 
$\mathbf{S}_x^T(t_i)$ and $\mathbf{S}_y(t_i)$ in $F_{xy}(t_i)$ come from data 
and $\mathcal{P}$ is known, the vector $\mathbf{X}_x$ can be obtained directly 
while the matrix $\mathcal{G}_x$ can be calculated from the strategy and 
payoff data. The vector $\mathbf{A}_x$ can thus be predicted based solely on 
the time series game data. Since the self-interaction terms $a_{xx}$ are not 
included in the vector $\mathbf{A}_x$ and the self-column 
$[F_{xx}(t_1),\cdots ,F_{xx}(t_M)]^T$ is excluded from the matrix 
$\mathcal{G}_x$, the computation required for compressive sensing can be 
reduced. In a similar fashion, the neighboring vectors of all other agents can
be predicted, yielding the network adjacency matrix 
\begin{displaymath}
\mathcal{A} = (\mathbf{A}_{1},\mathbf{A}_{2},\cdots,\mathbf{A}_{n})
\end{displaymath}
and, hence, the connection structure of the underlying social network. 

\subsection{Sparse optimization based on LASSO and applications in reconstructing complex networks with binary-state dynamics}

A generalization of the compressive sensing approach to inverse problems in 
nonlinear and complex dynamical systems is sparse optimization based on LASSO 
(Least Absolute Shrinkage and Selection Operator)~\cite{HSWD:2015,LSWGL:2017}.
In statistical and machine learning, LASSO is a regression method that embodies 
both variable selection and regularization to enhance the prediction accuracy 
and interpretability of the statistical model it 
produces~\cite{SS:1986,FHT:book}. In particular, LASSO incorporates an 
$L_1$-norm and an error control term to solve the sparse vector $\mathbf{a}$ 
according to the constraint $\mathcal{G}\cdot \mathbf{a}=\mathbf{X}$ [as in 
Eq.~(\ref{eq:L30_CS})] from a small amount of data by optimizing
\begin{equation} \label{eq:L30_lasso}
\min_{{\bf a}} \Bigl\lbrace \dfrac{1}{2M}\Vert \mathcal{G}\cdot \mathbf{a} - \mathbf{X}\Vert_2^2 + \lambda\Vert\mathbf{a}\Vert_1 \Bigr\rbrace \text{,}
\end{equation}
where $\Vert\mathbf{a}\Vert_1$ is the $L_1$ norm of $\mathbf{X}$ assuring the 
sparsity of the solution, the least squares term 
$\Vert\mathcal{G}\cdot \mathbf{a} - \mathbf{X}\Vert_2^2$ guarantees the 
robustness of the solution against noise in the data, and $\lambda$ is a 
nonnegative regularization parameter that affects the reconstruction 
performance in terms of the sparsity of the network, which can be determined 
by a cross-validation method~\cite{PVGMTGBPWD:2011}. The advantage of LASSO is 
similar to that of compressive sensing: the number $M$ of bases (measurements) 
needed can be much less than the length of $\mathbf{a}$. The LASSO-based 
sparse optimization method was successfully applied to reconstructing the 
structures of complex networks~\cite{HSWD:2015}. 

In data based reconstruction of complex networks, a difficult problem is when 
the nodal dynamical states are discrete and are of the binary type - a 
situation that arises commonly in nature, technology, and 
society~\cite{barrat2008}. In a networked system hosting binary nodal dynamics,
each node can be in one of the two possible states, e.g., being active or 
inactive in neuronal and gene regulatory networks~\cite{KRA:2010}, cooperation 
or defection in networks hosting evolutionary game dynamics~\cite{SF:2007}, 
being susceptible or infected in epidemic spreading on social and technological
networks~\cite{PSCVV:2015}, and two competing opinions in social
communities~\cite{SHS:2009}, etc. The interactions among the nodes are complex 
and a state change can be triggered either deterministically (e.g., depending 
on the states of their neighbors) or randomly. Indeed, deterministic and 
stochastic state changes can account for a variety of emergent phenomena, such 
as the outbreak of epidemic spreading~\cite{GGA:2013}, cooperation among 
selfish individuals~\cite{SP:2005}, oscillations in biological 
systems~\cite{KVK:2013}, power blackout~\cite{BPPSH:2010}, financial 
crisis~\cite{GDB:2013}, and phase transitions in natural systems~\cite{BV:2011}.
A variety of models have been introduced to gain insights into binary-state 
dynamics on complex networks~\cite{Newman:book}, such as the voter models
for competition of two opinions~\cite{SR:2005}, stochastic propagation
models for epidemic spreading~\cite{PSV:2001}, models of rumor diffusion
and adoption of new technologies~\cite{CFL:2009}, cascading failure 
models~\cite{BBB:2013}, Ising spin models for ferromagnetic
phase transition~\cite{KRB:book}, and evolutionary games for cooperation
and altruism~\cite{SSP:2008}. A general theoretical approach to dealing
with networks hosting binary state dynamics was developed~\cite{Gleeson:2013} 
based on the pair approximation and the master equations, providing a good 
understanding of the effect of the network structure on the emergent phenomena.

The problem of reconstructing complex networks with binary-state dynamics
is challenging, for three reasons. Firstly, the basic nodal dynamics are 
governed by the probability to transition between the two states - the 
switching probability of a node, which depends on the states of its neighbors 
according to a variety of functions that can be linear,
nonlinear, piecewise, or stochastic. If the function that governs the switching
probability is unknown, it would be difficult to obtain a 
solution of the reconstruction problem. Secondly, structural information is 
often hidden in the binary states of the nodes in an unknown manner and the 
dimension of the solution space can be high, rendering impractical 
(computationally prohibitive) brute-force enumeration of all possible network 
configurations. Thirdly, the presence of measurement noise, missing data, and 
stochastic effects in the switching probability make the reconstruction task 
even more challenging, calling for the development of effective methods that 
are robust against internal and external random effects.

in 2017, a general and robust framework for reconstructing complex networks 
based solely on the binary states of the nodes without any knowledge about the 
switching functions was developed~\cite{LSWGL:2017}. The idea was centered 
about developing a general method to linearize the switching functions from 
binary data. The data-based linearization method was demonstrated to be 
applicable to linear, nonlinear, piecewise, or stochastic switching functions. 
The method allows one to convert the network reconstruction problem into a 
sparse signal reconstruction problem for local structures associated with 
each node. In particular, because of the natural sparsity of complex networks,
LASSO was used~\cite{HSWD:2015,LSWGL:2017} to identify the neighbors of 
each node in the network from sparse binary data contaminated by noise. 
The linearization procedure was justified through a number of linear, nonlinear
and piecewise binary-state dynamics on a large number of model and real complex
networks. For all the models tested, universally high reconstruction accuracy
was achieved~\cite{LSWGL:2017} even for small data amount with noise. Because 
of its high accuracy, efficiency and robustness against noise and missing data,
the reconstruction framework can serve as a general solution to the inverse
problem of network reconstruction from binary-state time series, which is
key to articulating effective strategies to control complex networks with
binary state dynamics.

\subsection{Discovering models of spatiotemporal dynamical systems from data}

There was an early work on modeling and parameter estimation for coupled
oscillators and spatiotemporal dynamical systems based on the Kronecker
product presentation~\cite{YB:2007}. In recent years, sparse optimization
or learning has been applied to discovering models of spatiotemporal systems 
described by partial differential equations 
(PDEs)~\cite{RBPK:2017,LLYTVKY:2019,RG:2019,GRG:2019,RGG:2020}. PDEs for 
spatiotemporal systems in science and engineering have the general form:
\begin{equation} \label{eq:L20_PDE_general}
	\sum^{n}_{i=0} c_i \mathbf{f}_i(\mathbf{u},\partial_t\mathbf{u},\nabla\mathbf{u},\nabla^2\mathbf{u},\ldots) = 0,
\end{equation}
where $c_i$'s are constant coefficients and $\mathbf{f}_i$'s are vector 
functions of time and space derivatives of various orders of the vector 
field $\mathbf{u}$. Usually, symmetry and physical considerations can be used
to reduce the number of terms in Eq.~(\ref{eq:L20_PDE_general}) and to narrow 
down the possible functional forms of $\mathbf{f}_i$. From the measurement
(noisy) spatiotemporal data $\mathbf{u}$, sparse optimization can be used to 
remove the unnecessary terms and yield a model that contains a small number of 
terms~\cite{RBPK:2017,LLYTVKY:2019,RG:2019,GRG:2019,RGG:2020}.

In Ref.~\cite{RGG:2020}, the following procedure was devised. Consider a system 
that contains only a single term of the first-order time derivative of the 
vector field, written as
\begin{equation} \label{eq:L20_PDE_general_2}
\partial_t \mathbf{u} = \sum^{n}_{i=1} c_i \mathbf{f}_i(\mathbf{u},\nabla\mathbf{u},\nabla^2\mathbf{u},\ldots).
\end{equation}
To obtain a set of linear equations for sparse optimization, one multiplies a 
weight vector $\mathbf{w}$ to Eq.~(\ref{eq:L20_PDE_general_2}) and integrates 
both sides over a number of distinct spatiotemporal domains $\Omega_l$ 
($l=1,\ldots,L$) to convert Eq.~(\ref{eq:L20_PDE_general_2}) to
\begin{equation} \label{eq:L20_PDE_linear_form}
	\mathbf{q}_0 = \sum^{n}_{i=1} c_i \mathbf{q}_i = \mathcal{Q}\cdot\mathbf{c},
\end{equation}
where $\mathbf{c} \equiv (c_1,\ldots,c_n)^T$ is the coefficient vector to 
be determined from data, $\mathbf{q}_i$ is an $L$-dimensional column vector
with entries given by:
\begin{equation} \label{eq:L20_PDE_linear_form_q}
	q_i^l = \int_{\Omega_l} \mathbf{q}\cdot\mathbf{f}_i d\Omega,
\end{equation}
and $\mathcal{Q} \equiv (\mathbf{q}_1,\ldots,\mathbf{q}_n)$ constitutes 
the ``library'' of possible terms $\mathbf{q}_i$'s. Note that, the integral
in Eq.~(\ref{eq:L20_PDE_linear_form_q}) involves derivatives of the vector
field $\mathbf{u}$. If the measurements of $\mathbf{u}$ are noisy, evaluating 
these derivatives can be problematic. However, performing integration by parts
entails transferring the derivatives to the weight vector $\mathbf{w}$, which
can be chosen to be smooth. 

For $L \ge n$, an iterative sparse regression algorithm can be used to 
solve~\cite{RGG:2020} the coefficient vector in 
Eq.~(\ref{eq:L20_PDE_linear_form_q}). Each iteration involves the following
form of the solution that minimizes the residual in 
Eq.~(\ref{eq:L20_PDE_linear_form_q}): 
\begin{equation} \label{eq:L20_PDE_c_solution}
\bar{\mathbf{c}} = \mathcal{Q}^{+}\cdot \mathbf{q}_0,
\end{equation}
where $\mathcal{Q}^+$ is the pseudoinverse of the matrix $\mathcal{Q}$.      
Using an empirical thresholding procedure to eliminate dynamically irrelevant
terms $\mathbf{q}_i$, with the solved sparse coefficient vector $\mathbf{c}$, 
one can obtain a minimal PDE model from the measurements $\mathbf{u}$. 
The method was tested~\cite{RGG:2020} with the one-dimensional 
Kuramoto-Sivashinsky model that has in its solutions spatiotemporal chaos, 
the two-dimensional Navier-Stokes equation, and the reaction-diffusion 
equation.  

\section{Discussion} \label{sec:discussion}

The principle of exploiting sparse optimization such as compressive sensing
to find the equations of nonlinear dynamical systems from data was first 
articulated~\cite{WYLKG:2011} in 2011. The basic idea is to expand the 
equations (the velocity field for a continuous time dynamical system or the 
map function for a discrete time system) of the underlying system into a power 
series or a Fourier series of a finite number of terms and then to determine 
the vector of the expansion coefficients based on data through sparse 
optimization. The sparse optimization principle has been demonstrated to 
be effective for finding the governing equations of certain types of nonlinear
dynamical systems for inferring the detailed connection structures of complex 
dynamical networks such as oscillator networks and social networks hosting 
evolutionary game dynamics. In spite of the demonstrated success, limitation
and open questions remain.

A key requirement is that the coefficient vector to be determined
must be sparse. If the vector field or the map function contains a few power
series terms, such as the classical Lorenz~\cite{Lorenz:1963} or 
R\"{o}ssler~\cite{Rossler:1976} chaotic oscillator, or contains a few Fourier 
series terms, such as the standard map~\cite{CI:1973,Chirikov:1979}, then 
sparse optimization can be quite effective and computationally efficient for 
finding the system equations~\cite{WYLKG:2011}. However, if the vector field 
or the map function contains a large number of terms in its power series or 
Fourier series expansion so that the coefficient vector to be determined is 
dense, then the sparse optimization methodology will fail. One such example 
is the classical Ikeda map~\cite{Ikeda:1979,HJM:1985} that describes the 
propagation of a laser pulse in an optical cavity:
\begin{equation} \label{eq:L30_Ikeda}
\mathbf{F}(x,y) = \left( \begin{array}{c}
    a + b (x\cos{\phi} - y\sin{\phi} \\ 
        b (x\sin{\phi} + y\cos{\phi} 
        \end{array} \right),
\end{equation}
with the nonlinear phase variable $\phi$ given by
\begin{equation} \label{eq:L30_Ikeda_phi}
\phi \equiv p - \frac{k}{1 + x^2 + y^2},
\end{equation}
where $a$, $b$, $k$, and $p$ are parameters. It can be seen that both 
components of the map function contain an infinite number of power series
terms, rendering inapplicable sparse optimization for finding the system 
equations from data.

In the mathematical formulation of compressive sensing Eq.~(\ref{eq:L30_CS}),
a requirement is that the projection matrix $\mathcal{G}$ be random, e.g.,
Gaussian type of random matrices with no correlations among the matrix 
elements~\cite{CRT:2006a,CRT:2006b,Donoho:2006,Baraniuk:2007,CW:2008}. However,
in the power series formulation, e.g., Eq.~(\ref{eq:L30_CS_form_x}), 
the elements of the projection matrix are different combinations of the 
powers of the dynamical variables, which are correlated even for a chaotic
system. The demonstrated success in finding the system equations as reviewed
in this article thus has no mathematical guarantee. It may also be possible
that the ``workable'' domain of sparse optimization is larger than that 
guaranteed by rigorous mathematics. To possibly relax the conditions 
under which sparse optimization is effective remains an open mathematical 
issue.

Another difficulty with the application of sparse optimization to find system
equations from data is the need to collect time series from all dynamical
variables of the system. In real world applications, situations are common
where only a limited set of the intrinsic dynamical variables of the system are
externally accessible. The requirement of observing all dynamical variables 
thus represents a formidable obstacle to actual application of the sparse
optimization methods. This should be contrasted to the traditional delay 
coordinate embedding paradigm~\cite{Takens:1981,PCFS:1980} where, in principle,
measurements from a single dynamical variable are sufficient to reconstruct 
the phase space of the underlying system. The capability of the embedding 
paradigm in uncovering the topological and statistical properties of the 
dynamical invariant set responsible for the observed data notwithstanding, it
is unable to yield the system equations.       

In the past two or three years, machine learning has emerged as a promising 
paradigm for predicting the state evolution of nonlinear dynamical systems. 
In particular, a class of recurrent neural 
networks~\cite{Jaeger:2001,MNM:2002,JH:2004,MJ:2013}, the so-called reservoir
computing machines, have attracted considerable attention since 2017 as a 
powerful paradigm for model-free, fully data driven prediction of nonlinear 
and chaotic dynamical systems~\cite{HSRFG:2015,LBMUCJ:2017,PLHGO:2017,LPHGBO:2017,LHO:2018,PWFCHGO:2018,PHGLO:2018,Carroll:2018,NS:2018,ZP:2018,WYGZS:2019,GPG:2019,JL:2019,VPHSGOK:2019,FJZWL:2020,ZJQL:2020}. 
A typical reservoir computing machine consists of an input layer, a hidden
layer that is usually a complex dynamical network, and an output layer. Time 
series data from the dynamical system to be predicted are used to train the 
machine through a series of adjustments to the weights that connect the hidden 
layer with the output layer. Once the machine has been trained, it can predict
the state evolution of the target system for certain duration of time. A well
trained reservoir computing machine can thus be viewed as a ``replica'' of the 
target system, where temporal synchronization between the two can be 
maintained~\cite{FJZWL:2020}. While the machine learning approach does not 
yield the equations of the system, it can be used to predict the system 
behavior especially for those that do not meet the sparsity condition.


\begin{thebibliography}{152}%
\makeatletter
\providecommand \@ifxundefined [1]{%
 \@ifx{#1\undefined}
}%
\providecommand \@ifnum [1]{%
 \ifnum #1\expandafter \@firstoftwo
 \else \expandafter \@secondoftwo
 \fi
}%
\providecommand \@ifx [1]{%
 \ifx #1\expandafter \@firstoftwo
 \else \expandafter \@secondoftwo
 \fi
}%
\providecommand \natexlab [1]{#1}%
\providecommand \enquote  [1]{``#1''}%
\providecommand \bibnamefont  [1]{#1}%
\providecommand \bibfnamefont [1]{#1}%
\providecommand \citenamefont [1]{#1}%
\providecommand \href@noop [0]{\@secondoftwo}%
\providecommand \href [0]{\begingroup \@sanitize@url \@href}%
\providecommand \@href[1]{\@@startlink{#1}\@@href}%
\providecommand \@@href[1]{\endgroup#1\@@endlink}%
\providecommand \@sanitize@url [0]{\catcode `\\12\catcode `\$12\catcode
  `\&12\catcode `\#12\catcode `\^12\catcode `\_12\catcode `\%12\relax}%
\providecommand \@@startlink[1]{}%
\providecommand \@@endlink[0]{}%
\providecommand \url  [0]{\begingroup\@sanitize@url \@url }%
\providecommand \@url [1]{\endgroup\@href {#1}{\urlprefix }}%
\providecommand \urlprefix  [0]{URL }%
\providecommand \Eprint [0]{\href }%
\providecommand \doibase [0]{http://dx.doi.org/}%
\providecommand \selectlanguage [0]{\@gobble}%
\providecommand \bibinfo  [0]{\@secondoftwo}%
\providecommand \bibfield  [0]{\@secondoftwo}%
\providecommand \translation [1]{[#1]}%
\providecommand \BibitemOpen [0]{}%
\providecommand \bibitemStop [0]{}%
\providecommand \bibitemNoStop [0]{.\EOS\space}%
\providecommand \EOS [0]{\spacefactor3000\relax}%
\providecommand \BibitemShut  [1]{\csname bibitem#1\endcsname}%
\let\auto@bib@innerbib\@empty
\bibitem [{\citenamefont {Takens}(1981)}]{Takens:1981}%
  \BibitemOpen
  \bibfield  {author} {\bibinfo {author} {\bibfnamefont {F.}~\bibnamefont
  {Takens}},\ }\bibfield  {title} {\enquote {\bibinfo {title} {Detecting
  strange attractors in fluid turbulence},}\ }in\ \href@noop {} {\emph
  {\bibinfo {booktitle} {Dynamical Systems and Turbulence}}},\ \bibinfo
  {series} {Lecture Notes in Mathematics}, Vol.\ \bibinfo {volume} {898},\
  \bibinfo {editor} {edited by\ \bibinfo {editor} {\bibfnamefont
  {D.}~\bibnamefont {Rand}}\ and\ \bibinfo {editor} {\bibfnamefont {L.~S.}\
  \bibnamefont {Young}}}\ (\bibinfo  {publisher} {Springer-Verlag},\ \bibinfo
  {address} {Berlin},\ \bibinfo {year} {1981})\ pp.\ \bibinfo {pages}
  {366--381}\BibitemShut {NoStop}%
\bibitem [{\citenamefont {Packard}\ \emph {et~al.}(1980)\citenamefont
  {Packard}, \citenamefont {Crutchfield}, \citenamefont {Farmer},\ and\
  \citenamefont {Shaw}}]{PCFS:1980}%
  \BibitemOpen
  \bibfield  {author} {\bibinfo {author} {\bibfnamefont {N.~H.}\ \bibnamefont
  {Packard}}, \bibinfo {author} {\bibfnamefont {J.~P.}\ \bibnamefont
  {Crutchfield}}, \bibinfo {author} {\bibfnamefont {J.~D.}\ \bibnamefont
  {Farmer}}, \ and\ \bibinfo {author} {\bibfnamefont {R.~S.}\ \bibnamefont
  {Shaw}},\ }\bibfield  {title} {\enquote {\bibinfo {title} {Geometry from a
  time series},}\ }\href@noop {} {\bibfield  {journal} {\bibinfo  {journal}
  {Phys. Rev. Lett.}\ }\textbf {\bibinfo {volume} {45}},\ \bibinfo {pages}
  {712} (\bibinfo {year} {1980})}\BibitemShut {NoStop}%
\bibitem [{\citenamefont {Kantz}\ and\ \citenamefont
  {Schreiber}(1997)}]{KS:book}%
  \BibitemOpen
  \bibfield  {author} {\bibinfo {author} {\bibfnamefont {H.}~\bibnamefont
  {Kantz}}\ and\ \bibinfo {author} {\bibfnamefont {T.}~\bibnamefont
  {Schreiber}},\ }\href@noop {} {\emph {\bibinfo {title} {Nonlinear Time Series
  Analysis}}},\ \bibinfo {edition} {1st}\ ed.\ (\bibinfo  {publisher}
  {Cambridge University Press},\ \bibinfo {address} {Cambridge, UK},\ \bibinfo
  {year} {1997})\BibitemShut {NoStop}%
\bibitem [{\citenamefont {Hegger}\ \emph {et~al.}(2007)\citenamefont {Hegger},
  \citenamefont {Kantz},\ and\ \citenamefont {Schreiber}}]{Hegger:book}%
  \BibitemOpen
  \bibfield  {author} {\bibinfo {author} {\bibfnamefont {R.}~\bibnamefont
  {Hegger}}, \bibinfo {author} {\bibfnamefont {H.}~\bibnamefont {Kantz}}, \
  and\ \bibinfo {author} {\bibfnamefont {R.}~\bibnamefont {Schreiber}},\
  }\href@noop {} {\emph {\bibinfo {title} {TISEAN}}},\ \bibinfo {edition}
  {e-book}\ ed.,\ Hegger:book\ (\bibinfo  {publisher}
  {http:$//$www.mpipks-dresden.mpg.de$/$~tisean$/$TISEA$N_{}$3.01$/$
  index.html},\ \bibinfo {address} {Dresden},\ \bibinfo {year}
  {2007})\BibitemShut {NoStop}%
\bibitem [{\citenamefont {Grassberger}\ and\ \citenamefont
  {Procaccia}(1983)}]{GP:1983}%
  \BibitemOpen
  \bibfield  {author} {\bibinfo {author} {\bibfnamefont {P.}~\bibnamefont
  {Grassberger}}\ and\ \bibinfo {author} {\bibfnamefont {I.}~\bibnamefont
  {Procaccia}},\ }\bibfield  {title} {\enquote {\bibinfo {title} {Measuring the
  strangeness of strange attractors},}\ }\href@noop {} {\bibfield  {journal}
  {\bibinfo  {journal} {Physica D}\ }\textbf {\bibinfo {volume} {9}},\ \bibinfo
  {pages} {189} (\bibinfo {year} {1983})}\BibitemShut {NoStop}%
\bibitem [{\citenamefont {Grassberger}(1986)}]{Grass:1986}%
  \BibitemOpen
  \bibfield  {author} {\bibinfo {author} {\bibfnamefont {P.}~\bibnamefont
  {Grassberger}},\ }\bibfield  {title} {\enquote {\bibinfo {title} {Do climatic
  attractors exist?}}\ }\href@noop {} {\bibfield  {journal} {\bibinfo
  {journal} {Nature (London)}\ }\textbf {\bibinfo {volume} {323}},\ \bibinfo
  {pages} {609} (\bibinfo {year} {1986})}\BibitemShut {NoStop}%
\bibitem [{\citenamefont {Procaccia}(1988)}]{Proc:1988}%
  \BibitemOpen
  \bibfield  {author} {\bibinfo {author} {\bibfnamefont {I.}~\bibnamefont
  {Procaccia}},\ }\bibfield  {title} {\enquote {\bibinfo {title} {Complex or
  just complicated?}}\ }\href@noop {} {\bibfield  {journal} {\bibinfo
  {journal} {Nature (London)}\ }\textbf {\bibinfo {volume} {333}},\ \bibinfo
  {pages} {498} (\bibinfo {year} {1988})}\BibitemShut {NoStop}%
\bibitem [{\citenamefont {Osborne}\ and\ \citenamefont
  {Provenzale}(1989)}]{OP:1989}%
  \BibitemOpen
  \bibfield  {author} {\bibinfo {author} {\bibfnamefont {A.}~\bibnamefont
  {Osborne}}\ and\ \bibinfo {author} {\bibfnamefont {A.}~\bibnamefont
  {Provenzale}},\ }\bibfield  {title} {\enquote {\bibinfo {title} {Finite
  correlation dimension for stochastic systems with power-law spectra},}\
  }\href@noop {} {\bibfield  {journal} {\bibinfo  {journal} {Physica D}\
  }\textbf {\bibinfo {volume} {35}},\ \bibinfo {pages} {357} (\bibinfo {year}
  {1989})}\BibitemShut {NoStop}%
\bibitem [{\citenamefont {Lorenz}(1991)}]{Lorenz:1991}%
  \BibitemOpen
  \bibfield  {author} {\bibinfo {author} {\bibfnamefont {E.}~\bibnamefont
  {Lorenz}},\ }\bibfield  {title} {\enquote {\bibinfo {title} {Dimension of
  weather and climate attractors},}\ }\href@noop {} {\bibfield  {journal}
  {\bibinfo  {journal} {Nature (London)}\ }\textbf {\bibinfo {volume} {353}},\
  \bibinfo {pages} {241} (\bibinfo {year} {1991})}\BibitemShut {NoStop}%
\bibitem [{\citenamefont {Ding}\ \emph {et~al.}(1993)\citenamefont {Ding},
  \citenamefont {Grebogi}, \citenamefont {Ott}, \citenamefont {Sauer},\ and\
  \citenamefont {Yorke}}]{DGOSY:1993}%
  \BibitemOpen
  \bibfield  {author} {\bibinfo {author} {\bibfnamefont {M.}~\bibnamefont
  {Ding}}, \bibinfo {author} {\bibfnamefont {C.}~\bibnamefont {Grebogi}},
  \bibinfo {author} {\bibfnamefont {E.}~\bibnamefont {Ott}}, \bibinfo {author}
  {\bibfnamefont {T.~D.}\ \bibnamefont {Sauer}}, \ and\ \bibinfo {author}
  {\bibfnamefont {J.~A.}\ \bibnamefont {Yorke}},\ }\bibfield  {title} {\enquote
  {\bibinfo {title} {Plateau onset for correlation dimension: when does it
  occur?}}\ }\href@noop {} {\bibfield  {journal} {\bibinfo  {journal} {Phys.
  Rev. Lett.}\ }\textbf {\bibinfo {volume} {70}},\ \bibinfo {pages} {3872}
  (\bibinfo {year} {1993})}\BibitemShut {NoStop}%
\bibitem [{\citenamefont {Lai}\ \emph {et~al.}(1996)\citenamefont {Lai},
  \citenamefont {Lerner},\ and\ \citenamefont {Hayden}}]{LLH:1996}%
  \BibitemOpen
  \bibfield  {author} {\bibinfo {author} {\bibfnamefont {Y.-C.}\ \bibnamefont
  {Lai}}, \bibinfo {author} {\bibfnamefont {D.}~\bibnamefont {Lerner}}, \ and\
  \bibinfo {author} {\bibfnamefont {R.}~\bibnamefont {Hayden}},\ }\bibfield
  {title} {\enquote {\bibinfo {title} {An upper bound for the proper delay time
  in chaotic time series analysis},}\ }\href@noop {} {\bibfield  {journal}
  {\bibinfo  {journal} {Phys. Lett. A}\ }\textbf {\bibinfo {volume} {218}},\
  \bibinfo {pages} {30} (\bibinfo {year} {1996})}\BibitemShut {NoStop}%
\bibitem [{\citenamefont {Lai}\ and\ \citenamefont {Lerner}(1998)}]{LL:1998}%
  \BibitemOpen
  \bibfield  {author} {\bibinfo {author} {\bibfnamefont {Y.-C.}\ \bibnamefont
  {Lai}}\ and\ \bibinfo {author} {\bibfnamefont {D.}~\bibnamefont {Lerner}},\
  }\bibfield  {title} {\enquote {\bibinfo {title} {Effective scaling regime for
  computing the correlation dimension in chaotic time series analysis},}\
  }\href@noop {} {\bibfield  {journal} {\bibinfo  {journal} {Physica D}\
  }\textbf {\bibinfo {volume} {115}},\ \bibinfo {pages} {1} (\bibinfo {year}
  {1998})}\BibitemShut {NoStop}%
\bibitem [{\citenamefont {Wolf}\ \emph {et~al.}(1985)\citenamefont {Wolf},
  \citenamefont {Swift}, \citenamefont {Swinney},\ and\ \citenamefont
  {Vastano}}]{WSSV:1985}%
  \BibitemOpen
  \bibfield  {author} {\bibinfo {author} {\bibfnamefont {A.}~\bibnamefont
  {Wolf}}, \bibinfo {author} {\bibfnamefont {J.~B.}\ \bibnamefont {Swift}},
  \bibinfo {author} {\bibfnamefont {H.~L.}\ \bibnamefont {Swinney}}, \ and\
  \bibinfo {author} {\bibfnamefont {J.~A.}\ \bibnamefont {Vastano}},\
  }\bibfield  {title} {\enquote {\bibinfo {title} {{Determining Lyapunov
  exponents from a time series}},}\ }\href@noop {} {\bibfield  {journal}
  {\bibinfo  {journal} {Physica D}\ }\textbf {\bibinfo {volume} {16}},\
  \bibinfo {pages} {285} (\bibinfo {year} {1985})}\BibitemShut {NoStop}%
\bibitem [{\citenamefont {Sano}\ and\ \citenamefont
  {Sawada}(1985)}]{Sano:1985}%
  \BibitemOpen
  \bibfield  {author} {\bibinfo {author} {\bibfnamefont {M.}~\bibnamefont
  {Sano}}\ and\ \bibinfo {author} {\bibfnamefont {Y.}~\bibnamefont {Sawada}},\
  }\bibfield  {title} {\enquote {\bibinfo {title} {{Measurement of the Lyapunov
  spectrum from a chaotic time series}},}\ }\href@noop {} {\bibfield  {journal}
  {\bibinfo  {journal} {Phys. Rev. Lett.}\ }\textbf {\bibinfo {volume} {55}},\
  \bibinfo {pages} {1082} (\bibinfo {year} {1985})}\BibitemShut {NoStop}%
\bibitem [{\citenamefont {Eckmann}\ and\ \citenamefont
  {Ruelle}(1985)}]{ER:1985}%
  \BibitemOpen
  \bibfield  {author} {\bibinfo {author} {\bibfnamefont {J.~P.}\ \bibnamefont
  {Eckmann}}\ and\ \bibinfo {author} {\bibfnamefont {D.}~\bibnamefont
  {Ruelle}},\ }\bibfield  {title} {\enquote {\bibinfo {title} {Ergodic theory
  of chaos and strange attractors},}\ }\href@noop {} {\bibfield  {journal}
  {\bibinfo  {journal} {Rev. Mod. Phys.}\ }\textbf {\bibinfo {volume} {57}},\
  \bibinfo {pages} {617} (\bibinfo {year} {1985})}\BibitemShut {NoStop}%
\bibitem [{\citenamefont {Eckmann}\ \emph {et~al.}(1986)\citenamefont
  {Eckmann}, \citenamefont {Kamphorst}, \citenamefont {Ruelle},\ and\
  \citenamefont {Ciliberto}}]{EKRC:1986}%
  \BibitemOpen
  \bibfield  {author} {\bibinfo {author} {\bibfnamefont {J.~P.}\ \bibnamefont
  {Eckmann}}, \bibinfo {author} {\bibfnamefont {S.~O.}\ \bibnamefont
  {Kamphorst}}, \bibinfo {author} {\bibfnamefont {D.}~\bibnamefont {Ruelle}}, \
  and\ \bibinfo {author} {\bibfnamefont {S.}~\bibnamefont {Ciliberto}},\
  }\bibfield  {title} {\enquote {\bibinfo {title} {Liapunov exponents from time
  series},}\ }\href@noop {} {\bibfield  {journal} {\bibinfo  {journal} {Phys.
  Rev. A}\ }\textbf {\bibinfo {volume} {34}},\ \bibinfo {pages} {4971}
  (\bibinfo {year} {1986})}\BibitemShut {NoStop}%
\bibitem [{\citenamefont {Brown}\ \emph {et~al.}(1991)\citenamefont {Brown},
  \citenamefont {Bryant},\ and\ \citenamefont {Abarbanel}}]{BBA:1991}%
  \BibitemOpen
  \bibfield  {author} {\bibinfo {author} {\bibfnamefont {R.}~\bibnamefont
  {Brown}}, \bibinfo {author} {\bibfnamefont {P.}~\bibnamefont {Bryant}}, \
  and\ \bibinfo {author} {\bibfnamefont {H.~D.~I.}\ \bibnamefont {Abarbanel}},\
  }\bibfield  {title} {\enquote {\bibinfo {title} {{Computing the Lyapunov
  spectrum of a dynamical system from an observed time series}},}\ }\href@noop
  {} {\bibfield  {journal} {\bibinfo  {journal} {Phys. Rev. A}\ }\textbf
  {\bibinfo {volume} {43}},\ \bibinfo {pages} {2787} (\bibinfo {year}
  {1991})}\BibitemShut {NoStop}%
\bibitem [{\citenamefont {Sauer}\ \emph {et~al.}(1998)\citenamefont {Sauer},
  \citenamefont {Tempkin},\ and\ \citenamefont {Yorke}}]{STY:1998}%
  \BibitemOpen
  \bibfield  {author} {\bibinfo {author} {\bibfnamefont {T.~D.}\ \bibnamefont
  {Sauer}}, \bibinfo {author} {\bibfnamefont {J.~A.}\ \bibnamefont {Tempkin}},
  \ and\ \bibinfo {author} {\bibfnamefont {J.~A.}\ \bibnamefont {Yorke}},\
  }\bibfield  {title} {\enquote {\bibinfo {title} {{Spurious Lyapunov exponents
  in attractor reconstruction}},}\ }\href@noop {} {\bibfield  {journal}
  {\bibinfo  {journal} {Phys. Rev. Lett.}\ }\textbf {\bibinfo {volume} {81}},\
  \bibinfo {pages} {4341} (\bibinfo {year} {1998})}\BibitemShut {NoStop}%
\bibitem [{\citenamefont {Sauer}\ and\ \citenamefont {Yorke}(1999)}]{SY:1999}%
  \BibitemOpen
  \bibfield  {author} {\bibinfo {author} {\bibfnamefont {T.~D.}\ \bibnamefont
  {Sauer}}\ and\ \bibinfo {author} {\bibfnamefont {J.~A.}\ \bibnamefont
  {Yorke}},\ }\bibfield  {title} {\enquote {\bibinfo {title} {{Reconstructing
  the Jacobian from data with observational noise}},}\ }\href@noop {}
  {\bibfield  {journal} {\bibinfo  {journal} {Phys. Rev. Lett.}\ }\textbf
  {\bibinfo {volume} {83}},\ \bibinfo {pages} {1331} (\bibinfo {year}
  {1999})}\BibitemShut {NoStop}%
\bibitem [{\citenamefont {Lathrop}\ and\ \citenamefont
  {Kostelich}(1989)}]{LK:1989}%
  \BibitemOpen
  \bibfield  {author} {\bibinfo {author} {\bibfnamefont {D.~P.}\ \bibnamefont
  {Lathrop}}\ and\ \bibinfo {author} {\bibfnamefont {E.~J.}\ \bibnamefont
  {Kostelich}},\ }\bibfield  {title} {\enquote {\bibinfo {title}
  {Characterization of an experimental strange attractor by periodic orbits},}\
  }\href@noop {} {\bibfield  {journal} {\bibinfo  {journal} {Phys. Rev. A}\
  }\textbf {\bibinfo {volume} {40}},\ \bibinfo {pages} {4028} (\bibinfo {year}
  {1989})}\BibitemShut {NoStop}%
\bibitem [{\citenamefont {Badii}\ \emph {et~al.}(1994)\citenamefont {Badii},
  \citenamefont {Brun}, \citenamefont {Finardi}, \citenamefont {Flepp},
  \citenamefont {Holzner}, \citenamefont {Pariso}, \citenamefont {Reyl},\ and\
  \citenamefont {Simonet}}]{BBFFHPRS:1994}%
  \BibitemOpen
  \bibfield  {author} {\bibinfo {author} {\bibfnamefont {R.}~\bibnamefont
  {Badii}}, \bibinfo {author} {\bibfnamefont {E.}~\bibnamefont {Brun}},
  \bibinfo {author} {\bibfnamefont {M.}~\bibnamefont {Finardi}}, \bibinfo
  {author} {\bibfnamefont {L.}~\bibnamefont {Flepp}}, \bibinfo {author}
  {\bibfnamefont {R.}~\bibnamefont {Holzner}}, \bibinfo {author} {\bibfnamefont
  {J.}~\bibnamefont {Pariso}}, \bibinfo {author} {\bibfnamefont
  {C.}~\bibnamefont {Reyl}}, \ and\ \bibinfo {author} {\bibfnamefont
  {J.}~\bibnamefont {Simonet}},\ }\bibfield  {title} {\enquote {\bibinfo
  {title} {Progress in the analysis of experimental chaos through periodic
  orbits},}\ }\href@noop {} {\bibfield  {journal} {\bibinfo  {journal} {Rev.
  Mod. Phys.}\ }\textbf {\bibinfo {volume} {66}},\ \bibinfo {pages} {1389}
  (\bibinfo {year} {1994})}\BibitemShut {NoStop}%
\bibitem [{\citenamefont {Pierson}\ and\ \citenamefont {Moss}(1995)}]{PM:1995}%
  \BibitemOpen
  \bibfield  {author} {\bibinfo {author} {\bibfnamefont {D.}~\bibnamefont
  {Pierson}}\ and\ \bibinfo {author} {\bibfnamefont {F.}~\bibnamefont {Moss}},\
  }\bibfield  {title} {\enquote {\bibinfo {title} {Detecting periodic unstable
  points in noisy chaotic and limit-cycle attractors with applications to
  biology},}\ }\href@noop {} {\bibfield  {journal} {\bibinfo  {journal} {Phys.
  Rev. Lett.}\ }\textbf {\bibinfo {volume} {75}},\ \bibinfo {pages} {2124}
  (\bibinfo {year} {1995})}\BibitemShut {NoStop}%
\bibitem [{\citenamefont {Pei}\ and\ \citenamefont {Moss}(1996)}]{PM:1996a}%
  \BibitemOpen
  \bibfield  {author} {\bibinfo {author} {\bibfnamefont {X.}~\bibnamefont
  {Pei}}\ and\ \bibinfo {author} {\bibfnamefont {F.}~\bibnamefont {Moss}},\
  }\bibfield  {title} {\enquote {\bibinfo {title} {Characterization of
  low-dimensional dynamics in the crayfish caudal photoreceptor},}\ }\href@noop
  {} {\bibfield  {journal} {\bibinfo  {journal} {Nature (London)}\ }\textbf
  {\bibinfo {volume} {379}},\ \bibinfo {pages} {618} (\bibinfo {year}
  {1996})}\BibitemShut {NoStop}%
\bibitem [{\citenamefont {So}\ \emph {et~al.}(1996)\citenamefont {So},
  \citenamefont {Ott}, \citenamefont {Schiff}, \citenamefont {Kaplan},
  \citenamefont {Sauer},\ and\ \citenamefont {Grebogi}}]{SOSKSG:1996}%
  \BibitemOpen
  \bibfield  {author} {\bibinfo {author} {\bibfnamefont {P.}~\bibnamefont
  {So}}, \bibinfo {author} {\bibfnamefont {E.}~\bibnamefont {Ott}}, \bibinfo
  {author} {\bibfnamefont {S.~J.}\ \bibnamefont {Schiff}}, \bibinfo {author}
  {\bibfnamefont {D.~T.}\ \bibnamefont {Kaplan}}, \bibinfo {author}
  {\bibfnamefont {T.~D.}\ \bibnamefont {Sauer}}, \ and\ \bibinfo {author}
  {\bibfnamefont {C.}~\bibnamefont {Grebogi}},\ }\bibfield  {title} {\enquote
  {\bibinfo {title} {Detecting unstable periodic orbits in chaotic experimental
  data},}\ }\href@noop {} {\bibfield  {journal} {\bibinfo  {journal} {Phys.
  Rev. Lett.}\ }\textbf {\bibinfo {volume} {76}},\ \bibinfo {pages} {4705}
  (\bibinfo {year} {1996})}\BibitemShut {NoStop}%
\bibitem [{\citenamefont {Allie}\ and\ \citenamefont {Mees}(1997)}]{AM:1997}%
  \BibitemOpen
  \bibfield  {author} {\bibinfo {author} {\bibfnamefont {S.}~\bibnamefont
  {Allie}}\ and\ \bibinfo {author} {\bibfnamefont {A.}~\bibnamefont {Mees}},\
  }\bibfield  {title} {\enquote {\bibinfo {title} {Finding periodic points from
  short time series},}\ }\href@noop {} {\bibfield  {journal} {\bibinfo
  {journal} {Phys. Rev. E}\ }\textbf {\bibinfo {volume} {56}},\ \bibinfo
  {pages} {346} (\bibinfo {year} {1997})}\BibitemShut {NoStop}%
\bibitem [{\citenamefont {Pecora}\ \emph {et~al.}(1995)\citenamefont {Pecora},
  \citenamefont {Carroll},\ and\ \citenamefont {Heagy}}]{PCH:1995}%
  \BibitemOpen
  \bibfield  {author} {\bibinfo {author} {\bibfnamefont {L.~M.}\ \bibnamefont
  {Pecora}}, \bibinfo {author} {\bibfnamefont {T.~L.}\ \bibnamefont {Carroll}},
  \ and\ \bibinfo {author} {\bibfnamefont {J.~F.}\ \bibnamefont {Heagy}},\
  }\bibfield  {title} {\enquote {\bibinfo {title} {Statistics for mathematical
  properties of maps between time series embeddings},}\ }\href@noop {}
  {\bibfield  {journal} {\bibinfo  {journal} {Phys. Rev. E}\ }\textbf {\bibinfo
  {volume} {52}},\ \bibinfo {pages} {3420} (\bibinfo {year}
  {1995})}\BibitemShut {NoStop}%
\bibitem [{\citenamefont {Pecora}\ and\ \citenamefont
  {Carroll}(1996)}]{PC:1996}%
  \BibitemOpen
  \bibfield  {author} {\bibinfo {author} {\bibfnamefont {L.~M.}\ \bibnamefont
  {Pecora}}\ and\ \bibinfo {author} {\bibfnamefont {T.~L.}\ \bibnamefont
  {Carroll}},\ }\bibfield  {title} {\enquote {\bibinfo {title} {Discontinuous
  and nondifferentiable functions and dimension increase induced by filtering
  chaotic data},}\ }\href@noop {} {\bibfield  {journal} {\bibinfo  {journal}
  {Chaos}\ }\textbf {\bibinfo {volume} {6}},\ \bibinfo {pages} {432} (\bibinfo
  {year} {1996})}\BibitemShut {NoStop}%
\bibitem [{\citenamefont {Pecora}\ \emph {et~al.}(1997)\citenamefont {Pecora},
  \citenamefont {Carroll},\ and\ \citenamefont {Heagy}}]{PCH:1997}%
  \BibitemOpen
  \bibfield  {author} {\bibinfo {author} {\bibfnamefont {L.~M.}\ \bibnamefont
  {Pecora}}, \bibinfo {author} {\bibfnamefont {T.~L.}\ \bibnamefont {Carroll}},
  \ and\ \bibinfo {author} {\bibfnamefont {J.~F.}\ \bibnamefont {Heagy}},\
  }\bibfield  {title} {\enquote {\bibinfo {title} {Statistics for continuity
  and differentiability: an application to attractor reconstruction from time
  series},}\ }\href@noop {} {\bibfield  {journal} {\bibinfo  {journal} {Fields
  Inst. Commun.}\ }\textbf {\bibinfo {volume} {11}},\ \bibinfo {pages} {49}
  (\bibinfo {year} {1997})}\BibitemShut {NoStop}%
\bibitem [{\citenamefont {Goodridge}\ \emph {et~al.}(2001)\citenamefont
  {Goodridge}, \citenamefont {Pecora}, \citenamefont {Carroll},\ and\
  \citenamefont {Rachford}}]{GPCR:2001}%
  \BibitemOpen
  \bibfield  {author} {\bibinfo {author} {\bibfnamefont {C.~L.}\ \bibnamefont
  {Goodridge}}, \bibinfo {author} {\bibfnamefont {L.~M.}\ \bibnamefont
  {Pecora}}, \bibinfo {author} {\bibfnamefont {T.~L.}\ \bibnamefont {Carroll}},
  \ and\ \bibinfo {author} {\bibfnamefont {F.~J.}\ \bibnamefont {Rachford}},\
  }\bibfield  {title} {\enquote {\bibinfo {title} {Detecting functional
  relationships between simultaneous time series},}\ }\href@noop {} {\bibfield
  {journal} {\bibinfo  {journal} {Phys. Rev. E}\ }\textbf {\bibinfo {volume}
  {64}},\ \bibinfo {pages} {026221} (\bibinfo {year} {2001})}\BibitemShut
  {NoStop}%
\bibitem [{\citenamefont {Pecora}\ \emph {et~al.}(2007)\citenamefont {Pecora},
  \citenamefont {Moniz}, \citenamefont {Nichols},\ and\ \citenamefont
  {Carroll}}]{PMNC:2007}%
  \BibitemOpen
  \bibfield  {author} {\bibinfo {author} {\bibfnamefont {L.~M.}\ \bibnamefont
  {Pecora}}, \bibinfo {author} {\bibfnamefont {L.}~\bibnamefont {Moniz}},
  \bibinfo {author} {\bibfnamefont {J.}~\bibnamefont {Nichols}}, \ and\
  \bibinfo {author} {\bibfnamefont {T.}~\bibnamefont {Carroll}},\ }\bibfield
  {title} {\enquote {\bibinfo {title} {A unified approach to attractor
  reconstruction},}\ }\href@noop {} {\bibfield  {journal} {\bibinfo  {journal}
  {Chaos}\ }\textbf {\bibinfo {volume} {17}},\ \bibinfo {pages} {013110}
  (\bibinfo {year} {2007})}\BibitemShut {NoStop}%
\bibitem [{\citenamefont {Theiler}(1986)}]{Theiler:1986}%
  \BibitemOpen
  \bibfield  {author} {\bibinfo {author} {\bibfnamefont {J.}~\bibnamefont
  {Theiler}},\ }\bibfield  {title} {\enquote {\bibinfo {title} {Spurious
  dimension from correlation algorithms applied to limited time series data},}\
  }\href@noop {} {\bibfield  {journal} {\bibinfo  {journal} {Phys. Rev. A}\
  }\textbf {\bibinfo {volume} {34}},\ \bibinfo {pages} {2427} (\bibinfo {year}
  {1986})}\BibitemShut {NoStop}%
\bibitem [{\citenamefont {Liebert}\ and\ \citenamefont
  {Schuster}(1989)}]{LS:1989}%
  \BibitemOpen
  \bibfield  {author} {\bibinfo {author} {\bibfnamefont {W.}~\bibnamefont
  {Liebert}}\ and\ \bibinfo {author} {\bibfnamefont {H.~G.}\ \bibnamefont
  {Schuster}},\ }\bibfield  {title} {\enquote {\bibinfo {title} {Proper choice
  of the time-delay for the analysis of chaotic time-series},}\ }\href@noop {}
  {\bibfield  {journal} {\bibinfo  {journal} {Phys. Lett. A}\ }\textbf
  {\bibinfo {volume} {142}},\ \bibinfo {pages} {107} (\bibinfo {year}
  {1989})}\BibitemShut {NoStop}%
\bibitem [{\citenamefont {Liebert}\ \emph {et~al.}(1991)\citenamefont
  {Liebert}, \citenamefont {Pawelzik},\ and\ \citenamefont
  {Schuster}}]{LPS:1991}%
  \BibitemOpen
  \bibfield  {author} {\bibinfo {author} {\bibfnamefont {W.}~\bibnamefont
  {Liebert}}, \bibinfo {author} {\bibfnamefont {K.}~\bibnamefont {Pawelzik}}, \
  and\ \bibinfo {author} {\bibfnamefont {H.~G.}\ \bibnamefont {Schuster}},\
  }\bibfield  {title} {\enquote {\bibinfo {title} {Optimal embeddings of
  chaotic attractors from topological considerations},}\ }\href@noop {}
  {\bibfield  {journal} {\bibinfo  {journal} {Europhys. Lett.}\ }\textbf
  {\bibinfo {volume} {14}},\ \bibinfo {pages} {521} (\bibinfo {year}
  {1991})}\BibitemShut {NoStop}%
\bibitem [{\citenamefont {Buzug}\ and\ \citenamefont
  {Pfister}(1992)}]{BP:1992}%
  \BibitemOpen
  \bibfield  {author} {\bibinfo {author} {\bibfnamefont {T.}~\bibnamefont
  {Buzug}}\ and\ \bibinfo {author} {\bibfnamefont {G.}~\bibnamefont
  {Pfister}},\ }\bibfield  {title} {\enquote {\bibinfo {title} {Optimal delay
  time and embedding dimension for delay-time coordinates by analysis of the
  glocal static and local dynamic behavior of strange attractors},}\
  }\href@noop {} {\bibfield  {journal} {\bibinfo  {journal} {Phys. Rev. A}\
  }\textbf {\bibinfo {volume} {45}},\ \bibinfo {pages} {7073} (\bibinfo {year}
  {1992})}\BibitemShut {NoStop}%
\bibitem [{\citenamefont {Kember}\ and\ \citenamefont
  {Fowler}(1993)}]{KF:1993}%
  \BibitemOpen
  \bibfield  {author} {\bibinfo {author} {\bibfnamefont {G.}~\bibnamefont
  {Kember}}\ and\ \bibinfo {author} {\bibfnamefont {A.~C.}\ \bibnamefont
  {Fowler}},\ }\bibfield  {title} {\enquote {\bibinfo {title} {A
  correlation-function for choosing time delays in-phase portrait
  reconstructions},}\ }\href@noop {} {\bibfield  {journal} {\bibinfo  {journal}
  {Phys. Lett. A}\ }\textbf {\bibinfo {volume} {179}},\ \bibinfo {pages} {72}
  (\bibinfo {year} {1993})}\BibitemShut {NoStop}%
\bibitem [{\citenamefont {Rosenstein}\ \emph {et~al.}(1994)\citenamefont
  {Rosenstein}, \citenamefont {Collins},\ and\ \citenamefont
  {Luca}}]{RCdL:1994}%
  \BibitemOpen
  \bibfield  {author} {\bibinfo {author} {\bibfnamefont {M.~T.}\ \bibnamefont
  {Rosenstein}}, \bibinfo {author} {\bibfnamefont {J.~J.}\ \bibnamefont
  {Collins}}, \ and\ \bibinfo {author} {\bibfnamefont {C.~J.~D.}\ \bibnamefont
  {Luca}},\ }\bibfield  {title} {\enquote {\bibinfo {title} {Reconstruction
  expansion as a geometry-based framework for choosing proper delay times},}\
  }\href@noop {} {\bibfield  {journal} {\bibinfo  {journal} {Physica D}\
  }\textbf {\bibinfo {volume} {73}},\ \bibinfo {pages} {82} (\bibinfo {year}
  {1994})}\BibitemShut {NoStop}%
\bibitem [{\citenamefont {Sauer}\ \emph {et~al.}(1991)\citenamefont {Sauer},
  \citenamefont {Yorke},\ and\ \citenamefont {Casdagli}}]{SYC:1991}%
  \BibitemOpen
  \bibfield  {author} {\bibinfo {author} {\bibfnamefont {T.~D.}\ \bibnamefont
  {Sauer}}, \bibinfo {author} {\bibfnamefont {J.~A.}\ \bibnamefont {Yorke}}, \
  and\ \bibinfo {author} {\bibfnamefont {M.}~\bibnamefont {Casdagli}},\
  }\bibfield  {title} {\enquote {\bibinfo {title} {Embedology},}\ }\href@noop
  {} {\bibfield  {journal} {\bibinfo  {journal} {J. Stat. Phys.}\ }\textbf
  {\bibinfo {volume} {65}},\ \bibinfo {pages} {579} (\bibinfo {year}
  {1991})}\BibitemShut {NoStop}%
\bibitem [{\citenamefont {Lai}\ and\ \citenamefont {T\'{e}l}(2011)}]{LT:book}%
  \BibitemOpen
  \bibfield  {author} {\bibinfo {author} {\bibfnamefont {Y.-C.}\ \bibnamefont
  {Lai}}\ and\ \bibinfo {author} {\bibfnamefont {T.}~\bibnamefont {T\'{e}l}},\
  }\href@noop {} {\emph {\bibinfo {title} {Transient Chaos - Complex Dynamics
  on Finite Time Scales}}},\ \bibinfo {edition} {1st}\ ed.\ (\bibinfo
  {publisher} {Springer},\ \bibinfo {address} {New York},\ \bibinfo {year}
  {2011})\BibitemShut {NoStop}%
\bibitem [{\citenamefont {J\'anosi}\ \emph {et~al.}(1994)\citenamefont
  {J\'anosi}, \citenamefont {Flepp},\ and\ \citenamefont {T\'{e}l}}]{JFT:1994}%
  \BibitemOpen
  \bibfield  {author} {\bibinfo {author} {\bibfnamefont {J.}~\bibnamefont
  {J\'anosi}}, \bibinfo {author} {\bibfnamefont {L.}~\bibnamefont {Flepp}}, \
  and\ \bibinfo {author} {\bibfnamefont {T.}~\bibnamefont {T\'{e}l}},\
  }\bibfield  {title} {\enquote {\bibinfo {title} {{Exploring transient chaos
  in an NMR-laser experiment}},}\ }\href@noop {} {\bibfield  {journal}
  {\bibinfo  {journal} {Phys. Rev. Lett.}\ }\textbf {\bibinfo {volume} {73}},\
  \bibinfo {pages} {529} (\bibinfo {year} {1994})}\BibitemShut {NoStop}%
\bibitem [{\citenamefont {J\'anosi}\ and\ \citenamefont
  {T\'{e}l}(1994)}]{JT:1994}%
  \BibitemOpen
  \bibfield  {author} {\bibinfo {author} {\bibfnamefont {J.}~\bibnamefont
  {J\'anosi}}\ and\ \bibinfo {author} {\bibfnamefont {T.}~\bibnamefont
  {T\'{e}l}},\ }\bibfield  {title} {\enquote {\bibinfo {title} {Time-series
  analysis of transient chaos},}\ }\href@noop {} {\bibfield  {journal}
  {\bibinfo  {journal} {Phys. Rev. E}\ }\textbf {\bibinfo {volume} {49}},\
  \bibinfo {pages} {2756} (\bibinfo {year} {1994})}\BibitemShut {NoStop}%
\bibitem [{\citenamefont {Dhamala}\ \emph {et~al.}(2000)\citenamefont
  {Dhamala}, \citenamefont {Lai},\ and\ \citenamefont {Kostelich}}]{DLK:2000}%
  \BibitemOpen
  \bibfield  {author} {\bibinfo {author} {\bibfnamefont {M.}~\bibnamefont
  {Dhamala}}, \bibinfo {author} {\bibfnamefont {Y.-C.}\ \bibnamefont {Lai}}, \
  and\ \bibinfo {author} {\bibfnamefont {E.~J.}\ \bibnamefont {Kostelich}},\
  }\bibfield  {title} {\enquote {\bibinfo {title} {Detecting unstable periodic
  orbits from transient chaotic time series},}\ }\href@noop {} {\bibfield
  {journal} {\bibinfo  {journal} {Phys. Rev. E}\ }\textbf {\bibinfo {volume}
  {61}},\ \bibinfo {pages} {6485} (\bibinfo {year} {2000})}\BibitemShut
  {NoStop}%
\bibitem [{\citenamefont {Dhamala}\ \emph {et~al.}(2001)\citenamefont
  {Dhamala}, \citenamefont {Lai},\ and\ \citenamefont {Kostelich}}]{DLK:2001}%
  \BibitemOpen
  \bibfield  {author} {\bibinfo {author} {\bibfnamefont {M.}~\bibnamefont
  {Dhamala}}, \bibinfo {author} {\bibfnamefont {Y.-C.}\ \bibnamefont {Lai}}, \
  and\ \bibinfo {author} {\bibfnamefont {E.~J.}\ \bibnamefont {Kostelich}},\
  }\bibfield  {title} {\enquote {\bibinfo {title} {Analysis of transient
  chaotic time series},}\ }\href@noop {} {\bibfield  {journal} {\bibinfo
  {journal} {Phys. Rev. E}\ }\textbf {\bibinfo {volume} {64}},\ \bibinfo
  {pages} {056207} (\bibinfo {year} {2001})}\BibitemShut {NoStop}%
\bibitem [{\citenamefont {Triandaf}\ \emph {et~al.}(2003)\citenamefont
  {Triandaf}, \citenamefont {Bollt},\ and\ \citenamefont
  {Schwartz}}]{TBS:2003}%
  \BibitemOpen
  \bibfield  {author} {\bibinfo {author} {\bibfnamefont {I.}~\bibnamefont
  {Triandaf}}, \bibinfo {author} {\bibfnamefont {E.}~\bibnamefont {Bollt}}, \
  and\ \bibinfo {author} {\bibfnamefont {I.~B.}\ \bibnamefont {Schwartz}},\
  }\bibfield  {title} {\enquote {\bibinfo {title} {Approximating stable and
  unstable manifolds in experiments},}\ }\href@noop {} {\bibfield  {journal}
  {\bibinfo  {journal} {Phys. Rev. E}\ }\textbf {\bibinfo {volume} {67}},\
  \bibinfo {pages} {037201} (\bibinfo {year} {2003})}\BibitemShut {NoStop}%
\bibitem [{\citenamefont {Taylor}\ and\ \citenamefont
  {Campbell}(2007)}]{TC:2007}%
  \BibitemOpen
  \bibfield  {author} {\bibinfo {author} {\bibfnamefont {S.~R.}\ \bibnamefont
  {Taylor}}\ and\ \bibinfo {author} {\bibfnamefont {S.~A.}\ \bibnamefont
  {Campbell}},\ }\bibfield  {title} {\enquote {\bibinfo {title} {Approximating
  chaotic saddles in delay differential equations},}\ }\href@noop {} {\bibfield
   {journal} {\bibinfo  {journal} {Phys. Rev. E}\ }\textbf {\bibinfo {volume}
  {75}},\ \bibinfo {pages} {046215} (\bibinfo {year} {2007})}\BibitemShut
  {NoStop}%
\bibitem [{\citenamefont {Farmer}\ and\ \citenamefont
  {Sidorowich}(1987)}]{FS:1987}%
  \BibitemOpen
  \bibfield  {author} {\bibinfo {author} {\bibfnamefont {J.~D.}\ \bibnamefont
  {Farmer}}\ and\ \bibinfo {author} {\bibfnamefont {J.~J.}\ \bibnamefont
  {Sidorowich}},\ }\bibfield  {title} {\enquote {\bibinfo {title} {Predicting
  chaotic time series},}\ }\href {\doibase 10.1103/PhysRevLett.59.845}
  {\bibfield  {journal} {\bibinfo  {journal} {Phys. Rev. Lett.}\ }\textbf
  {\bibinfo {volume} {59}},\ \bibinfo {pages} {845} (\bibinfo {year}
  {1987})}\BibitemShut {NoStop}%
\bibitem [{\citenamefont {Casdagli}(1989)}]{Casdagli:1989}%
  \BibitemOpen
  \bibfield  {author} {\bibinfo {author} {\bibfnamefont {M.}~\bibnamefont
  {Casdagli}},\ }\bibfield  {title} {\enquote {\bibinfo {title} {Nonlinear
  prediction of chaotic time series},}\ }\href@noop {} {\bibfield  {journal}
  {\bibinfo  {journal} {Physica D}\ }\textbf {\bibinfo {volume} {35}},\
  \bibinfo {pages} {335} (\bibinfo {year} {1989})}\BibitemShut {NoStop}%
\bibitem [{\citenamefont {Sugihara}\ \emph {et~al.}(1990)\citenamefont
  {Sugihara}, \citenamefont {Grenfell}, \citenamefont {May}, \citenamefont
  {Chesson}, \citenamefont {Platt},\ and\ \citenamefont
  {Williamson}}]{SGMCPW:1990}%
  \BibitemOpen
  \bibfield  {author} {\bibinfo {author} {\bibfnamefont {G.}~\bibnamefont
  {Sugihara}}, \bibinfo {author} {\bibfnamefont {B.}~\bibnamefont {Grenfell}},
  \bibinfo {author} {\bibfnamefont {R.~M.}\ \bibnamefont {May}}, \bibinfo
  {author} {\bibfnamefont {P.}~\bibnamefont {Chesson}}, \bibinfo {author}
  {\bibfnamefont {H.~M.}\ \bibnamefont {Platt}}, \ and\ \bibinfo {author}
  {\bibfnamefont {M.}~\bibnamefont {Williamson}},\ }\bibfield  {title}
  {\enquote {\bibinfo {title} {Distinguishing error from chaos in ecological
  time series},}\ }\href@noop {} {\bibfield  {journal} {\bibinfo  {journal}
  {Phil. Trans. Roy. Soc. London B}\ }\textbf {\bibinfo {volume} {330}},\
  \bibinfo {pages} {235} (\bibinfo {year} {1990})}\BibitemShut {NoStop}%
\bibitem [{\citenamefont {Kurths}\ and\ \citenamefont
  {Ruzmaikin}(1990)}]{KR:1990}%
  \BibitemOpen
  \bibfield  {author} {\bibinfo {author} {\bibfnamefont {J.}~\bibnamefont
  {Kurths}}\ and\ \bibinfo {author} {\bibfnamefont {A.~A.}\ \bibnamefont
  {Ruzmaikin}},\ }\bibfield  {title} {\enquote {\bibinfo {title} {On
  forecasting the sunspot numbers},}\ }\href@noop {} {\bibfield  {journal}
  {\bibinfo  {journal} {Solar Phys.}\ }\textbf {\bibinfo {volume} {126}},\
  \bibinfo {pages} {407} (\bibinfo {year} {1990})}\BibitemShut {NoStop}%
\bibitem [{\citenamefont {Grassberger}\ and\ \citenamefont
  {Schreiber}(1990)}]{GS:1990}%
  \BibitemOpen
  \bibfield  {author} {\bibinfo {author} {\bibfnamefont {P.}~\bibnamefont
  {Grassberger}}\ and\ \bibinfo {author} {\bibfnamefont {T.}~\bibnamefont
  {Schreiber}},\ }\bibfield  {title} {\enquote {\bibinfo {title} {Nonlinear
  time sequence analysis},}\ }\href@noop {} {\bibfield  {journal} {\bibinfo
  {journal} {Int. J. Bif. Chaos}\ }\textbf {\bibinfo {volume} {1}},\ \bibinfo
  {pages} {521} (\bibinfo {year} {1990})}\BibitemShut {NoStop}%
\bibitem [{\citenamefont {Gouesbet}(1991)}]{Gouesbet:1991}%
  \BibitemOpen
  \bibfield  {author} {\bibinfo {author} {\bibfnamefont {G.}~\bibnamefont
  {Gouesbet}},\ }\bibfield  {title} {\enquote {\bibinfo {title}
  {{Reconstruction of standard and inverse vector fields equivalent to a
  R\"ossler system}},}\ }\href {\doibase 10.1103/PhysRevA.44.6264} {\bibfield
  {journal} {\bibinfo  {journal} {Phys. Rev. A}\ }\textbf {\bibinfo {volume}
  {44}},\ \bibinfo {pages} {6264} (\bibinfo {year} {1991})}\BibitemShut
  {NoStop}%
\bibitem [{\citenamefont {Tsonis}\ and\ \citenamefont
  {Elsner}(1992)}]{TE:1992}%
  \BibitemOpen
  \bibfield  {author} {\bibinfo {author} {\bibfnamefont {A.~A.}\ \bibnamefont
  {Tsonis}}\ and\ \bibinfo {author} {\bibfnamefont {J.~B.}\ \bibnamefont
  {Elsner}},\ }\bibfield  {title} {\enquote {\bibinfo {title} {Nonlinear
  prediction as a way of distinguishing chaos from random fractal sequences},}\
  }\href@noop {} {\bibfield  {journal} {\bibinfo  {journal} {Nature (London)}\
  }\textbf {\bibinfo {volume} {358}},\ \bibinfo {pages} {217} (\bibinfo {year}
  {1992})}\BibitemShut {NoStop}%
\bibitem [{\citenamefont {Baake}\ \emph {et~al.}(1992)\citenamefont {Baake},
  \citenamefont {Baake}, \citenamefont {Bock},\ and\ \citenamefont
  {Briggs}}]{BBBB:1992}%
  \BibitemOpen
  \bibfield  {author} {\bibinfo {author} {\bibfnamefont {E.}~\bibnamefont
  {Baake}}, \bibinfo {author} {\bibfnamefont {M.}~\bibnamefont {Baake}},
  \bibinfo {author} {\bibfnamefont {H.~G.}\ \bibnamefont {Bock}}, \ and\
  \bibinfo {author} {\bibfnamefont {K.~M.}\ \bibnamefont {Briggs}},\ }\bibfield
   {title} {\enquote {\bibinfo {title} {Fitting ordinary differential equations
  to chaotic data},}\ }\href {\doibase 10.1103/PhysRevA.45.5524} {\bibfield
  {journal} {\bibinfo  {journal} {Phys. Rev. A}\ }\textbf {\bibinfo {volume}
  {45}},\ \bibinfo {pages} {5524} (\bibinfo {year} {1992})}\BibitemShut
  {NoStop}%
\bibitem [{\citenamefont {Longtin}(1993)}]{Longtin:1993}%
  \BibitemOpen
  \bibfield  {author} {\bibinfo {author} {\bibfnamefont {A.}~\bibnamefont
  {Longtin}},\ }\bibfield  {title} {\enquote {\bibinfo {title} {Nonlinear
  forecasting of spike trains from sensory neurons},}\ }\href@noop {}
  {\bibfield  {journal} {\bibinfo  {journal} {Int. J. Bif. Chaos}\ }\textbf
  {\bibinfo {volume} {3}},\ \bibinfo {pages} {651} (\bibinfo {year}
  {1993})}\BibitemShut {NoStop}%
\bibitem [{\citenamefont {Murray}(1993)}]{Murray:1993}%
  \BibitemOpen
  \bibfield  {author} {\bibinfo {author} {\bibfnamefont {D.~B.}\ \bibnamefont
  {Murray}},\ }\bibfield  {title} {\enquote {\bibinfo {title} {Forecasting a
  chaotic time series using an improved metric for embedding space},}\
  }\href@noop {} {\bibfield  {journal} {\bibinfo  {journal} {Physica D}\
  }\textbf {\bibinfo {volume} {68}},\ \bibinfo {pages} {318} (\bibinfo {year}
  {1993})}\BibitemShut {NoStop}%
\bibitem [{\citenamefont {Sauer}(1994)}]{Sauer:1994}%
  \BibitemOpen
  \bibfield  {author} {\bibinfo {author} {\bibfnamefont {T.}~\bibnamefont
  {Sauer}},\ }\bibfield  {title} {\enquote {\bibinfo {title} {Reconstruction of
  dynamical systems from interspike intervals},}\ }\href {\doibase
  10.1103/PhysRevLett.72.3811} {\bibfield  {journal} {\bibinfo  {journal}
  {Phys. Rev. Lett.}\ }\textbf {\bibinfo {volume} {72}},\ \bibinfo {pages}
  {3811} (\bibinfo {year} {1994})}\BibitemShut {NoStop}%
\bibitem [{\citenamefont {Sugihara}(1994)}]{Sugihara:1994}%
  \BibitemOpen
  \bibfield  {author} {\bibinfo {author} {\bibfnamefont {G.}~\bibnamefont
  {Sugihara}},\ }\bibfield  {title} {\enquote {\bibinfo {title} {Nonlinear
  forecasting for the classification of natural time series},}\ }\href@noop {}
  {\bibfield  {journal} {\bibinfo  {journal} {Philos. T. Roy. Soc. A.}\
  }\textbf {\bibinfo {volume} {348}},\ \bibinfo {pages} {477} (\bibinfo {year}
  {1994})}\BibitemShut {NoStop}%
\bibitem [{\citenamefont {Finkenst\"{a}dt}\ and\ \citenamefont
  {Kuhbier}(1995)}]{FK:1995}%
  \BibitemOpen
  \bibfield  {author} {\bibinfo {author} {\bibfnamefont {B.}~\bibnamefont
  {Finkenst\"{a}dt}}\ and\ \bibinfo {author} {\bibfnamefont {P.}~\bibnamefont
  {Kuhbier}},\ }\bibfield  {title} {\enquote {\bibinfo {title} {Forecasting
  nonlinear economic time series: A simple test to accompany the nearest
  neighbor approach},}\ }\href@noop {} {\bibfield  {journal} {\bibinfo
  {journal} {Empiri. Econ.}\ }\textbf {\bibinfo {volume} {20}},\ \bibinfo
  {pages} {243} (\bibinfo {year} {1995})}\BibitemShut {NoStop}%
\bibitem [{\citenamefont {Parlitz}(1996)}]{Parlitz:1996}%
  \BibitemOpen
  \bibfield  {author} {\bibinfo {author} {\bibfnamefont {U.}~\bibnamefont
  {Parlitz}},\ }\bibfield  {title} {\enquote {\bibinfo {title} {Estimating
  model parameters from time series by autosynchronization},}\ }\href@noop {}
  {\bibfield  {journal} {\bibinfo  {journal} {Phys. Rev. Lett.}\ }\textbf
  {\bibinfo {volume} {76}},\ \bibinfo {pages} {1232} (\bibinfo {year}
  {1996})}\BibitemShut {NoStop}%
\bibitem [{\citenamefont {Schiff}\ \emph {et~al.}(1996)\citenamefont {Schiff},
  \citenamefont {So}, \citenamefont {Chang}, \citenamefont {Burke},\ and\
  \citenamefont {Sauer}}]{SSCBS:1996}%
  \BibitemOpen
  \bibfield  {author} {\bibinfo {author} {\bibfnamefont {S.~J.}\ \bibnamefont
  {Schiff}}, \bibinfo {author} {\bibfnamefont {P.}~\bibnamefont {So}}, \bibinfo
  {author} {\bibfnamefont {T.}~\bibnamefont {Chang}}, \bibinfo {author}
  {\bibfnamefont {R.~E.}\ \bibnamefont {Burke}}, \ and\ \bibinfo {author}
  {\bibfnamefont {T.}~\bibnamefont {Sauer}},\ }\bibfield  {title} {\enquote
  {\bibinfo {title} {Detecting dynamical interdependence and generalized
  synchrony through mutual prediction in a neural ensemble},}\ }\href {\doibase
  10.1103/PhysRevE.54.6708} {\bibfield  {journal} {\bibinfo  {journal} {Phys.
  Rev. E}\ }\textbf {\bibinfo {volume} {54}},\ \bibinfo {pages} {6708}
  (\bibinfo {year} {1996})}\BibitemShut {NoStop}%
\bibitem [{\citenamefont {Szpiro}(1997)}]{Szpiro:1997}%
  \BibitemOpen
  \bibfield  {author} {\bibinfo {author} {\bibfnamefont {G.~G.}\ \bibnamefont
  {Szpiro}},\ }\bibfield  {title} {\enquote {\bibinfo {title} {Forecasting
  chaotic time series with genetic algorithms},}\ }\href {\doibase
  10.1103/PhysRevE.55.2557} {\bibfield  {journal} {\bibinfo  {journal} {Phys.
  Rev. E}\ }\textbf {\bibinfo {volume} {55}},\ \bibinfo {pages} {2557}
  (\bibinfo {year} {1997})}\BibitemShut {NoStop}%
\bibitem [{\citenamefont {Hegger}\ \emph {et~al.}(1999)\citenamefont {Hegger},
  \citenamefont {Kantz},\ and\ \citenamefont {Schreiber}}]{HKS:1999}%
  \BibitemOpen
  \bibfield  {author} {\bibinfo {author} {\bibfnamefont {R.}~\bibnamefont
  {Hegger}}, \bibinfo {author} {\bibfnamefont {H.}~\bibnamefont {Kantz}}, \
  and\ \bibinfo {author} {\bibfnamefont {T.}~\bibnamefont {Schreiber}},\
  }\bibfield  {title} {\enquote {\bibinfo {title} {Practical implementation of
  nonlinear time series methods: The tisean package},}\ }\href@noop {}
  {\bibfield  {journal} {\bibinfo  {journal} {Chaos}\ }\textbf {\bibinfo
  {volume} {9}},\ \bibinfo {pages} {413} (\bibinfo {year} {1999})}\BibitemShut
  {NoStop}%
\bibitem [{\citenamefont {Hegger}\ \emph {et~al.}(2000)\citenamefont {Hegger},
  \citenamefont {Kantz}, \citenamefont {Matassini},\ and\ \citenamefont
  {Schreiber}}]{HKMS:2000}%
  \BibitemOpen
  \bibfield  {author} {\bibinfo {author} {\bibfnamefont {R.}~\bibnamefont
  {Hegger}}, \bibinfo {author} {\bibfnamefont {H.}~\bibnamefont {Kantz}},
  \bibinfo {author} {\bibfnamefont {L.}~\bibnamefont {Matassini}}, \ and\
  \bibinfo {author} {\bibfnamefont {T.}~\bibnamefont {Schreiber}},\ }\bibfield
  {title} {\enquote {\bibinfo {title} {Coping with nonstationarity by
  overembedding},}\ }\href {\doibase 10.1103/PhysRevLett.84.4092} {\bibfield
  {journal} {\bibinfo  {journal} {Phys. Rev. Lett.}\ }\textbf {\bibinfo
  {volume} {84}},\ \bibinfo {pages} {4092} (\bibinfo {year}
  {2000})}\BibitemShut {NoStop}%
\bibitem [{\citenamefont {Sello}(2001)}]{Sello:2001}%
  \BibitemOpen
  \bibfield  {author} {\bibinfo {author} {\bibfnamefont {S.}~\bibnamefont
  {Sello}},\ }\bibfield  {title} {\enquote {\bibinfo {title} {Solar cycle
  forecasting: a nonlinear dynamics approach},}\ }\href@noop {} {\bibfield
  {journal} {\bibinfo  {journal} {Astron. Astrophys.}\ }\textbf {\bibinfo
  {volume} {377}},\ \bibinfo {pages} {312} (\bibinfo {year}
  {2001})}\BibitemShut {NoStop}%
\bibitem [{\citenamefont {Matsumoto}\ \emph {et~al.}(2001)\citenamefont
  {Matsumoto}, \citenamefont {Nakajima}, \citenamefont {Saito}, \citenamefont
  {Sugi},\ and\ \citenamefont {Hamagishi}}]{MNSSH:2001}%
  \BibitemOpen
  \bibfield  {author} {\bibinfo {author} {\bibfnamefont {T.}~\bibnamefont
  {Matsumoto}}, \bibinfo {author} {\bibfnamefont {Y.}~\bibnamefont {Nakajima}},
  \bibinfo {author} {\bibfnamefont {M.}~\bibnamefont {Saito}}, \bibinfo
  {author} {\bibfnamefont {J.}~\bibnamefont {Sugi}}, \ and\ \bibinfo {author}
  {\bibfnamefont {H.}~\bibnamefont {Hamagishi}},\ }\bibfield  {title} {\enquote
  {\bibinfo {title} {Reconstructions and predictions of nonlinear dynamical
  systems: a hierarchical bayesian approach},}\ }\href@noop {} {\bibfield
  {journal} {\bibinfo  {journal} {IEEE Trans. Signal Proc.}\ }\textbf {\bibinfo
  {volume} {49}},\ \bibinfo {pages} {2138} (\bibinfo {year}
  {2001})}\BibitemShut {NoStop}%
\bibitem [{\citenamefont {Smith}(2002)}]{Smith:2002}%
  \BibitemOpen
  \bibfield  {author} {\bibinfo {author} {\bibfnamefont {L.~A.}\ \bibnamefont
  {Smith}},\ }\bibfield  {title} {\enquote {\bibinfo {title} {What might we
  learn from climate forecasts?}}\ }\href@noop {} {\bibfield  {journal}
  {\bibinfo  {journal} {Proc. Nat. Acad. Sci. (USA)}\ }\textbf {\bibinfo
  {volume} {19}},\ \bibinfo {pages} {2487} (\bibinfo {year}
  {2002})}\BibitemShut {NoStop}%
\bibitem [{\citenamefont {Judd}(2003)}]{Judd:2003}%
  \BibitemOpen
  \bibfield  {author} {\bibinfo {author} {\bibfnamefont {K.}~\bibnamefont
  {Judd}},\ }\bibfield  {title} {\enquote {\bibinfo {title} {Nonlinear state
  estimation, indistinguishable states, and the extended kalman filter},}\
  }\href@noop {} {\bibfield  {journal} {\bibinfo  {journal} {Physica D}\
  }\textbf {\bibinfo {volume} {183}},\ \bibinfo {pages} {273} (\bibinfo {year}
  {2003})}\BibitemShut {NoStop}%
\bibitem [{\citenamefont {Sauer}(2004)}]{Sauer:2004}%
  \BibitemOpen
  \bibfield  {author} {\bibinfo {author} {\bibfnamefont {T.~D.}\ \bibnamefont
  {Sauer}},\ }\bibfield  {title} {\enquote {\bibinfo {title} {Reconstruction of
  shared nonlinear dynamics in a network},}\ }\href {\doibase
  10.1103/PhysRevLett.93.198701} {\bibfield  {journal} {\bibinfo  {journal}
  {Phys. Rev. Lett.}\ }\textbf {\bibinfo {volume} {93}},\ \bibinfo {pages}
  {198701} (\bibinfo {year} {2004})}\BibitemShut {NoStop}%
\bibitem [{\citenamefont {Tao}\ \emph {et~al.}(2007)\citenamefont {Tao},
  \citenamefont {Zhang},\ and\ \citenamefont {Jiang}}]{TZJ:2007}%
  \BibitemOpen
  \bibfield  {author} {\bibinfo {author} {\bibfnamefont {C.}~\bibnamefont
  {Tao}}, \bibinfo {author} {\bibfnamefont {Y.}~\bibnamefont {Zhang}}, \ and\
  \bibinfo {author} {\bibfnamefont {J.~J.}\ \bibnamefont {Jiang}},\ }\bibfield
  {title} {\enquote {\bibinfo {title} {Estimating system parameters from
  chaotic time series with synchronization optimized by a genetic algorithm},}\
  }\href {\doibase 10.1103/PhysRevE.76.016209} {\bibfield  {journal} {\bibinfo
  {journal} {Phys. Rev. E}\ }\textbf {\bibinfo {volume} {76}},\ \bibinfo
  {pages} {016209} (\bibinfo {year} {2007})}\BibitemShut {NoStop}%
\bibitem [{\citenamefont {Crutchfield}\ and\ \citenamefont
  {McNamara}(1987)}]{CM:1987}%
  \BibitemOpen
  \bibfield  {author} {\bibinfo {author} {\bibfnamefont {J.~P.}\ \bibnamefont
  {Crutchfield}}\ and\ \bibinfo {author} {\bibfnamefont {B.}~\bibnamefont
  {McNamara}},\ }\bibfield  {title} {\enquote {\bibinfo {title} {Equations of
  motion from a data series},}\ }\href@noop {} {\bibfield  {journal} {\bibinfo
  {journal} {Complex Sys.}\ }\textbf {\bibinfo {volume} {1}},\ \bibinfo {pages}
  {417} (\bibinfo {year} {1987})}\BibitemShut {NoStop}%
\bibitem [{\citenamefont {Bollt}(2000)}]{Bollt:2000}%
  \BibitemOpen
  \bibfield  {author} {\bibinfo {author} {\bibfnamefont {E.~M.}\ \bibnamefont
  {Bollt}},\ }\bibfield  {title} {\enquote {\bibinfo {title} {Controlling chaos
  and the inverse frobenius-perron problem: global stabilization of arbitrary
  invariant measures},}\ }\href@noop {} {\bibfield  {journal} {\bibinfo
  {journal} {Int. J. Bif. Chaos}\ }\textbf {\bibinfo {volume} {10}},\ \bibinfo
  {pages} {1033} (\bibinfo {year} {2000})}\BibitemShut {NoStop}%
\bibitem [{\citenamefont {Yao}\ and\ \citenamefont {Bollt}(2007)}]{YB:2007}%
  \BibitemOpen
  \bibfield  {author} {\bibinfo {author} {\bibfnamefont {C.}~\bibnamefont
  {Yao}}\ and\ \bibinfo {author} {\bibfnamefont {E.~M.}\ \bibnamefont
  {Bollt}},\ }\bibfield  {title} {\enquote {\bibinfo {title} {Modeling and
  nonlinear parameter estimation with {Kronecker} product representation for
  coupled oscillators and spatiotemporal systems},}\ }\href@noop {} {\bibfield
  {journal} {\bibinfo  {journal} {Physica D}\ }\textbf {\bibinfo {volume}
  {227}},\ \bibinfo {pages} {78} (\bibinfo {year} {2007})}\BibitemShut
  {NoStop}%
\bibitem [{\citenamefont {Cand\`{e}s}\ \emph
  {et~al.}(2006{\natexlab{a}})\citenamefont {Cand\`{e}s}, \citenamefont
  {Romberg},\ and\ \citenamefont {Tao}}]{CRT:2006a}%
  \BibitemOpen
  \bibfield  {author} {\bibinfo {author} {\bibfnamefont {E.}~\bibnamefont
  {Cand\`{e}s}}, \bibinfo {author} {\bibfnamefont {J.}~\bibnamefont {Romberg}},
  \ and\ \bibinfo {author} {\bibfnamefont {T.}~\bibnamefont {Tao}},\ }\bibfield
   {title} {\enquote {\bibinfo {title} {Robust uncertainty principles: exact
  signal reconstruction from highly incomplete frequency information},}\
  }\href@noop {} {\bibfield  {journal} {\bibinfo  {journal} {IEEE Trans. Info.
  Theory}\ }\textbf {\bibinfo {volume} {52}},\ \bibinfo {pages} {489} (\bibinfo
  {year} {2006}{\natexlab{a}})}\BibitemShut {NoStop}%
\bibitem [{\citenamefont {Cand\`{e}s}\ \emph
  {et~al.}(2006{\natexlab{b}})\citenamefont {Cand\`{e}s}, \citenamefont
  {Romberg},\ and\ \citenamefont {Tao}}]{CRT:2006b}%
  \BibitemOpen
  \bibfield  {author} {\bibinfo {author} {\bibfnamefont {E.}~\bibnamefont
  {Cand\`{e}s}}, \bibinfo {author} {\bibfnamefont {J.}~\bibnamefont {Romberg}},
  \ and\ \bibinfo {author} {\bibfnamefont {T.}~\bibnamefont {Tao}},\ }\bibfield
   {title} {\enquote {\bibinfo {title} {Stable signal recovery from incomplete
  and inaccurate measurements},}\ }\href@noop {} {\bibfield  {journal}
  {\bibinfo  {journal} {Comm. Pure Appl. Math.}\ }\textbf {\bibinfo {volume}
  {59}},\ \bibinfo {pages} {1207} (\bibinfo {year}
  {2006}{\natexlab{b}})}\BibitemShut {NoStop}%
\bibitem [{\citenamefont {Cande\`s}(2006)}]{Candes:2006}%
  \BibitemOpen
  \bibfield  {author} {\bibinfo {author} {\bibfnamefont {E.}~\bibnamefont
  {Cande\`s}},\ }\bibfield  {title} {\enquote {\bibinfo {title} {Compressive
  sampling},}\ }in\ \href@noop {} {\emph {\bibinfo {booktitle} {Proceedings of
  the International Congress of Mathematicians}}},\ Vol.~\bibinfo {volume} {3}\
  (\bibinfo {organization} {Madrid, Spain},\ \bibinfo {year} {2006})\ pp.\
  \bibinfo {pages} {1433--1452}\BibitemShut {NoStop}%
\bibitem [{\citenamefont {Donoho}(2006)}]{Donoho:2006}%
  \BibitemOpen
  \bibfield  {author} {\bibinfo {author} {\bibfnamefont {D.}~\bibnamefont
  {Donoho}},\ }\bibfield  {title} {\enquote {\bibinfo {title} {Compressed
  sensing},}\ }\href@noop {} {\bibfield  {journal} {\bibinfo  {journal} {IEEE
  Trans. Info. Theory}\ }\textbf {\bibinfo {volume} {52}},\ \bibinfo {pages}
  {1289} (\bibinfo {year} {2006})}\BibitemShut {NoStop}%
\bibitem [{\citenamefont {Baraniuk}(2007)}]{Baraniuk:2007}%
  \BibitemOpen
  \bibfield  {author} {\bibinfo {author} {\bibfnamefont {R.~G.}\ \bibnamefont
  {Baraniuk}},\ }\bibfield  {title} {\enquote {\bibinfo {title} {Compressed
  sensing},}\ }\href@noop {} {\bibfield  {journal} {\bibinfo  {journal} {IEEE
  Signal Process. Mag.}\ }\textbf {\bibinfo {volume} {24}},\ \bibinfo {pages}
  {118} (\bibinfo {year} {2007})}\BibitemShut {NoStop}%
\bibitem [{\citenamefont {Cande\`s}\ and\ \citenamefont
  {Wakin}(2008)}]{CW:2008}%
  \BibitemOpen
  \bibfield  {author} {\bibinfo {author} {\bibfnamefont {E.}~\bibnamefont
  {Cande\`s}}\ and\ \bibinfo {author} {\bibfnamefont {M.}~\bibnamefont
  {Wakin}},\ }\bibfield  {title} {\enquote {\bibinfo {title} {An introduction
  to compressive sampling},}\ }\href@noop {} {\bibfield  {journal} {\bibinfo
  {journal} {IEEE Signal Process. Mag.}\ }\textbf {\bibinfo {volume} {25}},\
  \bibinfo {pages} {21} (\bibinfo {year} {2008})}\BibitemShut {NoStop}%
\bibitem [{\citenamefont {Wang}\ \emph
  {et~al.}(2011{\natexlab{a}})\citenamefont {Wang}, \citenamefont {Yang},
  \citenamefont {Lai}, \citenamefont {Kovanis},\ and\ \citenamefont
  {Grebogi}}]{WYLKG:2011}%
  \BibitemOpen
  \bibfield  {author} {\bibinfo {author} {\bibfnamefont {W.-X.}\ \bibnamefont
  {Wang}}, \bibinfo {author} {\bibfnamefont {R.}~\bibnamefont {Yang}}, \bibinfo
  {author} {\bibfnamefont {Y.-C.}\ \bibnamefont {Lai}}, \bibinfo {author}
  {\bibfnamefont {V.}~\bibnamefont {Kovanis}}, \ and\ \bibinfo {author}
  {\bibfnamefont {C.}~\bibnamefont {Grebogi}},\ }\bibfield  {title} {\enquote
  {\bibinfo {title} {Predicting catastrophes in nonlinear dynamical systems by
  compressive sensing},}\ }\href@noop {} {\bibfield  {journal} {\bibinfo
  {journal} {Phys. Rev. Lett.}\ }\textbf {\bibinfo {volume} {106}},\ \bibinfo
  {pages} {154101} (\bibinfo {year} {2011}{\natexlab{a}})}\BibitemShut
  {NoStop}%
\bibitem [{\citenamefont {Wang}\ \emph
  {et~al.}(2011{\natexlab{b}})\citenamefont {Wang}, \citenamefont {Lai},
  \citenamefont {Grebogi},\ and\ \citenamefont {Ye}}]{WLGY:2011}%
  \BibitemOpen
  \bibfield  {author} {\bibinfo {author} {\bibfnamefont {W.-X.}\ \bibnamefont
  {Wang}}, \bibinfo {author} {\bibfnamefont {Y.-C.}\ \bibnamefont {Lai}},
  \bibinfo {author} {\bibfnamefont {C.}~\bibnamefont {Grebogi}}, \ and\
  \bibinfo {author} {\bibfnamefont {J.-P.}\ \bibnamefont {Ye}},\ }\bibfield
  {title} {\enquote {\bibinfo {title} {Network reconstruction based on
  evolutionary-game data via compressive sensing},}\ }\href@noop {} {\bibfield
  {journal} {\bibinfo  {journal} {Phys. Rev. X}\ }\textbf {\bibinfo {volume}
  {1}},\ \bibinfo {pages} {021021} (\bibinfo {year}
  {2011}{\natexlab{b}})}\BibitemShut {NoStop}%
\bibitem [{\citenamefont {Wang}\ \emph
  {et~al.}(2011{\natexlab{c}})\citenamefont {Wang}, \citenamefont {Yang},
  \citenamefont {Lai}, \citenamefont {Kovanis},\ and\ \citenamefont
  {Harrison}}]{WYLKH:2011}%
  \BibitemOpen
  \bibfield  {author} {\bibinfo {author} {\bibfnamefont {W.-X.}\ \bibnamefont
  {Wang}}, \bibinfo {author} {\bibfnamefont {R.}~\bibnamefont {Yang}}, \bibinfo
  {author} {\bibfnamefont {Y.-C.}\ \bibnamefont {Lai}}, \bibinfo {author}
  {\bibfnamefont {V.}~\bibnamefont {Kovanis}}, \ and\ \bibinfo {author}
  {\bibfnamefont {M.~A.~F.}\ \bibnamefont {Harrison}},\ }\bibfield  {title}
  {\enquote {\bibinfo {title} {Time-series-based prediction of complex
  oscillator networks via compressive sensing},}\ }\href@noop {} {\bibfield
  {journal} {\bibinfo  {journal} {EPL (Europhys. Lett.)}\ }\textbf {\bibinfo
  {volume} {94}},\ \bibinfo {pages} {48006} (\bibinfo {year}
  {2011}{\natexlab{c}})}\BibitemShut {NoStop}%
\bibitem [{\citenamefont {Su}\ \emph {et~al.}(2012{\natexlab{a}})\citenamefont
  {Su}, \citenamefont {Ni}, \citenamefont {Wang},\ and\ \citenamefont
  {Lai}}]{SNWL:2012}%
  \BibitemOpen
  \bibfield  {author} {\bibinfo {author} {\bibfnamefont {R.-Q.}\ \bibnamefont
  {Su}}, \bibinfo {author} {\bibfnamefont {X.}~\bibnamefont {Ni}}, \bibinfo
  {author} {\bibfnamefont {W.-X.}\ \bibnamefont {Wang}}, \ and\ \bibinfo
  {author} {\bibfnamefont {Y.-C.}\ \bibnamefont {Lai}},\ }\bibfield  {title}
  {\enquote {\bibinfo {title} {Forecasting synchronizability of complex
  networks from data},}\ }\href@noop {} {\bibfield  {journal} {\bibinfo
  {journal} {Phys. Rev. E}\ }\textbf {\bibinfo {volume} {85}},\ \bibinfo
  {pages} {056220} (\bibinfo {year} {2012}{\natexlab{a}})}\BibitemShut
  {NoStop}%
\bibitem [{\citenamefont {Su}\ \emph {et~al.}(2012{\natexlab{b}})\citenamefont
  {Su}, \citenamefont {Wang},\ and\ \citenamefont {Lai}}]{SWL:2012}%
  \BibitemOpen
  \bibfield  {author} {\bibinfo {author} {\bibfnamefont {R.-Q.}\ \bibnamefont
  {Su}}, \bibinfo {author} {\bibfnamefont {W.-X.}\ \bibnamefont {Wang}}, \ and\
  \bibinfo {author} {\bibfnamefont {Y.-C.}\ \bibnamefont {Lai}},\ }\bibfield
  {title} {\enquote {\bibinfo {title} {Detecting hidden nodes in complex
  networks from time series},}\ }\href@noop {} {\bibfield  {journal} {\bibinfo
  {journal} {Phys. Rev. E}\ }\textbf {\bibinfo {volume} {85}},\ \bibinfo
  {pages} {065201} (\bibinfo {year} {2012}{\natexlab{b}})}\BibitemShut
  {NoStop}%
\bibitem [{\citenamefont {Su}\ \emph {et~al.}(2014{\natexlab{a}})\citenamefont
  {Su}, \citenamefont {Lai},\ and\ \citenamefont {Wang}}]{SWL:2014}%
  \BibitemOpen
  \bibfield  {author} {\bibinfo {author} {\bibfnamefont {R.-Q.}\ \bibnamefont
  {Su}}, \bibinfo {author} {\bibfnamefont {Y.-C.}\ \bibnamefont {Lai}}, \ and\
  \bibinfo {author} {\bibfnamefont {X.}~\bibnamefont {Wang}},\ }\bibfield
  {title} {\enquote {\bibinfo {title} {Identifying chaotic fitzhugh-nagumo
  neurons using compressive sensing},}\ }\href@noop {} {\bibfield  {journal}
  {\bibinfo  {journal} {Entropy}\ }\textbf {\bibinfo {volume} {16}},\ \bibinfo
  {pages} {3889} (\bibinfo {year} {2014}{\natexlab{a}})}\BibitemShut {NoStop}%
\bibitem [{\citenamefont {Su}\ \emph {et~al.}(2014{\natexlab{b}})\citenamefont
  {Su}, \citenamefont {Lai}, \citenamefont {Wang},\ and\ \citenamefont
  {Do}}]{SLWD:2014}%
  \BibitemOpen
  \bibfield  {author} {\bibinfo {author} {\bibfnamefont {R.-Q.}\ \bibnamefont
  {Su}}, \bibinfo {author} {\bibfnamefont {Y.-C.}\ \bibnamefont {Lai}},
  \bibinfo {author} {\bibfnamefont {X.}~\bibnamefont {Wang}}, \ and\ \bibinfo
  {author} {\bibfnamefont {Y.-H.}\ \bibnamefont {Do}},\ }\bibfield  {title}
  {\enquote {\bibinfo {title} {Uncovering hidden nodes in complex networks in
  the presence of noise},}\ }\href@noop {} {\bibfield  {journal} {\bibinfo
  {journal} {Sci. Rep.}\ }\textbf {\bibinfo {volume} {4}},\ \bibinfo {pages}
  {3944} (\bibinfo {year} {2014}{\natexlab{b}})}\BibitemShut {NoStop}%
\bibitem [{\citenamefont {Shen}\ \emph {et~al.}(2014)\citenamefont {Shen},
  \citenamefont {Wang}, \citenamefont {Fan}, \citenamefont {Di},\ and\
  \citenamefont {Lai}}]{SWFDL:2014}%
  \BibitemOpen
  \bibfield  {author} {\bibinfo {author} {\bibfnamefont {Z.}~\bibnamefont
  {Shen}}, \bibinfo {author} {\bibfnamefont {W.-X.}\ \bibnamefont {Wang}},
  \bibinfo {author} {\bibfnamefont {Y.}~\bibnamefont {Fan}}, \bibinfo {author}
  {\bibfnamefont {Z.}~\bibnamefont {Di}}, \ and\ \bibinfo {author}
  {\bibfnamefont {Y.-C.}\ \bibnamefont {Lai}},\ }\bibfield  {title} {\enquote
  {\bibinfo {title} {Reconstructing propagation networks with natural diversity
  and identifying hidden sources},}\ }\href@noop {} {\bibfield  {journal}
  {\bibinfo  {journal} {Nat. Commun.}\ }\textbf {\bibinfo {volume} {5}},\
  \bibinfo {pages} {4323} (\bibinfo {year} {2014})}\BibitemShut {NoStop}%
\bibitem [{\citenamefont {Su}\ \emph {et~al.}(2016)\citenamefont {Su},
  \citenamefont {Wang}, \citenamefont {Wang},\ and\ \citenamefont
  {Lai}}]{SWWL:2016}%
  \BibitemOpen
  \bibfield  {author} {\bibinfo {author} {\bibfnamefont {R.-Q.}\ \bibnamefont
  {Su}}, \bibinfo {author} {\bibfnamefont {W.-W.}\ \bibnamefont {Wang}},
  \bibinfo {author} {\bibfnamefont {X.}~\bibnamefont {Wang}}, \ and\ \bibinfo
  {author} {\bibfnamefont {Y.-C.}\ \bibnamefont {Lai}},\ }\bibfield  {title}
  {\enquote {\bibinfo {title} {Data based reconstruction of complex geospatial
  networks, nodal positioning, and detection of hidden node},}\ }\href@noop {}
  {\bibfield  {journal} {\bibinfo  {journal} {R. Soc. Open Sci.}\ }\textbf
  {\bibinfo {volume} {3}},\ \bibinfo {pages} {150577} (\bibinfo {year}
  {2016})}\BibitemShut {NoStop}%
\bibitem [{EMB()}]{EMB}%
  \BibitemOpen
  \href@noop {} {\bibinfo  {journal} {Professor E. M. Bollt from Clarkson
  University conceived the same idea of exploiting sparse optimization for
  discovering system equations from data (unpublished, private communication)}\
  }\BibitemShut {NoStop}%
\bibitem [{\citenamefont {AlMomani}\ \emph {et~al.}(2020)\citenamefont
  {AlMomani}, \citenamefont {Jie},\ and\ \citenamefont {Bollt}}]{ASB:2020}%
  \BibitemOpen
\bibfield  {journal} {  }\bibfield  {author} {\bibinfo {author} {\bibfnamefont
  {A.~A.~R.}\ \bibnamefont {AlMomani}}, \bibinfo {author} {\bibfnamefont
  {S.}~\bibnamefont {Jie}}, \ and\ \bibinfo {author} {\bibfnamefont {E.~M.}\
  \bibnamefont {Bollt}},\ }\bibfield  {title} {\enquote {\bibinfo {title} {How
  entropic regression beats the outliers problem in nonlinear system
  identification},}\ }\href@noop {} {\bibfield  {journal} {\bibinfo  {journal}
  {Chaos}\ }\textbf {\bibinfo {volume} {30}},\ \bibinfo {pages} {013107}
  (\bibinfo {year} {2020})}\BibitemShut {NoStop}%
\bibitem [{\citenamefont {Yang}\ \emph {et~al.}(2012)\citenamefont {Yang},
  \citenamefont {Lai},\ and\ \citenamefont {Grebogi}}]{YLG:2012}%
  \BibitemOpen
  \bibfield  {author} {\bibinfo {author} {\bibfnamefont {R.}~\bibnamefont
  {Yang}}, \bibinfo {author} {\bibfnamefont {Y.-C.}\ \bibnamefont {Lai}}, \
  and\ \bibinfo {author} {\bibfnamefont {C.}~\bibnamefont {Grebogi}},\
  }\bibfield  {title} {\enquote {\bibinfo {title} {Forecasting the future: is
  it possible for time-varying nonlinear dynamical systems?}}\ }\href@noop {}
  {\bibfield  {journal} {\bibinfo  {journal} {Chaos}\ }\textbf {\bibinfo
  {volume} {22}},\ \bibinfo {pages} {033119} (\bibinfo {year}
  {2012})}\BibitemShut {NoStop}%
\bibitem [{\citenamefont {Lorenz}(1963)}]{Lorenz:1963}%
  \BibitemOpen
  \bibfield  {author} {\bibinfo {author} {\bibfnamefont {E.~N.}\ \bibnamefont
  {Lorenz}},\ }\bibfield  {title} {\enquote {\bibinfo {title} {Deterministic
  nonperiodic flow},}\ }\href@noop {} {\bibfield  {journal} {\bibinfo
  {journal} {J. Atmos. Sci.}\ }\textbf {\bibinfo {volume} {20}},\ \bibinfo
  {pages} {130} (\bibinfo {year} {1963})}\BibitemShut {NoStop}%
\bibitem [{\citenamefont {R\"{o}ssler}(1976)}]{Rossler:1976}%
  \BibitemOpen
  \bibfield  {author} {\bibinfo {author} {\bibfnamefont {O.~E.}\ \bibnamefont
  {R\"{o}ssler}},\ }\bibfield  {title} {\enquote {\bibinfo {title} {Equation
  for continuous chaos},}\ }\href@noop {} {\bibfield  {journal} {\bibinfo
  {journal} {Phys. Lett. A}\ }\textbf {\bibinfo {volume} {57}},\ \bibinfo
  {pages} {397} (\bibinfo {year} {1976})}\BibitemShut {NoStop}%
\bibitem [{\citenamefont {Grebogi}\ \emph {et~al.}(1983)\citenamefont
  {Grebogi}, \citenamefont {Ott},\ and\ \citenamefont {Yorke}}]{GOY:1983a}%
  \BibitemOpen
  \bibfield  {author} {\bibinfo {author} {\bibfnamefont {C.}~\bibnamefont
  {Grebogi}}, \bibinfo {author} {\bibfnamefont {E.}~\bibnamefont {Ott}}, \ and\
  \bibinfo {author} {\bibfnamefont {J.~A.}\ \bibnamefont {Yorke}},\ }\bibfield
  {title} {\enquote {\bibinfo {title} {Crises, sudden changes in chaotic
  attractors and chaotic transients},}\ }\href@noop {} {\bibfield  {journal}
  {\bibinfo  {journal} {Physica D}\ }\textbf {\bibinfo {volume} {7}},\ \bibinfo
  {pages} {181} (\bibinfo {year} {1983})}\BibitemShut {NoStop}%
\bibitem [{\citenamefont {Dhamala}\ and\ \citenamefont {Lai}(1999)}]{DL:1999}%
  \BibitemOpen
  \bibfield  {author} {\bibinfo {author} {\bibfnamefont {M.}~\bibnamefont
  {Dhamala}}\ and\ \bibinfo {author} {\bibfnamefont {Y.-C.}\ \bibnamefont
  {Lai}},\ }\bibfield  {title} {\enquote {\bibinfo {title} {Controlling
  transient chaos in deterministic flows with applications to electrical power
  systems and ecology},}\ }\href@noop {} {\bibfield  {journal} {\bibinfo
  {journal} {Phys. Rev. E}\ }\textbf {\bibinfo {volume} {59}},\ \bibinfo
  {pages} {1646} (\bibinfo {year} {1999})}\BibitemShut {NoStop}%
\bibitem [{\citenamefont {McCann}\ and\ \citenamefont
  {Yodzis}(1994)}]{McY:1994}%
  \BibitemOpen
  \bibfield  {author} {\bibinfo {author} {\bibfnamefont {K.}~\bibnamefont
  {McCann}}\ and\ \bibinfo {author} {\bibfnamefont {P.}~\bibnamefont
  {Yodzis}},\ }\bibfield  {title} {\enquote {\bibinfo {title} {Nonlinear
  dynamics and population disappearances},}\ }\href@noop {} {\bibfield
  {journal} {\bibinfo  {journal} {Ame. Naturalist}\ }\textbf {\bibinfo {volume}
  {144}},\ \bibinfo {pages} {873} (\bibinfo {year} {1994})}\BibitemShut
  {NoStop}%
\bibitem [{\citenamefont {Hastings}\ \emph {et~al.}(2018)\citenamefont
  {Hastings}, \citenamefont {Abbott}, \citenamefont {Cuddington}, \citenamefont
  {Francis}, \citenamefont {Gellner}, \citenamefont {Lai}, \citenamefont
  {Morozov}, \citenamefont {Petrivskii}, \citenamefont {Scranton},\ and\
  \citenamefont {Zeeman}}]{HACFGLMPSZ:2018}%
  \BibitemOpen
  \bibfield  {author} {\bibinfo {author} {\bibfnamefont {A.}~\bibnamefont
  {Hastings}}, \bibinfo {author} {\bibfnamefont {K.~C.}\ \bibnamefont
  {Abbott}}, \bibinfo {author} {\bibfnamefont {K.}~\bibnamefont {Cuddington}},
  \bibinfo {author} {\bibfnamefont {T.}~\bibnamefont {Francis}}, \bibinfo
  {author} {\bibfnamefont {G.}~\bibnamefont {Gellner}}, \bibinfo {author}
  {\bibfnamefont {Y.-C.}\ \bibnamefont {Lai}}, \bibinfo {author} {\bibfnamefont
  {A.}~\bibnamefont {Morozov}}, \bibinfo {author} {\bibfnamefont
  {S.}~\bibnamefont {Petrivskii}}, \bibinfo {author} {\bibfnamefont
  {K.}~\bibnamefont {Scranton}}, \ and\ \bibinfo {author} {\bibfnamefont
  {M.~L.}\ \bibnamefont {Zeeman}},\ }\bibfield  {title} {\enquote {\bibinfo
  {title} {Transient phenomena in ecology},}\ }\href@noop {} {\bibfield
  {journal} {\bibinfo  {journal} {Science}\ }\textbf {\bibinfo {volume}
  {361}},\ \bibinfo {pages} {eaat6412} (\bibinfo {year} {2018})}\BibitemShut
  {NoStop}%
\bibitem [{\citenamefont {Wang}\ \emph {et~al.}(2016)\citenamefont {Wang},
  \citenamefont {Lai},\ and\ \citenamefont {Grebogi}}]{WLG:2016}%
  \BibitemOpen
  \bibfield  {author} {\bibinfo {author} {\bibfnamefont {W.}~\bibnamefont
  {Wang}}, \bibinfo {author} {\bibfnamefont {Y.-C.}\ \bibnamefont {Lai}}, \
  and\ \bibinfo {author} {\bibfnamefont {C.}~\bibnamefont {Grebogi}},\
  }\bibfield  {title} {\enquote {\bibinfo {title} {Data based identification
  and prediction of nonlinear and complex dynamical systems},}\ }\href@noop {}
  {\bibfield  {journal} {\bibinfo  {journal} {Phys. Rep.}\ }\textbf {\bibinfo
  {volume} {644}},\ \bibinfo {pages} {1} (\bibinfo {year} {2016})}\BibitemShut
  {NoStop}%
\bibitem [{\citenamefont {Newman}(2010)}]{Newman:book}%
  \BibitemOpen
  \bibfield  {author} {\bibinfo {author} {\bibfnamefont {M.~E.~J.}\
  \bibnamefont {Newman}},\ }\href@noop {} {\emph {\bibinfo {title} {Networks:
  An Introduction}}}\ (\bibinfo  {publisher} {Oxford University Press},\
  \bibinfo {address} {Oxford, UK},\ \bibinfo {year} {2010})\BibitemShut
  {NoStop}%
\bibitem [{\citenamefont {Nowak}\ and\ \citenamefont {May}(1992)}]{NM:1992}%
  \BibitemOpen
  \bibfield  {author} {\bibinfo {author} {\bibfnamefont {M.~A.}\ \bibnamefont
  {Nowak}}\ and\ \bibinfo {author} {\bibfnamefont {R.~M.}\ \bibnamefont
  {May}},\ }\bibfield  {title} {\enquote {\bibinfo {title} {Evolutionary games
  and spatial chaos},}\ }\href@noop {} {\bibfield  {journal} {\bibinfo
  {journal} {Nature (London)}\ }\textbf {\bibinfo {volume} {359}},\ \bibinfo
  {pages} {826} (\bibinfo {year} {1992})}\BibitemShut {NoStop}%
\bibitem [{\citenamefont {Hauert}\ and\ \citenamefont
  {Doebeli}(2004)}]{HD:2004}%
  \BibitemOpen
  \bibfield  {author} {\bibinfo {author} {\bibfnamefont {C.}~\bibnamefont
  {Hauert}}\ and\ \bibinfo {author} {\bibfnamefont {M.}~\bibnamefont
  {Doebeli}},\ }\bibfield  {title} {\enquote {\bibinfo {title} {Spatial
  structure often inhibits the evolution of cooperation in the snowdrift
  game},}\ }\href@noop {} {\bibfield  {journal} {\bibinfo  {journal} {Nature
  (London)}\ }\textbf {\bibinfo {volume} {428}},\ \bibinfo {pages} {643}
  (\bibinfo {year} {2004})}\BibitemShut {NoStop}%
\bibitem [{\citenamefont {Szab\'o}\ and\ \citenamefont
  {T\ifmmode~\mbox{\H{o}}\else \H{o}\fi{}ke}(1998)}]{ST:1998}%
  \BibitemOpen
  \bibfield  {author} {\bibinfo {author} {\bibfnamefont {G.}~\bibnamefont
  {Szab\'o}}\ and\ \bibinfo {author} {\bibfnamefont {C.}~\bibnamefont
  {T\ifmmode~\mbox{\H{o}}\else \H{o}\fi{}ke}},\ }\bibfield  {title} {\enquote
  {\bibinfo {title} {Evolutionary prisoner's dilemma game on a square
  lattice},}\ }\href {\doibase 10.1103/PhysRevE.58.69} {\bibfield  {journal}
  {\bibinfo  {journal} {Phys. Rev. E}\ }\textbf {\bibinfo {volume} {58}},\
  \bibinfo {pages} {69} (\bibinfo {year} {1998})}\BibitemShut {NoStop}%
\bibitem [{\citenamefont {Santos}\ \emph {et~al.}(2008)\citenamefont {Santos},
  \citenamefont {Santos},\ and\ \citenamefont {Pacheco}}]{SSP:2008}%
  \BibitemOpen
  \bibfield  {author} {\bibinfo {author} {\bibfnamefont {F.~C.}\ \bibnamefont
  {Santos}}, \bibinfo {author} {\bibfnamefont {M.~D.}\ \bibnamefont {Santos}},
  \ and\ \bibinfo {author} {\bibfnamefont {J.~M.}\ \bibnamefont {Pacheco}},\
  }\bibfield  {title} {\enquote {\bibinfo {title} {Social diversity promotes
  the emergence of cooperation in public goods games},}\ }\href@noop {}
  {\bibfield  {journal} {\bibinfo  {journal} {Nature}\ }\textbf {\bibinfo
  {volume} {454}},\ \bibinfo {pages} {213} (\bibinfo {year}
  {2008})}\BibitemShut {NoStop}%
\bibitem [{\citenamefont {Han}\ \emph {et~al.}(2015)\citenamefont {Han},
  \citenamefont {Shen}, \citenamefont {Wang},\ and\ \citenamefont
  {Di}}]{HSWD:2015}%
  \BibitemOpen
  \bibfield  {author} {\bibinfo {author} {\bibfnamefont {X.}~\bibnamefont
  {Han}}, \bibinfo {author} {\bibfnamefont {Z.}~\bibnamefont {Shen}}, \bibinfo
  {author} {\bibfnamefont {W.-X.}\ \bibnamefont {Wang}}, \ and\ \bibinfo
  {author} {\bibfnamefont {Z.}~\bibnamefont {Di}},\ }\bibfield  {title}
  {\enquote {\bibinfo {title} {Robust reconstruction of complex networks from
  sparse data},}\ }\href {\doibase 10.1103/PhysRevLett.114.028701} {\bibfield
  {journal} {\bibinfo  {journal} {Phys. Rev. Lett.}\ }\textbf {\bibinfo
  {volume} {114}},\ \bibinfo {pages} {028701} (\bibinfo {year}
  {2015})}\BibitemShut {NoStop}%
\bibitem [{\citenamefont {Li}\ \emph {et~al.}(2017)\citenamefont {Li},
  \citenamefont {Shen}, \citenamefont {Wang}, \citenamefont {Grebogi},\ and\
  \citenamefont {Lai}}]{LSWGL:2017}%
  \BibitemOpen
  \bibfield  {author} {\bibinfo {author} {\bibfnamefont {J.}~\bibnamefont
  {Li}}, \bibinfo {author} {\bibfnamefont {Z.}~\bibnamefont {Shen}}, \bibinfo
  {author} {\bibfnamefont {W.-X.}\ \bibnamefont {Wang}}, \bibinfo {author}
  {\bibfnamefont {C.}~\bibnamefont {Grebogi}}, \ and\ \bibinfo {author}
  {\bibfnamefont {Y.-C.}\ \bibnamefont {Lai}},\ }\bibfield  {title} {\enquote
  {\bibinfo {title} {Universal data-based method for reconstructing complex
  networks with binary-state dynamics},}\ }\href {\doibase
  10.1103/PhysRevE.95.032303} {\bibfield  {journal} {\bibinfo  {journal} {Phys.
  Rev. E}\ }\textbf {\bibinfo {volume} {95}},\ \bibinfo {pages} {032303}
  (\bibinfo {year} {2017})}\BibitemShut {NoStop}%
\bibitem [{\citenamefont {Santosa}\ and\ \citenamefont
  {Symes}(1986)}]{SS:1986}%
  \BibitemOpen
  \bibfield  {author} {\bibinfo {author} {\bibfnamefont {F.}~\bibnamefont
  {Santosa}}\ and\ \bibinfo {author} {\bibfnamefont {W.~W.}\ \bibnamefont
  {Symes}},\ }\bibfield  {title} {\enquote {\bibinfo {title} {Linear inversion
  of band-limited reflection seismograms},}\ }\href@noop {} {\bibfield
  {journal} {\bibinfo  {journal} {SIAM J. Sci. Stat. Comp.}\ }\textbf {\bibinfo
  {volume} {7}},\ \bibinfo {pages} {1307} (\bibinfo {year} {1986})}\BibitemShut
  {NoStop}%
\bibitem [{\citenamefont {Friedman}\ \emph {et~al.}(2001)\citenamefont
  {Friedman}, \citenamefont {Hastie},\ and\ \citenamefont
  {Tibshirani}}]{FHT:book}%
  \BibitemOpen
  \bibfield  {author} {\bibinfo {author} {\bibfnamefont {J.}~\bibnamefont
  {Friedman}}, \bibinfo {author} {\bibfnamefont {T.}~\bibnamefont {Hastie}}, \
  and\ \bibinfo {author} {\bibfnamefont {R.}~\bibnamefont {Tibshirani}},\
  }\href@noop {} {\emph {\bibinfo {title} {The Elements of Statistical
  Learning}}}\ (\bibinfo  {publisher} {Springer},\ \bibinfo {address}
  {Berlin},\ \bibinfo {year} {2001})\BibitemShut {NoStop}%
\bibitem [{\citenamefont {Pedregosa}\ \emph {et~al.}(2011)\citenamefont
  {Pedregosa}, \citenamefont {Varoquaux}, \citenamefont {Gramfort},
  \citenamefont {Michel}, \citenamefont {Thirion}, \citenamefont {Grisel},
  \citenamefont {Blondel}, \citenamefont {Prettenhofer}, \citenamefont {Weiss},
  \citenamefont {Dubourg} \emph {et~al.}}]{PVGMTGBPWD:2011}%
  \BibitemOpen
  \bibfield  {author} {\bibinfo {author} {\bibfnamefont {F.}~\bibnamefont
  {Pedregosa}}, \bibinfo {author} {\bibfnamefont {G.}~\bibnamefont
  {Varoquaux}}, \bibinfo {author} {\bibfnamefont {A.}~\bibnamefont {Gramfort}},
  \bibinfo {author} {\bibfnamefont {V.}~\bibnamefont {Michel}}, \bibinfo
  {author} {\bibfnamefont {B.}~\bibnamefont {Thirion}}, \bibinfo {author}
  {\bibfnamefont {O.}~\bibnamefont {Grisel}}, \bibinfo {author} {\bibfnamefont
  {M.}~\bibnamefont {Blondel}}, \bibinfo {author} {\bibfnamefont
  {P.}~\bibnamefont {Prettenhofer}}, \bibinfo {author} {\bibfnamefont
  {R.}~\bibnamefont {Weiss}}, \bibinfo {author} {\bibfnamefont
  {V.}~\bibnamefont {Dubourg}},  \emph {et~al.},\ }\bibfield  {title} {\enquote
  {\bibinfo {title} {Scikit-learn: Machine learning in python},}\ }\href@noop
  {} {\bibfield  {journal} {\bibinfo  {journal} {J. Mach. Learn. Res.}\
  }\textbf {\bibinfo {volume} {12}},\ \bibinfo {pages} {2825} (\bibinfo {year}
  {2011})}\BibitemShut {NoStop}%
\bibitem [{\citenamefont {Barrat}\ \emph {et~al.}(2008)\citenamefont {Barrat},
  \citenamefont {Barthelemy},\ and\ \citenamefont {Vespignani}}]{barrat2008}%
  \BibitemOpen
  \bibfield  {author} {\bibinfo {author} {\bibfnamefont {A.}~\bibnamefont
  {Barrat}}, \bibinfo {author} {\bibfnamefont {M.}~\bibnamefont {Barthelemy}},
  \ and\ \bibinfo {author} {\bibfnamefont {A.}~\bibnamefont {Vespignani}},\
  }\href@noop {} {\emph {\bibinfo {title} {Dynamical Processes on Complex
  Networks}}}\ (\bibinfo  {publisher} {Cambridge University Press},\ \bibinfo
  {year} {2008})\BibitemShut {NoStop}%
\bibitem [{\citenamefont {Kumar}\ \emph {et~al.}(2010)\citenamefont {Kumar},
  \citenamefont {Rotter},\ and\ \citenamefont {Aertsen}}]{KRA:2010}%
  \BibitemOpen
  \bibfield  {author} {\bibinfo {author} {\bibfnamefont {A.}~\bibnamefont
  {Kumar}}, \bibinfo {author} {\bibfnamefont {S.}~\bibnamefont {Rotter}}, \
  and\ \bibinfo {author} {\bibfnamefont {A.}~\bibnamefont {Aertsen}},\
  }\bibfield  {title} {\enquote {\bibinfo {title} {Spiking activity propagation
  in neuronal networks: reconciling different perspectives on neural coding},}\
  }\href@noop {} {\bibfield  {journal} {\bibinfo  {journal} {Nature Rev.
  Neuro.}\ }\textbf {\bibinfo {volume} {11}},\ \bibinfo {pages} {615} (\bibinfo
  {year} {2010})}\BibitemShut {NoStop}%
\bibitem [{\citenamefont {Szab{\'o}}\ and\ \citenamefont
  {Fath}(2007)}]{SF:2007}%
  \BibitemOpen
  \bibfield  {author} {\bibinfo {author} {\bibfnamefont {G.}~\bibnamefont
  {Szab{\'o}}}\ and\ \bibinfo {author} {\bibfnamefont {G.}~\bibnamefont
  {Fath}},\ }\bibfield  {title} {\enquote {\bibinfo {title} {Evolutionary games
  on graphs},}\ }\href@noop {} {\bibfield  {journal} {\bibinfo  {journal}
  {Phys. Rep.}\ }\textbf {\bibinfo {volume} {446}},\ \bibinfo {pages} {97}
  (\bibinfo {year} {2007})}\BibitemShut {NoStop}%
\bibitem [{\citenamefont {Pastor-Satorras}\ \emph {et~al.}(2015)\citenamefont
  {Pastor-Satorras}, \citenamefont {Castellano}, \citenamefont {Van~Mieghem},\
  and\ \citenamefont {Vespignani}}]{PSCVV:2015}%
  \BibitemOpen
  \bibfield  {author} {\bibinfo {author} {\bibfnamefont {R.}~\bibnamefont
  {Pastor-Satorras}}, \bibinfo {author} {\bibfnamefont {C.}~\bibnamefont
  {Castellano}}, \bibinfo {author} {\bibfnamefont {P.}~\bibnamefont
  {Van~Mieghem}}, \ and\ \bibinfo {author} {\bibfnamefont {A.}~\bibnamefont
  {Vespignani}},\ }\bibfield  {title} {\enquote {\bibinfo {title} {Epidemic
  processes in complex networks},}\ }\href {\doibase 10.1103/RevModPhys.87.925}
  {\bibfield  {journal} {\bibinfo  {journal} {Rev. Mod. Phys.}\ }\textbf
  {\bibinfo {volume} {87}},\ \bibinfo {pages} {925} (\bibinfo {year}
  {2015})}\BibitemShut {NoStop}%
\bibitem [{\citenamefont {Shao}\ \emph {et~al.}(2009)\citenamefont {Shao},
  \citenamefont {Havlin},\ and\ \citenamefont {Stanley}}]{SHS:2009}%
  \BibitemOpen
  \bibfield  {author} {\bibinfo {author} {\bibfnamefont {J.}~\bibnamefont
  {Shao}}, \bibinfo {author} {\bibfnamefont {S.}~\bibnamefont {Havlin}}, \ and\
  \bibinfo {author} {\bibfnamefont {H.~E.}\ \bibnamefont {Stanley}},\
  }\bibfield  {title} {\enquote {\bibinfo {title} {Dynamic opinion model and
  invasion percolation},}\ }\href {\doibase 10.1103/PhysRevLett.103.018701}
  {\bibfield  {journal} {\bibinfo  {journal} {Phys. Rev. Lett.}\ }\textbf
  {\bibinfo {volume} {103}},\ \bibinfo {pages} {018701} (\bibinfo {year}
  {2009})}\BibitemShut {NoStop}%
\bibitem [{\citenamefont {Granell}\ \emph {et~al.}(2013)\citenamefont
  {Granell}, \citenamefont {G{\'o}mez},\ and\ \citenamefont
  {Arenas}}]{GGA:2013}%
  \BibitemOpen
  \bibfield  {author} {\bibinfo {author} {\bibfnamefont {C.}~\bibnamefont
  {Granell}}, \bibinfo {author} {\bibfnamefont {S.}~\bibnamefont {G{\'o}mez}},
  \ and\ \bibinfo {author} {\bibfnamefont {A.}~\bibnamefont {Arenas}},\
  }\bibfield  {title} {\enquote {\bibinfo {title} {Dynamical interplay between
  awareness and epidemic spreading in multiplex networks},}\ }\href@noop {}
  {\bibfield  {journal} {\bibinfo  {journal} {Phys. Rev. Lett.}\ }\textbf
  {\bibinfo {volume} {111}},\ \bibinfo {pages} {128701} (\bibinfo {year}
  {2013})}\BibitemShut {NoStop}%
\bibitem [{\citenamefont {Santos}\ and\ \citenamefont
  {Pacheco}(2005)}]{SP:2005}%
  \BibitemOpen
  \bibfield  {author} {\bibinfo {author} {\bibfnamefont {F.~C.}\ \bibnamefont
  {Santos}}\ and\ \bibinfo {author} {\bibfnamefont {J.~M.}\ \bibnamefont
  {Pacheco}},\ }\bibfield  {title} {\enquote {\bibinfo {title} {Scale-free
  networks provide a unifying framework for the emergence of cooperation},}\
  }\href@noop {} {\bibfield  {journal} {\bibinfo  {journal} {Phys. Rev. Lett.}\
  }\textbf {\bibinfo {volume} {95}},\ \bibinfo {pages} {098104} (\bibinfo
  {year} {2005})}\BibitemShut {NoStop}%
\bibitem [{\citenamefont {Koseska}\ \emph {et~al.}(2013)\citenamefont
  {Koseska}, \citenamefont {Volkov},\ and\ \citenamefont {Kurths}}]{KVK:2013}%
  \BibitemOpen
  \bibfield  {author} {\bibinfo {author} {\bibfnamefont {A.}~\bibnamefont
  {Koseska}}, \bibinfo {author} {\bibfnamefont {E.}~\bibnamefont {Volkov}}, \
  and\ \bibinfo {author} {\bibfnamefont {J.}~\bibnamefont {Kurths}},\
  }\bibfield  {title} {\enquote {\bibinfo {title} {Oscillation quenching
  mechanisms: Amplitude vs. oscillation death},}\ }\href@noop {} {\bibfield
  {journal} {\bibinfo  {journal} {Phys. Rep.}\ }\textbf {\bibinfo {volume}
  {531}},\ \bibinfo {pages} {173} (\bibinfo {year} {2013})}\BibitemShut
  {NoStop}%
\bibitem [{\citenamefont {Buldyrev}\ \emph {et~al.}(2010)\citenamefont
  {Buldyrev}, \citenamefont {Parshani}, \citenamefont {Paul}, \citenamefont
  {Stanley},\ and\ \citenamefont {Havlin}}]{BPPSH:2010}%
  \BibitemOpen
  \bibfield  {author} {\bibinfo {author} {\bibfnamefont {S.~V.}\ \bibnamefont
  {Buldyrev}}, \bibinfo {author} {\bibfnamefont {R.}~\bibnamefont {Parshani}},
  \bibinfo {author} {\bibfnamefont {G.}~\bibnamefont {Paul}}, \bibinfo {author}
  {\bibfnamefont {H.~E.}\ \bibnamefont {Stanley}}, \ and\ \bibinfo {author}
  {\bibfnamefont {S.}~\bibnamefont {Havlin}},\ }\bibfield  {title} {\enquote
  {\bibinfo {title} {Catastrophic cascade of failures in interdependent
  networks},}\ }\href@noop {} {\bibfield  {journal} {\bibinfo  {journal}
  {Nature}\ }\textbf {\bibinfo {volume} {464}},\ \bibinfo {pages} {1025}
  (\bibinfo {year} {2010})}\BibitemShut {NoStop}%
\bibitem [{\citenamefont {Galbiati}\ \emph {et~al.}(2013)\citenamefont
  {Galbiati}, \citenamefont {Delpini},\ and\ \citenamefont
  {Battiston}}]{GDB:2013}%
  \BibitemOpen
  \bibfield  {author} {\bibinfo {author} {\bibfnamefont {M.}~\bibnamefont
  {Galbiati}}, \bibinfo {author} {\bibfnamefont {D.}~\bibnamefont {Delpini}}, \
  and\ \bibinfo {author} {\bibfnamefont {S.}~\bibnamefont {Battiston}},\
  }\bibfield  {title} {\enquote {\bibinfo {title} {The power to control},}\
  }\href@noop {} {\bibfield  {journal} {\bibinfo  {journal} {Nat. Phys.}\
  }\textbf {\bibinfo {volume} {9}},\ \bibinfo {pages} {126} (\bibinfo {year}
  {2013})}\BibitemShut {NoStop}%
\bibitem [{\citenamefont {Balcan}\ and\ \citenamefont
  {Vespignani}(2011)}]{BV:2011}%
  \BibitemOpen
  \bibfield  {author} {\bibinfo {author} {\bibfnamefont {D.}~\bibnamefont
  {Balcan}}\ and\ \bibinfo {author} {\bibfnamefont {A.}~\bibnamefont
  {Vespignani}},\ }\bibfield  {title} {\enquote {\bibinfo {title} {Phase
  transitions in contagion processes mediated by recurrent mobility
  patterns},}\ }\href@noop {} {\bibfield  {journal} {\bibinfo  {journal} {Nat.
  Phys.}\ }\textbf {\bibinfo {volume} {7}},\ \bibinfo {pages} {581} (\bibinfo
  {year} {2011})}\BibitemShut {NoStop}%
\bibitem [{\citenamefont {Sood}\ and\ \citenamefont {Redner}(2005)}]{SR:2005}%
  \BibitemOpen
  \bibfield  {author} {\bibinfo {author} {\bibfnamefont {V.}~\bibnamefont
  {Sood}}\ and\ \bibinfo {author} {\bibfnamefont {S.}~\bibnamefont {Redner}},\
  }\bibfield  {title} {\enquote {\bibinfo {title} {Voter model on heterogeneous
  graphs},}\ }\href@noop {} {\bibfield  {journal} {\bibinfo  {journal} {Phys.
  Rev. Lett.}\ }\textbf {\bibinfo {volume} {94}},\ \bibinfo {pages} {178701}
  (\bibinfo {year} {2005})}\BibitemShut {NoStop}%
\bibitem [{\citenamefont {Pastor-Satorras}\ and\ \citenamefont
  {Vespignani}(2001)}]{PSV:2001}%
  \BibitemOpen
  \bibfield  {author} {\bibinfo {author} {\bibfnamefont {R.}~\bibnamefont
  {Pastor-Satorras}}\ and\ \bibinfo {author} {\bibfnamefont {A.}~\bibnamefont
  {Vespignani}},\ }\bibfield  {title} {\enquote {\bibinfo {title} {Epidemic
  spreading in scale-free networks},}\ }\href@noop {} {\bibfield  {journal}
  {\bibinfo  {journal} {Phys. Rev. Lett.}\ }\textbf {\bibinfo {volume} {86}},\
  \bibinfo {pages} {3200} (\bibinfo {year} {2001})}\BibitemShut {NoStop}%
\bibitem [{\citenamefont {Castellano}\ \emph {et~al.}(2009)\citenamefont
  {Castellano}, \citenamefont {Fortunato},\ and\ \citenamefont
  {Loreto}}]{CFL:2009}%
  \BibitemOpen
  \bibfield  {author} {\bibinfo {author} {\bibfnamefont {C.}~\bibnamefont
  {Castellano}}, \bibinfo {author} {\bibfnamefont {S.}~\bibnamefont
  {Fortunato}}, \ and\ \bibinfo {author} {\bibfnamefont {V.}~\bibnamefont
  {Loreto}},\ }\bibfield  {title} {\enquote {\bibinfo {title} {Statistical
  physics of social dynamics},}\ }\href@noop {} {\bibfield  {journal} {\bibinfo
   {journal} {Rev. Mod. Phys.}\ }\textbf {\bibinfo {volume} {81}},\ \bibinfo
  {pages} {591} (\bibinfo {year} {2009})}\BibitemShut {NoStop}%
\bibitem [{\citenamefont {Bashan}\ \emph {et~al.}(2013)\citenamefont {Bashan},
  \citenamefont {Berezin}, \citenamefont {Buldyrev},\ and\ \citenamefont
  {Havlin}}]{BBB:2013}%
  \BibitemOpen
  \bibfield  {author} {\bibinfo {author} {\bibfnamefont {A.}~\bibnamefont
  {Bashan}}, \bibinfo {author} {\bibfnamefont {Y.}~\bibnamefont {Berezin}},
  \bibinfo {author} {\bibfnamefont {S.~V.}\ \bibnamefont {Buldyrev}}, \ and\
  \bibinfo {author} {\bibfnamefont {S.}~\bibnamefont {Havlin}},\ }\bibfield
  {title} {\enquote {\bibinfo {title} {The extreme vulnerability of
  interdependent spatially embedded networks},}\ }\href@noop {} {\bibfield
  {journal} {\bibinfo  {journal} {Nat. Phys.}\ }\textbf {\bibinfo {volume}
  {9}},\ \bibinfo {pages} {667} (\bibinfo {year} {2013})}\BibitemShut {NoStop}%
\bibitem [{\citenamefont {Krapivsky}\ \emph {et~al.}(2010)\citenamefont
  {Krapivsky}, \citenamefont {Redner},\ and\ \citenamefont
  {Ben-Naim}}]{KRB:book}%
  \BibitemOpen
  \bibfield  {author} {\bibinfo {author} {\bibfnamefont {P.~L.}\ \bibnamefont
  {Krapivsky}}, \bibinfo {author} {\bibfnamefont {S.}~\bibnamefont {Redner}}, \
  and\ \bibinfo {author} {\bibfnamefont {E.}~\bibnamefont {Ben-Naim}},\
  }\href@noop {} {\emph {\bibinfo {title} {A Kinetic View of Statistical
  Physics}}}\ (\bibinfo  {publisher} {Cambridge University Press},\ \bibinfo
  {year} {2010})\BibitemShut {NoStop}%
\bibitem [{\citenamefont {Gleeson}(2013)}]{Gleeson:2013}%
  \BibitemOpen
  \bibfield  {author} {\bibinfo {author} {\bibfnamefont {J.~P.}\ \bibnamefont
  {Gleeson}},\ }\bibfield  {title} {\enquote {\bibinfo {title} {Binary-state
  dynamics on complex networks: pair approximation and beyond},}\ }\href@noop
  {} {\bibfield  {journal} {\bibinfo  {journal} {Phys. Rev. X}\ }\textbf
  {\bibinfo {volume} {3}},\ \bibinfo {pages} {021004} (\bibinfo {year}
  {2013})}\BibitemShut {NoStop}%
\bibitem [{\citenamefont {Rudy}\ \emph {et~al.}(2017)\citenamefont {Rudy},
  \citenamefont {Brunton}, \citenamefont {Proctor},\ and\ \citenamefont
  {Kutz}}]{RBPK:2017}%
  \BibitemOpen
  \bibfield  {author} {\bibinfo {author} {\bibfnamefont {S.}~\bibnamefont
  {Rudy}}, \bibinfo {author} {\bibfnamefont {S.~L.}\ \bibnamefont {Brunton}},
  \bibinfo {author} {\bibfnamefont {J.~L.}\ \bibnamefont {Proctor}}, \ and\
  \bibinfo {author} {\bibfnamefont {J.~N.}\ \bibnamefont {Kutz}},\ }\bibfield
  {title} {\enquote {\bibinfo {title} {Data-driven discovery of partial
  differential equations},}\ }\href@noop {} {\bibfield  {journal} {\bibinfo
  {journal} {Sci. Adv.}\ }\textbf {\bibinfo {volume} {3}},\ \bibinfo {pages}
  {e1602614} (\bibinfo {year} {2017})}\BibitemShut {NoStop}%
\bibitem [{\citenamefont {Li}\ \emph {et~al.}(2019)\citenamefont {Li},
  \citenamefont {Li}, \citenamefont {Yue}, \citenamefont {Tang}, \citenamefont
  {Voss}, \citenamefont {Kurths},\ and\ \citenamefont {Y.}}]{LLYTVKY:2019}%
  \BibitemOpen
  \bibfield  {author} {\bibinfo {author} {\bibfnamefont {X.}~\bibnamefont
  {Li}}, \bibinfo {author} {\bibfnamefont {L.}~\bibnamefont {Li}}, \bibinfo
  {author} {\bibfnamefont {Z.}~\bibnamefont {Yue}}, \bibinfo {author}
  {\bibfnamefont {X.}~\bibnamefont {Tang}}, \bibinfo {author} {\bibfnamefont
  {H.~U.}\ \bibnamefont {Voss}}, \bibinfo {author} {\bibfnamefont
  {J.}~\bibnamefont {Kurths}}, \ and\ \bibinfo {author} {\bibfnamefont
  {Y.}~\bibnamefont {Y.}},\ }\bibfield  {title} {\enquote {\bibinfo {title}
  {Sparse learning of partial differential equations with structured dictionary
  matrix},}\ }\href@noop {} {\bibfield  {journal} {\bibinfo  {journal} {Chaos}\
  }\textbf {\bibinfo {volume} {29}},\ \bibinfo {pages} {043130} (\bibinfo
  {year} {2019})}\BibitemShut {NoStop}%
\bibitem [{\citenamefont {Reinbold}\ and\ \citenamefont
  {Grigoriev}(2019)}]{RG:2019}%
  \BibitemOpen
  \bibfield  {author} {\bibinfo {author} {\bibfnamefont {P.~A.~K.}\
  \bibnamefont {Reinbold}}\ and\ \bibinfo {author} {\bibfnamefont {R.~O.}\
  \bibnamefont {Grigoriev}},\ }\bibfield  {title} {\enquote {\bibinfo {title}
  {Data-driven discovery of partial differential equation models with latent
  variables},}\ }\href {\doibase 10.1103/PhysRevE.100.022219} {\bibfield
  {journal} {\bibinfo  {journal} {Phys. Rev. E}\ }\textbf {\bibinfo {volume}
  {100}},\ \bibinfo {pages} {022219} (\bibinfo {year} {2019})}\BibitemShut
  {NoStop}%
\bibitem [{\citenamefont {Gurevich}\ \emph {et~al.}(2019)\citenamefont
  {Gurevich}, \citenamefont {Reinbold},\ and\ \citenamefont
  {Grigoriev}}]{GRG:2019}%
  \BibitemOpen
  \bibfield  {author} {\bibinfo {author} {\bibfnamefont {D.~R.}\ \bibnamefont
  {Gurevich}}, \bibinfo {author} {\bibfnamefont {P.~A.~K.}\ \bibnamefont
  {Reinbold}}, \ and\ \bibinfo {author} {\bibfnamefont {R.~O.}\ \bibnamefont
  {Grigoriev}},\ }\bibfield  {title} {\enquote {\bibinfo {title} {Robust and
  optimal sparse regression for nonlinear pde models},}\ }\href@noop {}
  {\bibfield  {journal} {\bibinfo  {journal} {Chaos}\ }\textbf {\bibinfo
  {volume} {29}},\ \bibinfo {pages} {103113} (\bibinfo {year}
  {2019})}\BibitemShut {NoStop}%
\bibitem [{\citenamefont {Reinbold}\ \emph {et~al.}(2020)\citenamefont
  {Reinbold}, \citenamefont {Gurevich},\ and\ \citenamefont
  {Grigoriev}}]{RGG:2020}%
  \BibitemOpen
  \bibfield  {author} {\bibinfo {author} {\bibfnamefont {P.~A.~K.}\
  \bibnamefont {Reinbold}}, \bibinfo {author} {\bibfnamefont {D.~R.}\
  \bibnamefont {Gurevich}}, \ and\ \bibinfo {author} {\bibfnamefont {R.~O.}\
  \bibnamefont {Grigoriev}},\ }\bibfield  {title} {\enquote {\bibinfo {title}
  {Using noisy or incomplete data to discover models of spatiotemporal
  dynamics},}\ }\href {\doibase 10.1103/PhysRevE.101.010203} {\bibfield
  {journal} {\bibinfo  {journal} {Phys. Rev. E}\ }\textbf {\bibinfo {volume}
  {101}},\ \bibinfo {pages} {010203} (\bibinfo {year} {2020})}\BibitemShut
  {NoStop}%
\bibitem [{\citenamefont {Chirikov}\ and\ \citenamefont
  {Izraelev}(1973)}]{CI:1973}%
  \BibitemOpen
  \bibfield  {author} {\bibinfo {author} {\bibfnamefont {B.~V.}\ \bibnamefont
  {Chirikov}}\ and\ \bibinfo {author} {\bibfnamefont {F.~M.}\ \bibnamefont
  {Izraelev}},\ }\bibfield  {title} {\enquote {\bibinfo {title} {Some numerical
  experiments with a nonlinear mapping: Stochastic component},}\ }\href@noop {}
  {\bibfield  {journal} {\bibinfo  {journal} {Colloques. Int. du CNRS}\
  }\textbf {\bibinfo {volume} {229}},\ \bibinfo {pages} {409} (\bibinfo {year}
  {1973})}\BibitemShut {NoStop}%
\bibitem [{\citenamefont {Chirikov}(1979)}]{Chirikov:1979}%
  \BibitemOpen
  \bibfield  {author} {\bibinfo {author} {\bibfnamefont {B.~V.}\ \bibnamefont
  {Chirikov}},\ }\bibfield  {title} {\enquote {\bibinfo {title} {A universal
  instability of many-dimensional oscillator systems},}\ }\href@noop {}
  {\bibfield  {journal} {\bibinfo  {journal} {Phys. Rep.}\ }\textbf {\bibinfo
  {volume} {52}},\ \bibinfo {pages} {263} (\bibinfo {year} {1979})}\BibitemShut
  {NoStop}%
\bibitem [{\citenamefont {Ikeda}(1979)}]{Ikeda:1979}%
  \BibitemOpen
  \bibfield  {author} {\bibinfo {author} {\bibfnamefont {K.}~\bibnamefont
  {Ikeda}},\ }\bibfield  {title} {\enquote {\bibinfo {title} {Multiple-valued
  stationary state and its instability of the transmitted light by a ring
  cavity system},}\ }\href@noop {} {\bibfield  {journal} {\bibinfo  {journal}
  {Opt. Commun.}\ }\textbf {\bibinfo {volume} {30}},\ \bibinfo {pages} {257}
  (\bibinfo {year} {1979})}\BibitemShut {NoStop}%
\bibitem [{\citenamefont {Hammel}\ \emph {et~al.}(1985)\citenamefont {Hammel},
  \citenamefont {Jones},\ and\ \citenamefont {Moloney}}]{HJM:1985}%
  \BibitemOpen
  \bibfield  {author} {\bibinfo {author} {\bibfnamefont {S.~M.}\ \bibnamefont
  {Hammel}}, \bibinfo {author} {\bibfnamefont {C.~K. R.~T.}\ \bibnamefont
  {Jones}}, \ and\ \bibinfo {author} {\bibfnamefont {J.~V.}\ \bibnamefont
  {Moloney}},\ }\bibfield  {title} {\enquote {\bibinfo {title} {Global
  dynamical behavior of the optical field in a ring cavity},}\ }\href@noop {}
  {\bibfield  {journal} {\bibinfo  {journal} {J. Opt. Soc. Am. B}\ }\textbf
  {\bibinfo {volume} {2}},\ \bibinfo {pages} {552} (\bibinfo {year}
  {1985})}\BibitemShut {NoStop}%
\bibitem [{\citenamefont {Jaeger}(2001)}]{Jaeger:2001}%
  \BibitemOpen
  \bibfield  {author} {\bibinfo {author} {\bibfnamefont {H.}~\bibnamefont
  {Jaeger}},\ }\bibfield  {title} {\enquote {\bibinfo {title} {The “echo
  state” approach to analysing and training recurrent neural networks-with an
  erratum note},}\ }\href@noop {} {\bibfield  {journal} {\bibinfo  {journal}
  {Bonn, Germany: German National Research Center for Information Technology
  GMD Technical Report}\ }\textbf {\bibinfo {volume} {148}},\ \bibinfo {pages}
  {13} (\bibinfo {year} {2001})}\BibitemShut {NoStop}%
\bibitem [{\citenamefont {Mass}\ \emph {et~al.}(2002)\citenamefont {Mass},
  \citenamefont {Nachtschlaeger},\ and\ \citenamefont {Markram}}]{MNM:2002}%
  \BibitemOpen
  \bibfield  {author} {\bibinfo {author} {\bibfnamefont {W.}~\bibnamefont
  {Mass}}, \bibinfo {author} {\bibfnamefont {T.}~\bibnamefont
  {Nachtschlaeger}}, \ and\ \bibinfo {author} {\bibfnamefont {H.}~\bibnamefont
  {Markram}},\ }\bibfield  {title} {\enquote {\bibinfo {title} {Real-time
  computing without stable states: A new framework for neural computation based
  on perturbations},}\ }\href@noop {} {\bibfield  {journal} {\bibinfo
  {journal} {Neur. Comp.}\ }\textbf {\bibinfo {volume} {14}},\ \bibinfo {pages}
  {2531} (\bibinfo {year} {2002})}\BibitemShut {NoStop}%
\bibitem [{\citenamefont {Jaeger}\ and\ \citenamefont {Haas}(2004)}]{JH:2004}%
  \BibitemOpen
  \bibfield  {author} {\bibinfo {author} {\bibfnamefont {H.}~\bibnamefont
  {Jaeger}}\ and\ \bibinfo {author} {\bibfnamefont {H.}~\bibnamefont {Haas}},\
  }\bibfield  {title} {\enquote {\bibinfo {title} {Harnessing nonlinearity:
  Predicting chaotic systems and saving energy in wireless communication},}\
  }\href@noop {} {\bibfield  {journal} {\bibinfo  {journal} {Science}\ }\textbf
  {\bibinfo {volume} {304}},\ \bibinfo {pages} {78} (\bibinfo {year}
  {2004})}\BibitemShut {NoStop}%
\bibitem [{\citenamefont {Manjunath}\ and\ \citenamefont
  {Jaeger}(2013)}]{MJ:2013}%
  \BibitemOpen
  \bibfield  {author} {\bibinfo {author} {\bibfnamefont {G.}~\bibnamefont
  {Manjunath}}\ and\ \bibinfo {author} {\bibfnamefont {H.}~\bibnamefont
  {Jaeger}},\ }\bibfield  {title} {\enquote {\bibinfo {title} {Echo state
  property linked to an input: Exploring a fundamental characteristic of
  recurrent neural networks},}\ }\href@noop {} {\bibfield  {journal} {\bibinfo
  {journal} {Neur. Comp.}\ }\textbf {\bibinfo {volume} {25}},\ \bibinfo {pages}
  {671} (\bibinfo {year} {2013})}\BibitemShut {NoStop}%
\bibitem [{\citenamefont {Haynes}\ \emph {et~al.}(2015)\citenamefont {Haynes},
  \citenamefont {Soriano}, \citenamefont {Rosin}, \citenamefont {Fischer},\
  and\ \citenamefont {Gauthier}}]{HSRFG:2015}%
  \BibitemOpen
  \bibfield  {author} {\bibinfo {author} {\bibfnamefont {N.~D.}\ \bibnamefont
  {Haynes}}, \bibinfo {author} {\bibfnamefont {M.~C.}\ \bibnamefont {Soriano}},
  \bibinfo {author} {\bibfnamefont {D.~P.}\ \bibnamefont {Rosin}}, \bibinfo
  {author} {\bibfnamefont {I.}~\bibnamefont {Fischer}}, \ and\ \bibinfo
  {author} {\bibfnamefont {D.~J.}\ \bibnamefont {Gauthier}},\ }\bibfield
  {title} {\enquote {\bibinfo {title} {Reservoir computing with a single
  time-delay autonomous {Boolean} node},}\ }\href {\doibase
  10.1103/PhysRevE.91.020801} {\bibfield  {journal} {\bibinfo  {journal} {Phys.
  Rev. E}\ }\textbf {\bibinfo {volume} {91}},\ \bibinfo {pages} {020801}
  (\bibinfo {year} {2015})}\BibitemShut {NoStop}%
\bibitem [{\citenamefont {Larger}\ \emph {et~al.}(2017)\citenamefont {Larger},
  \citenamefont {Bayl\'on-Fuentes}, \citenamefont {Martinenghi}, \citenamefont
  {Udaltsov}, \citenamefont {Chembo},\ and\ \citenamefont
  {Jacquot}}]{LBMUCJ:2017}%
  \BibitemOpen
  \bibfield  {author} {\bibinfo {author} {\bibfnamefont {L.}~\bibnamefont
  {Larger}}, \bibinfo {author} {\bibfnamefont {A.}~\bibnamefont
  {Bayl\'on-Fuentes}}, \bibinfo {author} {\bibfnamefont {R.}~\bibnamefont
  {Martinenghi}}, \bibinfo {author} {\bibfnamefont {V.~S.}\ \bibnamefont
  {Udaltsov}}, \bibinfo {author} {\bibfnamefont {Y.~K.}\ \bibnamefont
  {Chembo}}, \ and\ \bibinfo {author} {\bibfnamefont {M.}~\bibnamefont
  {Jacquot}},\ }\bibfield  {title} {\enquote {\bibinfo {title} {High-speed
  photonic reservoir computing using a time-delay-based architecture: Million
  words per second classification},}\ }\href {\doibase
  10.1103/PhysRevX.7.011015} {\bibfield  {journal} {\bibinfo  {journal} {Phys.
  Rev. X}\ }\textbf {\bibinfo {volume} {7}},\ \bibinfo {pages} {011015}
  (\bibinfo {year} {2017})}\BibitemShut {NoStop}%
\bibitem [{\citenamefont {Pathak}\ \emph {et~al.}(2017)\citenamefont {Pathak},
  \citenamefont {Lu}, \citenamefont {Hunt}, \citenamefont {Girvan},\ and\
  \citenamefont {Ott}}]{PLHGO:2017}%
  \BibitemOpen
  \bibfield  {author} {\bibinfo {author} {\bibfnamefont {J.}~\bibnamefont
  {Pathak}}, \bibinfo {author} {\bibfnamefont {Z.}~\bibnamefont {Lu}}, \bibinfo
  {author} {\bibfnamefont {B.}~\bibnamefont {Hunt}}, \bibinfo {author}
  {\bibfnamefont {M.}~\bibnamefont {Girvan}}, \ and\ \bibinfo {author}
  {\bibfnamefont {E.}~\bibnamefont {Ott}},\ }\bibfield  {title} {\enquote
  {\bibinfo {title} {Using machine learning to replicate chaotic attractors and
  calculate {Lyapunov} exponents from data},}\ }\href@noop {} {\bibfield
  {journal} {\bibinfo  {journal} {Chaos}\ }\textbf {\bibinfo {volume} {27}},\
  \bibinfo {pages} {121102} (\bibinfo {year} {2017})}\BibitemShut {NoStop}%
\bibitem [{\citenamefont {Lu}\ \emph {et~al.}(2017)\citenamefont {Lu},
  \citenamefont {Pathak}, \citenamefont {Hunt}, \citenamefont {Girvan},
  \citenamefont {Brockett},\ and\ \citenamefont {Ott}}]{LPHGBO:2017}%
  \BibitemOpen
  \bibfield  {author} {\bibinfo {author} {\bibfnamefont {Z.}~\bibnamefont
  {Lu}}, \bibinfo {author} {\bibfnamefont {J.}~\bibnamefont {Pathak}}, \bibinfo
  {author} {\bibfnamefont {B.}~\bibnamefont {Hunt}}, \bibinfo {author}
  {\bibfnamefont {M.}~\bibnamefont {Girvan}}, \bibinfo {author} {\bibfnamefont
  {R.}~\bibnamefont {Brockett}}, \ and\ \bibinfo {author} {\bibfnamefont
  {E.}~\bibnamefont {Ott}},\ }\bibfield  {title} {\enquote {\bibinfo {title}
  {Reservoir observers: Model-free inference of unmeasured variables in chaotic
  systems},}\ }\href@noop {} {\bibfield  {journal} {\bibinfo  {journal}
  {Chaos}\ }\textbf {\bibinfo {volume} {27}},\ \bibinfo {pages} {041102}
  (\bibinfo {year} {2017})}\BibitemShut {NoStop}%
\bibitem [{\citenamefont {Lu}\ \emph {et~al.}(2018)\citenamefont {Lu},
  \citenamefont {Hunt},\ and\ \citenamefont {Ott}}]{LHO:2018}%
  \BibitemOpen
  \bibfield  {author} {\bibinfo {author} {\bibfnamefont {Z.}~\bibnamefont
  {Lu}}, \bibinfo {author} {\bibfnamefont {B.~R.}\ \bibnamefont {Hunt}}, \ and\
  \bibinfo {author} {\bibfnamefont {E.}~\bibnamefont {Ott}},\ }\bibfield
  {title} {\enquote {\bibinfo {title} {Attractor reconstruction by machine
  learning},}\ }\href@noop {} {\bibfield  {journal} {\bibinfo  {journal}
  {Chaos}\ }\textbf {\bibinfo {volume} {28}},\ \bibinfo {pages} {061104}
  (\bibinfo {year} {2018})}\BibitemShut {NoStop}%
\bibitem [{\citenamefont {Pathak}\ \emph
  {et~al.}(2018{\natexlab{a}})\citenamefont {Pathak}, \citenamefont {Wilner},
  \citenamefont {Fussell}, \citenamefont {Chandra}, \citenamefont {Hunt},
  \citenamefont {Girvan}, \citenamefont {Lu},\ and\ \citenamefont
  {Ott}}]{PWFCHGO:2018}%
  \BibitemOpen
  \bibfield  {author} {\bibinfo {author} {\bibfnamefont {J.}~\bibnamefont
  {Pathak}}, \bibinfo {author} {\bibfnamefont {A.}~\bibnamefont {Wilner}},
  \bibinfo {author} {\bibfnamefont {R.}~\bibnamefont {Fussell}}, \bibinfo
  {author} {\bibfnamefont {S.}~\bibnamefont {Chandra}}, \bibinfo {author}
  {\bibfnamefont {B.}~\bibnamefont {Hunt}}, \bibinfo {author} {\bibfnamefont
  {M.}~\bibnamefont {Girvan}}, \bibinfo {author} {\bibfnamefont
  {Z.}~\bibnamefont {Lu}}, \ and\ \bibinfo {author} {\bibfnamefont
  {E.}~\bibnamefont {Ott}},\ }\bibfield  {title} {\enquote {\bibinfo {title}
  {Hybrid forecasting of chaotic processes: Using machine learning in
  conjunction with a knowledge-based model},}\ }\href@noop {} {\bibfield
  {journal} {\bibinfo  {journal} {Chaos}\ }\textbf {\bibinfo {volume} {28}},\
  \bibinfo {pages} {041101} (\bibinfo {year} {2018}{\natexlab{a}})}\BibitemShut
  {NoStop}%
\bibitem [{\citenamefont {Pathak}\ \emph
  {et~al.}(2018{\natexlab{b}})\citenamefont {Pathak}, \citenamefont {Hunt},
  \citenamefont {Girvan}, \citenamefont {Lu},\ and\ \citenamefont
  {Ott}}]{PHGLO:2018}%
  \BibitemOpen
  \bibfield  {author} {\bibinfo {author} {\bibfnamefont {J.}~\bibnamefont
  {Pathak}}, \bibinfo {author} {\bibfnamefont {B.}~\bibnamefont {Hunt}},
  \bibinfo {author} {\bibfnamefont {M.}~\bibnamefont {Girvan}}, \bibinfo
  {author} {\bibfnamefont {Z.}~\bibnamefont {Lu}}, \ and\ \bibinfo {author}
  {\bibfnamefont {E.}~\bibnamefont {Ott}},\ }\bibfield  {title} {\enquote
  {\bibinfo {title} {Model-free prediction of large spatiotemporally chaotic
  systems from data: A reservoir computing approach},}\ }\href {\doibase
  10.1103/PhysRevLett.120.024102} {\bibfield  {journal} {\bibinfo  {journal}
  {Phys. Rev. Lett.}\ }\textbf {\bibinfo {volume} {120}},\ \bibinfo {pages}
  {024102} (\bibinfo {year} {2018}{\natexlab{b}})}\BibitemShut {NoStop}%
\bibitem [{\citenamefont {Carroll}(2018)}]{Carroll:2018}%
  \BibitemOpen
  \bibfield  {author} {\bibinfo {author} {\bibfnamefont {T.~L.}\ \bibnamefont
  {Carroll}},\ }\bibfield  {title} {\enquote {\bibinfo {title} {Using reservoir
  computers to distinguish chaotic signals},}\ }\href {\doibase
  10.1103/PhysRevE.98.052209} {\bibfield  {journal} {\bibinfo  {journal} {Phys.
  Rev. E}\ }\textbf {\bibinfo {volume} {98}},\ \bibinfo {pages} {052209}
  (\bibinfo {year} {2018})}\BibitemShut {NoStop}%
\bibitem [{\citenamefont {Nakai}\ and\ \citenamefont {Saiki}(2018)}]{NS:2018}%
  \BibitemOpen
  \bibfield  {author} {\bibinfo {author} {\bibfnamefont {K.}~\bibnamefont
  {Nakai}}\ and\ \bibinfo {author} {\bibfnamefont {Y.}~\bibnamefont {Saiki}},\
  }\bibfield  {title} {\enquote {\bibinfo {title} {Machine-learning inference
  of fluid variables from data using reservoir computing},}\ }\href {\doibase
  10.1103/PhysRevE.98.023111} {\bibfield  {journal} {\bibinfo  {journal} {Phys.
  Rev. E}\ }\textbf {\bibinfo {volume} {98}},\ \bibinfo {pages} {023111}
  (\bibinfo {year} {2018})}\BibitemShut {NoStop}%
\bibitem [{\citenamefont {Roland}\ and\ \citenamefont
  {Parlitz}(2018)}]{ZP:2018}%
  \BibitemOpen
  \bibfield  {author} {\bibinfo {author} {\bibfnamefont {Z.~S.}\ \bibnamefont
  {Roland}}\ and\ \bibinfo {author} {\bibfnamefont {U.}~\bibnamefont
  {Parlitz}},\ }\bibfield  {title} {\enquote {\bibinfo {title} {Observing
  spatio-temporal dynamics of excitable media using reservoir computing},}\
  }\href@noop {} {\bibfield  {journal} {\bibinfo  {journal} {Chaos}\ }\textbf
  {\bibinfo {volume} {28}},\ \bibinfo {pages} {043118} (\bibinfo {year}
  {2018})}\BibitemShut {NoStop}%
\bibitem [{\citenamefont {Weng}\ \emph {et~al.}(2019)\citenamefont {Weng},
  \citenamefont {Yang}, \citenamefont {Gu}, \citenamefont {Zhang},\ and\
  \citenamefont {Small}}]{WYGZS:2019}%
  \BibitemOpen
  \bibfield  {author} {\bibinfo {author} {\bibfnamefont {T.}~\bibnamefont
  {Weng}}, \bibinfo {author} {\bibfnamefont {H.}~\bibnamefont {Yang}}, \bibinfo
  {author} {\bibfnamefont {C.}~\bibnamefont {Gu}}, \bibinfo {author}
  {\bibfnamefont {J.}~\bibnamefont {Zhang}}, \ and\ \bibinfo {author}
  {\bibfnamefont {M.}~\bibnamefont {Small}},\ }\bibfield  {title} {\enquote
  {\bibinfo {title} {Synchronization of chaotic systems and their
  machine-learning models},}\ }\href {\doibase 10.1103/PhysRevE.99.042203}
  {\bibfield  {journal} {\bibinfo  {journal} {Phys. Rev. E}\ }\textbf {\bibinfo
  {volume} {99}},\ \bibinfo {pages} {042203} (\bibinfo {year}
  {2019})}\BibitemShut {NoStop}%
\bibitem [{\citenamefont {Griffith}\ \emph {et~al.}(2019)\citenamefont
  {Griffith}, \citenamefont {Pomerance},\ and\ \citenamefont
  {Gauthier}}]{GPG:2019}%
  \BibitemOpen
  \bibfield  {author} {\bibinfo {author} {\bibfnamefont {A.}~\bibnamefont
  {Griffith}}, \bibinfo {author} {\bibfnamefont {A.}~\bibnamefont {Pomerance}},
  \ and\ \bibinfo {author} {\bibfnamefont {D.~J.}\ \bibnamefont {Gauthier}},\
  }\bibfield  {title} {\enquote {\bibinfo {title} {Forecasting chaotic systems
  with very low connectivity reservoir computers},}\ }\href@noop {} {\bibfield
  {journal} {\bibinfo  {journal} {Chaos}\ }\textbf {\bibinfo {volume} {29}},\
  \bibinfo {pages} {123108} (\bibinfo {year} {2019})}\BibitemShut {NoStop}%
\bibitem [{\citenamefont {Jiang}\ and\ \citenamefont {Lai}(2019)}]{JL:2019}%
  \BibitemOpen
  \bibfield  {author} {\bibinfo {author} {\bibfnamefont {J.}~\bibnamefont
  {Jiang}}\ and\ \bibinfo {author} {\bibfnamefont {Y.-C.}\ \bibnamefont
  {Lai}},\ }\bibfield  {title} {\enquote {\bibinfo {title} {Irrelevance of
  linear controllability to nonlinear dynamical networks},}\ }\href@noop {}
  {\bibfield  {journal} {\bibinfo  {journal} {Nat. Commun.}\ }\textbf {\bibinfo
  {volume} {10}},\ \bibinfo {pages} {3961} (\bibinfo {year}
  {2019})}\BibitemShut {NoStop}%
\bibitem [{\citenamefont {Vlachas}\ \emph {et~al.}(2019)\citenamefont
  {Vlachas}, \citenamefont {Pathak}, \citenamefont {Hunt}, \citenamefont
  {Sapsis}, \citenamefont {Girvan}, \citenamefont {Ott},\ and\ \citenamefont
  {Koumoutsakos}}]{VPHSGOK:2019}%
  \BibitemOpen
  \bibfield  {author} {\bibinfo {author} {\bibfnamefont {P.~R.}\ \bibnamefont
  {Vlachas}}, \bibinfo {author} {\bibfnamefont {J.}~\bibnamefont {Pathak}},
  \bibinfo {author} {\bibfnamefont {B.~R.}\ \bibnamefont {Hunt}}, \bibinfo
  {author} {\bibfnamefont {T.~P.}\ \bibnamefont {Sapsis}}, \bibinfo {author}
  {\bibfnamefont {M.}~\bibnamefont {Girvan}}, \bibinfo {author} {\bibfnamefont
  {E.}~\bibnamefont {Ott}}, \ and\ \bibinfo {author} {\bibfnamefont
  {P.}~\bibnamefont {Koumoutsakos}},\ }\bibfield  {title} {\enquote {\bibinfo
  {title} {Forecasting of spatio-temporal chaotic dynamics with recurrent
  neural networks: A comparative study of reservoir computing and
  backpropagation algorithms},}\ }\href@noop {} {\bibfield  {journal} {\bibinfo
   {journal} {arXiv preprint arXiv:1910.05266}\ } (\bibinfo {year}
  {2019})}\BibitemShut {NoStop}%
\bibitem [{\citenamefont {Fan}\ \emph {et~al.}(2020)\citenamefont {Fan},
  \citenamefont {Jiang}, \citenamefont {Zhang}, \citenamefont {Wang},\ and\
  \citenamefont {Lai}}]{FJZWL:2020}%
  \BibitemOpen
  \bibfield  {author} {\bibinfo {author} {\bibfnamefont {H.}~\bibnamefont
  {Fan}}, \bibinfo {author} {\bibfnamefont {J.}~\bibnamefont {Jiang}}, \bibinfo
  {author} {\bibfnamefont {C.}~\bibnamefont {Zhang}}, \bibinfo {author}
  {\bibfnamefont {X.}~\bibnamefont {Wang}}, \ and\ \bibinfo {author}
  {\bibfnamefont {Y.-C.}\ \bibnamefont {Lai}},\ }\bibfield  {title} {\enquote
  {\bibinfo {title} {Long-term prediction of chaotic systems with machine
  learning},}\ }\href {\doibase 10.1103/PhysRevResearch.2.012080} {\bibfield
  {journal} {\bibinfo  {journal} {Phys. Rev. Research}\ }\textbf {\bibinfo
  {volume} {2}},\ \bibinfo {pages} {012080} (\bibinfo {year}
  {2020})}\BibitemShut {NoStop}%
\bibitem [{\citenamefont {Zhang}\ \emph {et~al.}(2020)\citenamefont {Zhang},
  \citenamefont {Jiang}, \citenamefont {Qu},\ and\ \citenamefont
  {Lai}}]{ZJQL:2020}%
  \BibitemOpen
  \bibfield  {author} {\bibinfo {author} {\bibfnamefont {C.}~\bibnamefont
  {Zhang}}, \bibinfo {author} {\bibfnamefont {J.}~\bibnamefont {Jiang}},
  \bibinfo {author} {\bibfnamefont {S.-X.}\ \bibnamefont {Qu}}, \ and\ \bibinfo
  {author} {\bibfnamefont {Y.-C.}\ \bibnamefont {Lai}},\ }\bibfield  {title}
  {\enquote {\bibinfo {title} {Predicting phase and sensing phase coherence in
  chaotic systems with machine learning},}\ }\href@noop {} {\bibfield
  {journal} {\bibinfo  {journal} {Chaos}\ }\textbf {\bibinfo {volume} {30}},\
  \bibinfo {pages} {083114} (\bibinfo {year} {2020})}\BibitemShut {NoStop}%
\end{thebibliography}

%
\end{document}